\documentclass{article}

\usepackage{arxiv}
\usepackage{appendix}

\pdfoutput=1 

\usepackage[utf8]{inputenc} 
\usepackage[T1]{fontenc}    
\usepackage{hyperref}       
\usepackage{url}            
\usepackage{booktabs}       
\usepackage{amsfonts}       
\usepackage{nicefrac}       
\usepackage{microtype}      
\usepackage{doi}

\usepackage{amssymb}
\usepackage{amsmath}
\usepackage{mathtools} 
\usepackage{subfig}
\usepackage{xcolor}
\usepackage{media9} 
\usepackage{graphicx, animate}
\usepackage{units} 
\usepackage{esvect}
\usepackage{comment}
\usepackage{cite} 
\usepackage{float}
\usepackage{bm}

\usepackage[labelfont={bf,sf},font={small}, labelsep=space]{caption}
\usepackage{algorithm}
\usepackage{algpseudocode}

\usepackage{listings}

\usepackage{todonotes}

\newcommand{\R}{\mathbb{R}}
\newcommand{\divergence}{\mathop\mathrm{div}}

\usepackage[markup=underlined
,addedmarkup=colored
,commandnameprefix=ifneeded]{changes}

\makeatletter
\@namedef{Changes@AuthorColor}{magenta}
\colorlet{Changes@Color}{magenta}
\makeatother

\title{Hybrid PDE-ODE Models for Efficient Simulation of Infection Spread in Epidemiology}

\date{} 					

\author{ 
{
Kristina Kehrer}
\\
	Zuse Institute Berlin\\
	Berlin, 14195 \\
	\texttt{kehrer@zib.de} \\
	\And
 {
  Martin Weiser} \\
	Zuse Institute Berlin\\
	Berlin, 14195 \\
	\texttt{weiser@zib.de} \\
\And
 {
Tim Conrad} \\
	Zuse Institute Berlin\\
	Berlin, 14195 \\
	\texttt{conrad@zib.de} \\
}

\renewcommand{\headeright}{Paper}

\begin{document}
\maketitle
\begin{abstract} 
This paper introduces a novel hybrid model combining Partial Differential Equations (PDEs) and Ordinary Differential Equations (ODEs) to simulate infectious disease dynamics across geographic regions. By leveraging the spatial detail of PDEs and the computational efficiency of ODEs, the model enables rapid evaluation of public health interventions. Applied to synthetic environments and real-world scenarios in Lombardy, Italy, and Berlin, Germany, the model highlights how interactions between PDE and ODE regions affect infection dynamics, especially in high-density areas. 
Key findings reveal that the placement of model boundaries in densely populated regions can lead to inaccuracies in infection spread, suggesting that boundaries should be positioned in areas of lower population density to better reflect transmission dynamics. Additionally, regions with low population density hinder infection flow, indicating a need for incorporating, e.g., jumps in the model to enhance its predictive capabilities.
Results indicate that the hybrid model achieves a balance between computational speed and accuracy, making it a valuable tool for policymakers in real-time decision-making and scenario analysis in epidemiology and potentially in other fields requiring similar modeling approaches.
\end{abstract}


\section{Introduction}

In the rapidly evolving field of mathematical modeling, the ability to accurately and efficiently simulate large complex systems is of great importance. This is particularly true in the field of epidemiology, where understanding the dynamics of disease spread can have significant implications for public health and policy decisions. In this field, multiple modeling approaches exist, such as agent-based models (ABMs) \cite{Muller2021ABM}, models based on partial differential equations (PDEs) \cite{Foutel2022PDE}, and ordinary differential equation (ODE) models \cite{Wulkow2021}.

While ABMs are well-suited for simulating detailed actions and interactions at the individual level, they can be computationally very demanding, especially when simulating large populations over extended periods or aiming at policy optimization~\cite{NiemannUramWolfDjurdjevacconradWeiser2024}. This reduces their suitability for real-time analysis or situations in which fast computational processing is essential. 

PDE models, on the other hand, complement the high-resolution insights provided by ABMs. PDE models are generally less computationally demanding compared to ABMs and thus allow simulations of larger geographical areas, in a spatially continuous environment. They provide valuable insights into the general transmission patterns by considering factors like population density and mobility, which enable modelers to gain valuable information about the modeled systems on a regional or national level, such as possible areas of high infection rates. However, solving large PDE systems also has some considerable computational complexity -- which can become a significant limitation, particularly when dealing with large environments and extensive scenarios on a large scale.

ODE models, in contrast, emphasize the temporal dynamics of disease transmission, offering computational efficiency but lacking the spatial resolution of PDE models and the detailed granularity of ABMs. This makes them extremely valuable in situations when rapid analysis is needed. Nevertheless, this benefit is accompanied by the drawback of decreased accuracy in representing spatially diverse situations.

While each of these modeling approaches — ABMs, PDEs, and ODEs — plays a crucial role in understanding epidemiological dynamics, integrating these models to leverage their distinct advantages remains challenging. In this paper, we introduce a hybrid PDE-ODE model that seeks to merge these approaches efficiently, providing a new perspective for enhanced epidemiological modeling.

The spread of infections has been explored using both machine learning \cite{Gaskin2024Neural,Rodríguez2023EINNs,Kharazmi2021Identifiability} and finite difference solvers \cite{ARIF2024100660,cmes.2023.028946,math11234794}. In this study, we opted for finite element methods due to their effectiveness in accurately resolving spatially complex regions.
Coupled PDE-ODE models have been used to model various phenomena, such as tumor development\cite{Romero2021}, traffic control\cite{Liard2023} and others, alongside specialized solvers \cite{Farrell2019}. However, not many applications are known in the context of disease dynamics. In our context, the innovative aspect of this coupling is the incorporation of the ODE model as a space-invariant component within the PDE framework. This integration allows for the exchange of boundary information between the two models at each time step, guaranteeing a dynamic yet consistent interplay between local and global disease dynamics. The resulting hybrid PDE-ODE model offers a faster alternative to a full-PDE model while maintaining an acceptably low loss of accuracy.

In summary, this paper introduces a new approach for coupling partial differential equation (PDE) and ordinary differential equation (ODE) models, in this case tailored for epidemiological applications but with potential implications beyond this field. By achieving a balance between the detailed spatial representation of PDEs and the computational efficiency of ODEs, this approach represents a significant advancement in the field of mathematical modeling. In the field of epidemiology, this provides a powerful tool for policymakers to rapidly and accurately assess the impact of various interventions and scenarios in the management of infectious diseases.~\\
~\\
The source code is available at this URL: \url{https://zenodo.org/records/11121452}.

\section{Full-PDE Model} 


In this section, we provide a comprehensive formulation of the hybrid PDE-ODE model, designed for an epidemiological use-case, in particular for the study of infection spreading. We begin by detailing the high-fidelity PDE model (Sec.~\ref{sec:model-equations}) that captures the spatial dynamics of the disease, including the definition of initial values (Sec.~\ref{sec:initial-values}) to accurately represent the starting conditions of the disease spread. We then proceed to discuss the approximation of regions of lesser interest using ODEs (Sec.~\ref{sec:model-reduction}), and finally, we present the methods for seamlessly coupling these models. This hybrid approach not only retains the spatial resolution critical for accurate simulation but also leverages the computational simplicity of ODEs, resulting in an efficient yet powerful tool for analyzing epidemiological data. This structured approach ensures that the model can simulate the progression of infectious diseases with both spatial detail and computational efficiency.

\subsection{Model Equations} \label{sec:model-equations}
We consider a spatial spreading model of SEIR type in a domain $\Omega\subset\R^2$, inspired by the following papers \cite{Quelle5,viguerie2020diffusion}. We employ a reaction-diffusion equation modeling both the SEIR dynamics and spatial diffusion of infectious individuals. The choice of the SEIR model is motivated by its ability to capture the essential stages of disease progression (susceptible, exposed, infectious, and removed), providing a detailed and realistic framework for understanding the dynamics of disease transmission within a population.

Let $s,e,i,r:\Omega\times[0,T]\to\R$ denote the fraction of population density of susceptible, exposed, infectious, and removed (recovered and immune or deceased) subjects, such that $s+e+i+r=1$. Their evolution is given by
\begin{equation} \label{PDE-model}
\begin{aligned}
    \partial_t s &= \divergence (D \nabla s) - \left(1-\frac{A}{n+n_0}\right) s (\beta_e e + \beta_i i) && \text{in $\Omega$,} \\
    \partial_t e &= \divergence(D \nabla e) + \left(1-\frac{A}{n+n_0}\right)s(\beta_e e + \beta_i i) - \sigma e - \phi_e e && \text{in $\Omega$,} \\
    \partial_t i &= \divergence(D\nabla i) + \sigma e - \phi_i i && \text{in $\Omega$,} \\
    \partial_t r &= \divergence (D \nabla r) + \phi_i i + \phi_e e && \text{in $\Omega$,} \\
    \nu^T D \nabla s &=  \nu^T D \nabla e = \nu^T D \nabla i = \nu^T D \nabla r = 0  && \text{on $\partial\Omega$}
\end{aligned}
\end{equation}
with $0<n \in L^\infty(\Omega)$ denoting the given population density and $\nu$ being the unit outer normal to $\Omega$. These equations (\ref{PDE-model}) represent the fractional densities of the population in each group, postulated directly without first formulating equations for the total population sizes.
The diffusion coefficient $D>0$ describes the spatial mixing of the subpopulations and may, in general, depend on the spatial position. The constants $\beta_e > 0$ and  $\beta_i > 0$ describe the infection rates, while the recovery rates are represented by $\phi_e$ and $\phi_i$. The rate of becoming infectious is denoted as $\sigma$.
A description of the parameters can be found in Table \ref{table: parameters}.
The homogeneous Neumann boundary conditions ensure that the fractions add up to one for all times.

The Allee effect plays a crucial role in (cell) population biology, explaining the positive correlation between per capita growth rate and (cell) population density \cite{Johnson2019Cancer,KumarSinha+2023Dynamics}. It can be expressed as $(1-A/n)$, where $A$ represents the Allee threshold. If the population density $n$ falls below this threshold $A$, the net growth rate becomes negative, indicating potential extinction. As the population density approaches the threshold, growth slows down but remains positive. Once the population density surpasses the threshold, the Allee effect diminishes, making it practically irrelevant when $A \ll n$.
In \cite{Quelle5,viguerie2020diffusion}, the Allee effect was instead used to model the infection rate as a function of population density, where it represents a lower transmission probability in less densely populated regions. 
This can simply be referred to as a nonlinear infection rate \cite{TURNES2014epidemic} and is similar to incidence functions \cite{Quelle3}. 
The influence of population density on the infection rate is controversial \cite{Bhadra2020Impact,Hamidi2020Density,Wong2020Spreading}.  
Nevertheless, we adopt the effect and bound it from below by $\nicefrac{1}{3}$. This ensures that the effect is at most three times lower in sparsely populated areas than in regions with high population density. To enforce the lower bound, we apply a shift $n_0$ in the Allee term.
We choose $n_0$ such that the 
factor is bounded from below 
\begin{align*}
    1-\frac{A}{n+n_0} \geq \frac{1}{3} \qquad
\Leftrightarrow \qquad  n+n_0 \geq \frac{3}{2} A 
\end{align*}
for all $n\in \mathbb{R}^+$. 
Hence, we choose $n_0 \coloneqq \nicefrac{3}{2} A$.
\begin{table}[H] 
\centering %
\begin{tabular*}{0.53\textwidth}{c  l} 
\toprule
 parameter & definition  \\
  \midrule %
  $s(x,t)$ & fraction of population density of susceptibles \\ 
  $e(x,t)$ & fraction of population density of exposed \\ 
  $i(x,t)$ & fraction of population density of infectious \\ 
  $r(x,t)$ & fraction of population density of removed \\ 
  $n(x,t)$ & population density \\ 
  $\beta_i$  & symptomatic infection rate \\ 
  $\beta_e$  & asymptomatic infection rate  \\ 
  $\sigma$  & inverse of incubation period \\
  $\phi_e$  & asymptomatic recovery rate \\
  $\phi_i$  & infected recovery rate \\
  $D$  & diffusion coefficient \\ 
\bottomrule 
\end{tabular*}  %
\caption{Parameters of the model (\ref{PDE-model}). 
} 
\label{table: parameters}%
\end{table}%
\noindent
%
%
\subsection{Initial Conditions} \label{sec:initial-values} 
Next, we address the challenge of defining initial starting values and estimating distribution patterns for our model. Due to the typical unavailability of precise spatial densities for exposed, infected, and removed populations, we use aggregated data, which is usually available at broader geographical levels such as districts, subregions, or provinces. The following outlines our approach for setting initial conditions for fractional densities within the PDE system. This step is important to ensure that our simulations start with an accurate and realistic representation of the modeled disease's spatial distribution -- based on the available data.

Let $\Omega = \bigcup_{l=1}^L \Omega^{(l)}$ be a decomposition of $\Omega$ into provinces $\Omega^{(l)}$ for which aggregated numbers of infectious and of removed individuals are given.

The initial conditions are computed using the real (smoothed) 
number $\hat{\mathcal{I}}^{(l)}$  of infectious and  $\hat{\mathcal{R}}^{(l)}$ of removed individuals for each subdomain $\Omega^{(l)}$, $l=1,\dots,L$.
Our given data for each subdomain $\Omega^{(l)}$
is
\begin{align*}
    &\hat{\mathcal{N}}^{(l)} &&\text{(total population number)}, \\
    &\hat{\mathcal{I}}^{(l)} &&\text{(number of infectious individuals)}, \\
    &\hat{\mathcal{R}}^{(l)} &&\text{(number of removed individuals)}. \\
\end{align*}
Further, we set
\begin{align*}
    &\hat{\mathcal{E}}^{(l)} = \Sigma_{\mathcal{E}} \hat{\mathcal{I}}^{(l)} &&\text{(number of exposed)}, \\
    &\hat{\mathcal{S}}^{(l)} = \hat{\mathcal{N}}_l - (1+ \Sigma_{\mathcal{E}}) \hat{\mathcal{I}}^{(l)} -  \hat{\mathcal{R}}^{(l)} &&\text{(number of susceptibles)},
\end{align*}
where $\Sigma_{\mathcal{E}}$ scales the number of 
infectious to estimate the number of exposed individuals. 
We use isotropic Gaussian functions $G^{(l)}(x) = \tfrac{1}{2\pi \sigma^2}\exp\Bigl(-\tfrac{\|x-x_0^{(l)}\|^2}{2\sigma^2}\Bigr)$ located at the centers $x_0^{(l)}$ of the provinces $\Omega^{(l)}\subset\Omega$ for representing density and density fractions as
\begin{equation*}
n(x) = \sum_{\tilde{l}=1}^{\tilde{L}} w_{n}^{(\tilde{l})} G^{(\tilde{l})}(x),
\quad s(x,0) = \sum_{\tilde{l}=1}^{\tilde{L}} w_{s}^{(\tilde{l})} G^{(\tilde{l})}(x), 
\quad i(x,0) = \sum_{\tilde{l}=1}^{\tilde{L}} w_{i}^{(\tilde{l})} G^{(\tilde{l})}(x), 
\quad\text{and}
\quad e(x,0) = \sum_{\tilde{l}=1}^{\tilde{L}} w_{e}^{(\tilde{l})} G^{(\tilde{l})}(x). 
\end{equation*}
By incorporating spatial variation from the start, we ensure that the Allee effect influences the population dynamics, allowing us to capture its impact on the infection spread more effectively throughout the domain.
The coefficients $w^{(\tilde{l})}$ are determined such that the integral of the densities over each subdomain area $\Omega^{(l)}$ is equal to the given target data 
$\hat{\mathcal{N}}_l,\hat{\mathcal{S}}^{(l)},\hat{\mathcal{I}}^{(l)},\hat{\mathcal{E}}^{(l)}$ 
at time $t=0$, i.e., 
\begin{align*}
&\int_{\Omega^{(l)}} n(x) \; dx = \hat{\mathcal{N}}_l, 
&&\int_{\Omega^{(l)}} s(x,0)n(x) \; dx = \hat{\mathcal{N}}_l - (1+ \Sigma_{\mathcal{E}}) \hat{\mathcal{I}}^{(l)} - \hat{\mathcal{R}}^{(l)}, \\
    &\int_{\Omega^{(l)}} i(x,0)n(x) \; dx = \hat{\mathcal{I}}^{(l)}, 
    &&\int_{\Omega^{(l)}} e(x,0)n(x) \; dx = \Sigma_{\mathcal{E}} \hat{\mathcal{I}}^{(l)}, \\
\end{align*}
leading to linear equation systems of dimension $\tilde{L}$. As $n$ is required for computing the initial densities $s,i$, and $e$, it is computed first.

\subsection{Hybrid PDE-ODE Model}\label{sec:model-reduction}
In the following section, we describe how we define our hybrid model based on the PDE model by segmenting the modeling domain into a PDE part and a -- less complex -- ODE part. Our approach is inspired by previous work where similar concepts have been used for models simulating microbial sensing \cite{Quelle2} and airborne pathogen movement in an SIR framework \cite{Quelle1}. Our approach adopts these concepts to simplify computational demands without compromising the essential dynamics of disease spread. We utilize integro-differential equations for the coupling, enabling us to maintain high-resolution modeling where it is crucial while simplifying other areas to enhance computational efficiency. This will ultimately allow us to efficiently model large-scale epidemiological phenomena by reducing complexity in areas where high resolution is less critical. We detail the mathematical underpinnings of this reduction, including the seamless integration and boundary conditions that ensure continuity and consistency across the coupled models.

We first consider a non-overlapping domain decomposition \cite{ALVES2021Domain,RICE1997Analysis}.
We assume the domain $\Omega$ comprises a region $\Omega_1\subset\R^2$ of interest to be modeled in spatial detail, and an adjacent region $\Omega_2\subset\R^2$ of less interest, where a coarser model is sufficient. In the region of interest, we employ a reaction-diffusion equation modeling both the SEIR dynamics and spatial diffusion of infectious individuals. On the domain $\Omega_2$, we assume homogeneous mixing and formulate a straightforward ODE model. While homogeneous mixing is an unreasonable assumption for $\Omega_2$, we adopt it to analyze the extent of error it introduces, particularly its impact on the region of interest, $\Omega_1$. The models are coupled across their joint boundary $\Gamma :=\partial\Omega_1 \cap \partial\Omega_2$.
This model coupling aims to achieve significantly reduced computational costs compared to implementing the PDE model on the entire domain $\Omega \coloneqq \mathop \mathrm{int}\overline{\Omega_1 \cup\Omega_2}$.

Let us first consider the structure of the full-PDE solution on $\Omega$. Across the interface $\Gamma$, separating $\Omega_1$ from $\Omega_2$, continuity of densities and continuity of fluxes hold almost everywhere. Introducing $s_1 = s|_{\Omega_1}$ and $s_2 = s|_{\Omega_2}$ this reads
\begin{align}
 s_1 &= s_2 && \text{a.e. on $\Gamma$}, \label{eq:continuity-density}\\
 \nu_1^T D \nabla s_1 &= \nu_1^T D \nabla s_2 && \text{a.e. on $\Gamma$}
\end{align}
with $\nu_1$ being the outer unit normal of $\Omega_1$. The same holds for $e$, $i$, and $r$.


Extending $s_1$ and $s_2$ by zero to the whole domain $\Omega$,  we can write the relative density of susceptibles in $\Omega$ as $s=s_1 + s_2$. Multiplying the first PDE of model (\ref{PDE-model}) by a test function $\varphi$ and integrating over $\Omega_2$ results in
\begin{align*}
    \int_{\Omega_2} \dot s \varphi \; dx 
    &= \int_{\Omega_2} \bigl(\divergence (D \nabla s) - \left(1-\frac{A}{n+n_0}\right) s (\beta_e e + \beta_i i) \bigr) \varphi \; dx \\
\Leftrightarrow    \int_{\Omega_2} \dot s_2 \varphi \; dx 
    &= \int_{\Omega_2} \bigl(\divergence (D \nabla s_2) - \left(1-\frac{A}{n_2+n_0}\right) s_2 (\beta_e e_2 + \beta_i i_2) \bigr) \varphi \; dx \\
    &= \int_{\Omega_2} - \nabla \varphi^T D \nabla s_2 - \left(1-\frac{A}{n_2+n_0}\right) s_2 (\beta_e e_2 + \beta_i i_2) \varphi \; dx + \int_{\Gamma} \varphi  \nu_2^T D \nabla s_2 \; do \\
    &= \int_{\Omega_2} - \nabla \varphi^T D \nabla s_2 - \left(1-\frac{A}{n_2+n_0}\right) s_2 (\beta_e e_2 + \beta_i i_2) \varphi \; dx - \int_{\Gamma} \varphi  \nu_1^T D \nabla s_1 \; do.
\end{align*}
Further, let us decompose $s_2 = \bar{s}_2 + \tilde{s}_2$ into a spatially constant part $\bar{s}_2 = \frac{1}{|\Omega_2|} \int_{\Omega_2} s_2 \; dx$ and a fluctuating part $\tilde s_2$ with mean zero. 

Testing with $\varphi = \frac{1}{|\Omega_2|}$ 
yields
\begin{align*}
    && \frac{1}{|\Omega_2|} \int_{\Omega_2} \dot s_2 \; dx 
    &= \frac{1}{|\Omega_2|} \int_{\Omega_2} - \left(1-\frac{A}{\bar{n}_2+\tilde{n}_2+n_0}\right) s_2 (\beta_e e_2 + \beta_i i_2) \; dx - \frac{1}{|\Omega_2|} \int_{\Gamma}\nu_1^T D \nabla s_1 \; do \\
    \Leftrightarrow && \dot{\bar{s}}_2 
    =& \frac{1}{|\Omega_2|}  \biggl( \int_{\Omega_2} 
    - \bigl(1-A ( \bar{n}_2 + \underbrace{\tilde{n}_2}_{\mathclap{\text{small assuming homogeneous mixing}}} +n_0)^{-1}\bigr) \bigl(
    \bar{s}_2(\beta_e \bar{e}_2 +\beta_i \bar{i}_2) 
    + \underbrace{\tilde{s}_2 (\beta_e\bar{e}_2 + \beta_i\bar{i}_2) 
    + \bar{s}_2 (\beta_e\tilde{e}_2 + \beta_i\tilde{i}_2)}_{=0 \text{ if integrated}}
    \\
   && &+ \underbrace{\tilde{s}_2 (\beta_e\tilde{e}_2 + \beta_i\tilde{i}_2)}_{\mathclap{\text{small assuming homogeneous mixing}}} 
   \bigr) \; dx - \int_{\Gamma}\nu_1^T D \nabla s_1 \; do \biggr) \\
    && \approx& - \bar{s}_2 (\beta_e \bar{e}_2 +\beta_i \bar{i}_2) \frac{1}{|\Omega_2|} \int_{\Omega_2} \left(1-\frac{A}{\bar{n}_2+n_0}\right) \; dx
    - \frac{1}{|\Omega_2|} \int_{\Gamma}\nu_1^T D \nabla s_1 \; do \\
    && =& - \left(1-\frac{A}{\bar{n}_2+n_0}\right)\bar{s}_2 (\beta_e \bar{e}_2 +\beta_i \bar{i}_2) 
    - \frac{1}{|\Omega_2|} \int_{\Gamma}\nu_1^T D \nabla s_1 \; do.
\end{align*}
Thus, we obtain an ODE for $\bar{s}_2$ with a right hand side contribution from the flux of the PDE solution $s_1$ across $\Gamma$. The same derivation applies to $i$ and $e$.
Consequently, the mean subpopulation density fractions in $\Omega_2$ evolve according to the ODE model
\begin{align*}
 \dot{\bar{s}}_2 &\approx - \left(1-\frac{A}{\bar{n}_2+n_0}\right) \bar{s}_2 (\beta_e \bar{e}_2 +\beta_i \bar{i}_2) - \frac{1}{|\Omega_2|} \int_\Gamma  \nu_1^T D\nabla s_1 \; do, \\
 \dot{\bar{e}}_2 &\approx \left(1-\frac{A}{\bar{n}_2+n_0}\right)\bar{s}_2 (\beta_e \bar{e}_2 +\beta_i \bar{i}_2) - \sigma \bar{e}_2 - \phi_e \bar{e}_2 - \frac{1}{|\Omega_2|}\int_\Gamma \nu_1^T D\nabla e_1 \; do, \\
 \dot{\bar{i}}_2 &\approx \sigma \bar{e}_2 - \phi_i \bar{i}_2 - \frac{1}{|\Omega_2|}\int_\Gamma \nu_1^T D\nabla i_1 \; do, \\
 \dot{\bar{r}}_2 &\approx \phi_i \bar{i}_2 + \phi_e \bar{e}_2 - \frac{1}{|\Omega_2|} \int_\Gamma \nu_1^T D\nabla r_1 \; do.
\end{align*}
The statement in \cite{Guglielmi2022Delay} about the Allee term not being included in the ODE model due to the lack of spatial variation in population density overlooks an important aspect. The Allee term has a maximum value of 1 and is typically smaller, even negative in their case. Therefore, it is crucial to consider the Allee term as a mean across the area in the ODE model, ensuring that the infection numbers in the ODE model do not explode and are more comparable to the results of the PDE model.

Note that we already assumed $\tilde{s}_2$, $\tilde{e}_2$, $\tilde{i}_2$ and $\tilde{r}_2$ to be approximately zero, since they become vanishingly small as $D$ becomes large. 
Further, assuming that the total number of susceptibles $S_2(0)$, exposed $E_2(0)$, infectious $I_2(0)$ and removed $R_2(0)$
of the ODE model at time $t=0$ are given, and therefore the (constant) total population size $N_2$, we can set the initial conditions as 
\begin{align*}
    \int_{\Omega_2} s_2(x,0)n_2(x) \; dx = S_2(0) &&
    \Leftrightarrow &&\int_{\Omega_2} \bar{s}_2(0)n_2(x) \; dx \approx S_2(0) &&
    \Leftrightarrow && \bar{s}_2(0) \approx \frac{S_2(0)}{N_2} .
\end{align*}

Turning to $\Omega_1$, the density continuity~\eqref{eq:continuity-density} implies the Dirichlet boundary condition $s_1 = s_2$ a.e. on $\Gamma$ of the hybrid model.
For the sake of completeness, the hybrid model is presented in the following, under the assumption of homogeneous mixing:
\begin{equation} \label{Eq:hybrid_model}
    \begin{aligned} 
        \partial_t s_1 &= \divergence (D \nabla s_1) - \left(1-\frac{A}{n_1+n_0}\right) s_1 (\beta_e e_1 + \beta_i i_1) && \text{in $\Omega_1 \times ]0,T[$,} \\
        \partial_t e_1 &= \divergence(D \nabla e_1) + \left(1-\frac{A}{n_1+n_0}\right)s_1(\beta_e e_1 + \beta_i i_1) - \sigma e_1 - \phi_e e_1 && \text{in $\Omega_1 \times ]0,T[$,} \\
        \partial_t i_1 &= \divergence(D\nabla i_1) + \sigma e_1 - \phi_i i_1 && \text{in $\Omega_1 \times ]0,T[$,}\\
         \partial_t r_1 &= \divergence(D\nabla r_1) +  \phi_e e_1 + \phi_i i_1 && \text{in $\Omega_1 \times ]0,T[$,}\\
        \dot{s}_2 &= - \left(1-\frac{A}{n_2+n_0}\right) s_2 (\beta_e e_2 +\beta_i i_2) - \frac{1}{|\Omega_2|} \int_\Gamma  \nu_1^T D\nabla s_1 \; do && \text{in $]0,T[$,} \\
         \dot{e}_2 &= \left(1-\frac{A}{n_2+n_0}\right)s_2 (\beta_e e_2 +\beta_i i_2) - \sigma e_2 - \phi_e e_2 - \frac{1}{|\Omega_2|}\int_\Gamma \nu_1^T D\nabla e_1 \; do && \text{in $]0,T[$,} \\
         \dot{i}_2 &= \sigma e_2 - \phi_i i_2 - \frac{1}{|\Omega_2|}\int_\Gamma \nu_1^T D\nabla i_1 \; do && \text{in $]0,T[$,} \\
         \dot{r}_2 &= \phi_i i_2 + \phi_e e_2 - \frac{1}{|\Omega_2|} \int_\Gamma \nu_1^T D\nabla r_1 \; do && \text{in $]0,T[$,} \\
         s_1 &= s_2, \qquad e_1 = e_2, \qquad i_1 = i_2, \qquad r_1 = r_2 && \text{a.e. on $\Gamma \times ]0,T[$,} \\ 
         \nu_1^T D \nabla s_1 &= \nu_1^T D \nabla e_1 = \nu_1^T D \nabla i_1 = \nu_1^T D \nabla r_1 = 0 && \text{on $(\partial \Omega_1 \setminus \Gamma) \times ]0,T[$,} \\ 
         s_1(x,0) &= \sum_{\tilde{l}=1}^{\tilde{L}} w_{s}^{(\tilde{l})} G^{(\tilde{l})}(x), \qquad
         e_1(x,0) = \sum_{\tilde{l}=1}^{\tilde{L}} w_{e}^{(\tilde{l})} G^{(\tilde{l})}(x) && \text{in $\Omega_1$,} \\
         i_1(x,0) &= \sum_{\tilde{l}=1}^{\tilde{L}} w_{i}^{(\tilde{l})} G^{(\tilde{l})}(x), \qquad 
         r_1(x,0) = \sum_{\tilde{l}=1}^{\tilde{L}} w_{r}^{(\tilde{l})} G^{(\tilde{l})}(x)  && \text{in $\Omega_1$,} \\
         s_2(0) &= \frac{S_2(0)}{N_2}, \qquad e_2(0) = \frac{E_2(0)}{N_2}, \qquad i_2(0) = \frac{I_2(0)}{N_2}, \qquad r_2(0) = \frac{R_2(0)}{N_2}.
    \end{aligned}
\end{equation}

\section{Parameter Identification}

In the following, we outline the methodology for calibrating the model parameters to align closely with observed epidemiological data. This process is crucial for ensuring that our simulations are able to accurately reflect the dynamics of a given infectious disease spread. The fitting is based on the (slightly modified) Levenberg-Marquardt algorithm. We focus on identifying the parameters of the models using the PDE model, where the number of infectious individuals serves as our primary target. We specifically analyze the fraction of population density of infectious individuals $i(x,t;p)$ with parameters $p =\begin{pmatrix}
    \sigma & \phi_e & \phi_i & \beta_i 
\end{pmatrix}^T$, assuming $\beta_e=\beta_i$ like in \cite{Quelle5,viguerie2020diffusion}. 
Note that the fraction of population density of infectious $i$ depends only on three out of four parameters explicitly, but on all parameters implicitly.
Further, let $\hat{\mathcal{I}}_{k}^{(l)}$ be the target value of 
the number of infectious individuals
\begin{align*}
    \mathcal{I}^{(l)}(t_k) = \int_{\Omega^{(l)}} i(x, t_k)n(x) \; dx \in \mathbb{R}
\end{align*} 
at time $t_k$, $k=1,\dots,K$ for subdomain $\Omega^{(l)}$ with $l=1,\dots,L$, where
$i(x, t_k)$ is the fraction of population density of infectious at location $x$ at time $t_k$ with $k=1,\dots,K$. \\ \\ 
Applying the Levenberg-Marquardt algorithm (see Algorithm \ref{alg:Levenberg_Marquardt} in Appendix  \ref{appendix:Levenberg_Marquardt}, similar to \cite{Gavin2019Levenberg,Quelle10,Zhao2023Levenberg}), we have to solve
\begin{align*}  
    \bigl( \text{DF}(p^{(l)})^T \text{DF}(p^{(l)}) + \lambda_1 I \bigr) \Delta p^{(l)} = - \text{DF}(p^{(l)})^T F(p^{(l)}) 
\end{align*}
with identity matrix $I$, first damping coefficient $\lambda_1 \geq 0$, 
new parameters
\begin{align*}
    p^{(l+1)} &=p^{(l)} + \lambda_2 \Delta p^{(l)}
\end{align*}
with second damping coefficient $\lambda_2 \in (0,1]$, 
error function 
\begin{align*}
F(p) &= \bigl( \mathcal{I}^{(l)}(t_k;p) - \hat{\mathcal{I}}_{k}^{(l)} \bigr)_{k=1, l=1}^{K,L} \in \mathbb{R}^{KL} 
\end{align*}
and their Jacobian matrix
\begin{align*}
\text{DF}(p) &= \biggl( \frac{\partial \mathcal{I}^{(l)}(t_k;p_m)}{\partial p_m}  \biggr)_{k=1, l=1, m=1}^{K,L,M} \in \mathbb{R}^{KL \times M},
\end{align*}
where $p = \bigl( p_m \bigr)_{m=1}^M$ are our $M=4$ parameters we want to optimize. 
The damping coefficient $\lambda_1$ is chosen so that the matrix $\bigl( \text{DF}(p^{(l)})^T \text{DF}(p^{(l)}) + \lambda_1 I \bigr)$ is strictly diagonally dominant. Since the matrix $\bigl( \text{DF}(p^{(l)})^T \text{DF}(p^{(l)}) + \lambda_1 I \bigr)$ is symmetric and has real values, we therefore get a positive definite matrix. For $\lambda_1=0$ we apply the Gauss-Newton method and for $\lambda_1 \gg 0$ the gradient descent method \cite{Quelle10}.
\\ \\
Now, we want to determine $\frac{\partial \mathcal{I}^{(l)}(t_k;p_m)}{\partial p_m}$. We know that
\begin{align*}
    \frac{\partial \mathcal{I}^{(l)}(t_k;p_m)}{\partial p_m}
    &= \int_{\Omega^{(l)}}  \frac{\partial \bigl(i(x, t_k;p_m)n(x;p_m) \bigr)}{\partial p_m}\; dx \\
    &=\int_{\Omega^{(l)}}  n(x;p_m)\frac{\partial i(x, t_k;p_m)}{\partial p_m}  + i(x, t_k;p_m)\frac{\partial n(x;p_m)}{\partial p_m} \; dx  \\
    &= \int_{\Omega^{(l)}}  n(x;p_m)\frac{\partial i(x, t_k;p_m)}{\partial p_m}\; dx,
\end{align*}
i.e., we need to determine the derivative of the fractional population density $i$ with respect to the parameters $p_m \; \forall m=1,\dots,M$.
Recall our PDE for infectious 
\begin{align*} 
    \partial_t i &= \divergence \bigl( D \nabla i  \bigr) + \sigma e - \phi_i i \\ 
    &= \divergence \bigl( D \nabla i \bigr) + p_1 e - p_3 i. 
\end{align*} 

Let $f_{i}(\divergence (D\nabla y_n)_{n=1}^N,y;p) \coloneqq \divergence \bigl( D \nabla i \bigr) + p_1 e - p_3 i$
be the right hand side of the PDE for infectious
with the solution of the PDE model $y= \bigl( y_n \bigr)_{n=1}^N = \begin{pmatrix}
    s &  e &  i 
\end{pmatrix}^T \in \mathbb{R}^3$.
Differentiating $f_{i}(\divergence (D\nabla y_n)_{n=1}^N, y;p)$ with respect to $p_m$, $m \in \{1,\dots,4\}$ results in
\begin{align*} 
    \partial_{p_m} \bigl( \partial_t i \bigr) &= \partial_{p_m} \bigl( f_{i}(\divergence (D\nabla y_n)_{n=1}^N, y;p) \bigr) \\
    &= \sum_{n=1}^N \partial_{\divergence (D\nabla y_n)} f_{i}(\divergence (D\nabla y_n)_{n=1}^N, y;p) \cdot \partial_{p_m} \bigl( \divergence (D\nabla y_n) \bigr) \\
    &\qquad +\partial_{y_n} f_{i}(\divergence (D\nabla y_n)_{n=1}^N, y;p) \cdot \partial_{p_m} y_n + \partial_{p_m} f_{i}(\divergence (D\nabla y_n)_{n=1}^N, y;p)  \in \mathbb{R}
\end{align*}
for $m=1,\dots,M$
with
\begin{align*} 
\partial_{p_m} \bigl( \partial_t i \bigr) &=
    \sum_{n=1}^N \partial_{\divergence (D\nabla y_n)} f_{i} \cdot \partial_{p_m} \bigl( \divergence (D\nabla y_n) \bigr) 
    + \partial_{y_n} f_{i} \cdot \partial_{p_m} y_n + \partial_{p_m} f_{i} \\
    &= \underbrace{\partial_{\divergence (D\nabla i)} f_{i}}_{=1} \cdot \partial_{p_m} \bigl( \divergence (D\nabla i) \bigr) + 
    \underbrace{\partial_{s} f_{i}}_{=0} \cdot \partial_{p_m} s +
    \underbrace{\partial_{e} f_{i}}_{=p_1} \cdot \partial_{p_m} e +
    \underbrace{\partial_{i} f_{i}}_{=-p_3} \cdot \partial_{p_m} i + \partial_{p_m} f_{i}
    \\
    &= \divergence (D \nabla \bigl(\partial_{p_m} i \bigr)) 
    + p_1 \partial_{p_m} e - p_3 \partial_{p_m} i +
    \delta_{m,1}e - \delta_{m,3}i
\end{align*}
in terms of Kronecker's delta $\delta_{j,k}$.

Now, considering the differential equation for exposed 
\begin{align*}
    \partial_t e = f_{e} &= \divergence (D \nabla e) 
    + \biggl(1-\frac{A}{n_l+n_0} \biggr) s (\beta_i i + \beta_e e) - \sigma e - \phi_e e \\
    &= \divergence (D \nabla e) + \biggl(1-\frac{A}{n_l+n_0} \biggr) p_4 s (i+e) - p_1 e - p_2 e,
\end{align*}
we have to solve the following DEs 
\begin{align*}
\partial_{p_m}  \bigl( \partial_t e \bigr)
&= \divergence(D\nabla\partial_{p_m} e) + \left(1-\frac{A}{n_l+n_0}\right)\big(\delta_{m,4} s (e+i) + p_4 (e+i)\partial_{p_m}s + p_4 s (\partial_{p_m}e+\partial_{p_m}i)\big) \\
&\quad - p_1\partial_{p_m}e - p_2\partial_{p_m}e - \delta_{m,1}e - \delta_{m,2}e.
\end{align*}

Moving on to the DE for susceptibles
\begin{align*}
    \partial_t s = f_{s} &= \divergence(D \nabla s) - \biggl(1-\frac{A}{n_l+n_0} \biggr) s (\beta_i i + \beta_e e) \\
    &= \divergence(D \nabla s) - \biggl(1-\frac{A}{n_l+n_0} \biggr) p_4 s (i + e),
\end{align*}
we must solve the subsequent DE
\begin{align*}
    \partial_{p_m} \bigl( \partial_t s \bigr) 
    &= \divergence(D\nabla\partial_{p_m}s) - \bigl(\divergence (D \nabla y_n) \bigr)\big(\delta_{m,4} s(e+i) + p_4 (e+i) \partial_{p_m}s + p_4s(\partial_{p_m}e + \partial_{p_m}i)\big).
\end{align*}

As the initial conditions of the PDE system (\ref{PDE-model}) do not depend on the parameters $p$, the initial conditions of the variational equations \cite{dierkes2013parameter} or sensitivity equations \cite{viguerie2020diffusion} with solution
$\begin{pmatrix}
    \partial_{p_m} s & \partial_{p_m} e & \partial_{p_m} i
\end{pmatrix}^T$
are zero.
The boundary conditions of the sensitivity equations are the same as for the PDE system (\ref{PDE-model}).
We solve the system of PDEs and variational equations simultaneously. This procedure is often termed as \textit{forward sensitivity analysis} \cite{mester2022differential}.
%

\section{Implementation}

In this section, we detail the technical processes and computational methods, and tools underlying the implementation of our hybrid epidemiological model.

To start, ensuring positive initial conditions for the PDE system is foundational. We impose the constraint $w_{y_n} \geq 0$ for $y_n = s,e,i,r$ and solve the corresponding constraint optimization problem using Sequential Least Squares Programming (SLSQP, as implemented in SciPy 1.12.0 using the \textit{scipy.optimize.minimize} function). 

Building on this initial setting, we proceed to calibrate the initial parameters, which is crucial for capturing the nuanced dynamics of disease spread. This calibration is achieved by fitting the fractional infectious density of the derived ODE model by using different patches, each patch representing one subdomain. Notably, inter-patch transmission is excluded from this phase to simplify the fitting process. This fitting step is done by minimizing the norm between target values and our numerical solution. The numerical solution is computed using finite differences with time step size $\Delta t = 0.1$, again by using SLSQP. 

Switching now to C++, our approach to solving the differential equations employs a linear implicit Euler method for temporal discretization. For the spatial discretization, we utilize a first-order finite element method.

Further refining the model, parameter identification for the PDE components employs a slightly modified version of the Levenberg-Marquardt algorithm, as detailed in Appendix \ref{appendix:Levenberg_Marquardt}. This step aims to minimize the weighted error $e_w^{(l)} = ||w_F(\Omega) F(p^{(l)})||_2$ where the weights vary by domain importance, amplifying the significance of Very Important Domains (VIDs), such as Milan (see below). The weight vector $w_F(\Omega)$ is defined as: 

\begin{align*}
    w_F(\Omega) = \begin{cases} 
    2, & \text{if }  \Omega    \text{ is VID,} \\ 
     1, & \text{if }  \Omega\text{ is not VID.}
    \end{cases}
\end{align*} 

This weighted error strategy allows for simultaneous optimization of all parameters and also a focused adjustment of specific parameters, such as $\beta_e$. This approach is inspired by \cite{Quelle5,viguerie2020diffusion}. The optimization continues iteratively until the error falls below a predefined threshold or the maximum number of iterations is reached. Then, we proceed with the subsequent temporal segments. 

To support these computational efforts, the triangulation of the modeling domain was performed using Triangle \cite{Quelle7}, while the full-PDE and the hybrid PDE-ODE systems were solved using the finite element method implemented in the Kaskade7 software \cite{GoetschelSchielaWeiser2020}, using the Dune interface. Finally, the visualization of our results was conducted using ParaView \cite{Quelle6}, a tool well-suited for rendering and analyzing complex simulation data.

\section{Numerical Results}

In this section, we focus on a detailed exploration of our hybrid model's efficacy through a series of numerical experiments. Our aim is to evaluate the performance of the hybrid model in comparison to the conventional full-PDE model across a spectrum of scenarios, ranging from controlled synthetic environments to complex real-world settings. To lay the groundwork, we begin by selecting a simplistic geographical area as our testing ground, where we vary the division of the domain into PDE and ODE segments (see section \ref{sec:Experiments Based on Synthetic Data}). Here, we aim not just to compare the hybrid model with its full-PDE counterpart, but also to illustrate the trade-offs between computational efficiency and solution accuracy (see section \ref{sec:rectangle_accuracy_efficiency}).
This initial phase is crucial for quantifying the model's stability and its capacity to maintain precision across different configurations, thereby providing us with a foundational understanding of the model's capabilities.

Moving beyond synthetic experiments, we advance to the evaluation of the models against empirical data, bringing into focus the real-world applicability and validity of our modeling approach. Our investigation extends across two distinct locales: Lombardy, Italy, and Berlin, Germany for which real-world data is available (see section \ref{sec:Real-World Examples}). This validation approach not only tests the models' predictive precision but also their adaptability and responsiveness to real-world infectious disease dynamics. The results of these experiments will help us understand how our hybrid modeling approach can be used in the complex field of epidemiological forecasting and analysis.

\subsection{Experiments Based on Synthetic Data} 
\label{sec:Experiments Based on Synthetic Data}
In this set of experiments, we investigate the relationship between solution accuracy and computational efficiency within a partitioned rectangular domain $\Omega = (0,2) \times (0,1)$, 
dividing it into PDE and ODE regions with varying percentages. Our aim is to uncover potential trade-offs in achieving accuracy while maintaining computational efficiency. Additionally, we explore extreme cases within a 50-50 domain division, assessing the model's resilience when infected individuals are concentrated in one region and evaluating its response to boundary challenges. These experiments provide valuable insights into the model's behavior under different configurations, contributing to a comprehensive understanding of the system's performance.

\subsubsection{Trade-offs Between Solution Accuracy and Computational Efficiency} 
\label{sec:rectangle_accuracy_efficiency}
%
The PDE model's solution is assumed to represent an idealized benchmark against which the results from the ODE model and the hybrid PDE-ODE model are compared. Considering the entire domain for evaluation, any error arising from the PDE domain or the region $\Omega_1\subset\R^2$ of interest, modeled in spatial detail, is already incorporated through the coupling.
The error measurement is conducted as an absolute error, calculated as the sum of deviations from the benchmark.
The accuracy is determined using the formula
\begin{align*}
    \text{accuracy} = 1 - \frac{||\text{deviation}||_2}{\text{maximum possible deviation}}. 
\end{align*}
Here, the solution of the PDE model achieves an accuracy of 1, while an ODE model's solution attains an accuracy of 0, effectively reflecting the maximum possible deviation. \\ \\

\begin{table}[H] 
\centering %
\begin{tabular*}{0.7\textwidth}{
@{\extracolsep{\fill}} *{4}{c}} 
\multicolumn{1}{c}{percentage of ODE-domain}    &
\multicolumn{3}{c}{number of}    \\
\cmidrule{1-1}
\cmidrule{2-4}
           & vertices & edges & triangles \\
  \hline
  \midrule %
  \textcolor{white}{0}$0$   & $25,669$  & $71,004$ & $45,336$ \\
  \textcolor{white}{0}$12.5$ & $22,918$ & $63,251$ & $40,334$ \\
  $25$   & $20,113$  & $55,336$ & $35,224$ \\
  \textcolor{white}{0}$37.5$ & $17,336$ & $47,505$ & $30,170$ \\
  $50$   & $14,561$ & $39,680$ & $25,120$ \\
  \textcolor{white}{0}$62.5$ & $11,753$ & $31,756$ & $20,004$ \\
  $75$   & \textcolor{white}{0}$8,941$ & $23,820$ & $14,880$ \\
  \textcolor{white}{0}$87.5$ & \textcolor{white}{0}$6,177$ & $16,028$ & \textcolor{white}{0}$9,852$ \\
  $100$ \textcolor{white}{0} & $0$ \textcolor{white}{00} & $0$ \textcolor{white}{00}& $0$ \textcolor{white}{00}\\
\bottomrule 
\end{tabular*}  %
\caption{Information on mesh of PDE-domain of the hybrid model (\ref{Eq:hybrid_model}) and full-PDE model (\ref{PDE-model}). 
} 
\label{table:mesh_info_rectangle}%
\end{table}%

We want to partition the domain vertically and vary the sizes for the PDE-domain and the ODE-domain. We simultaneously adjust the proportion of infected and recovered individuals in their respective regions, along with the total population, during the proportional partitioning of the area. We scale the infected individuals to represent the exposed category. The fractional densities $s$, $e$, $i$ are initially set as constants across the region, and for the population density, we use Gaussian functions that remain constant along the x-axis with a variance of $0.1$. 
With this choice of the Gaussian function, we prevent disadvantages that may arise when partitioning the domain.
In Table \ref{table:parameters_rectangle}, you will find the parameters specific to the rectangular domain across various time intervals.
\begin{table}[H] 
\centering %
\begin{tabular*}{\textwidth}{
@{\extracolsep{\fill}} *{5}{c}} 
\multicolumn{1}{c}{parameter}    &
\multicolumn{4}{c}{value}    \\
\cmidrule{1-1}
\cmidrule{2-5}
           & time index 0 & time index 10 & time index 23 & time index 29 \\
  \hline
  \midrule %
    $\sigma$              & $2.6676e-02$  
  & - & - & $1.8747e-01$ \\
  $\phi_e$              & $0.0000e+00$ 
  & - & - & $3.2655e-11$ \\
  $\phi_i$              & $2.3310e-01$  
  & - & - & $2.1907e+00$ \\
  $\beta_i$, $\beta_e$  & $4.4202e-01$ & $2.0588e-01$ & $6.0352e-08$ & $2.2168e-01$ \\
  $D$  & $\tfrac{0.0435}{111.3^2} \cdot 1.0e+02 $ 
  & $\tfrac{0.0198}{111.3^2}  \cdot 1.0e+2$ & $\tfrac{0.0090}{111.3^2} \cdot 1.0e+2$ & $\tfrac{0.0075}{111.3^2} \cdot 1.0e+2$ \\
  $A$                   & $3.0e+07$ & - & - & - \\
  $n_0$                 & $4.5e+07$ & - & - & - \\
\bottomrule 
\end{tabular*}  %
\caption{Parameters of PDE model (\ref{PDE-model}) and of the hybrid model (\ref{Eq:hybrid_model}) for rectangular domain. 
The symbol “-” means that the parameter has not changed over time. } 
\label{table:parameters_rectangle}%
\end{table}%
\noindent
The infectious density over time can be observed in the following Fig. \ref{fig:rectangle_infectious_density_PDE}. Given that the fractional densities $s$, $e$, $i$ are initially constant and the population density $n$ follows a Gaussian distribution, the resulting infectious density $i \cdot n$ will also follow a Gaussian distribution for $t=0$. The Gaussian-shaped spread of infectious disease is evident due to the constant Gaussian-shaped population density.
\begin{figure}[h!]
        \includegraphics[width=0.23\textwidth]{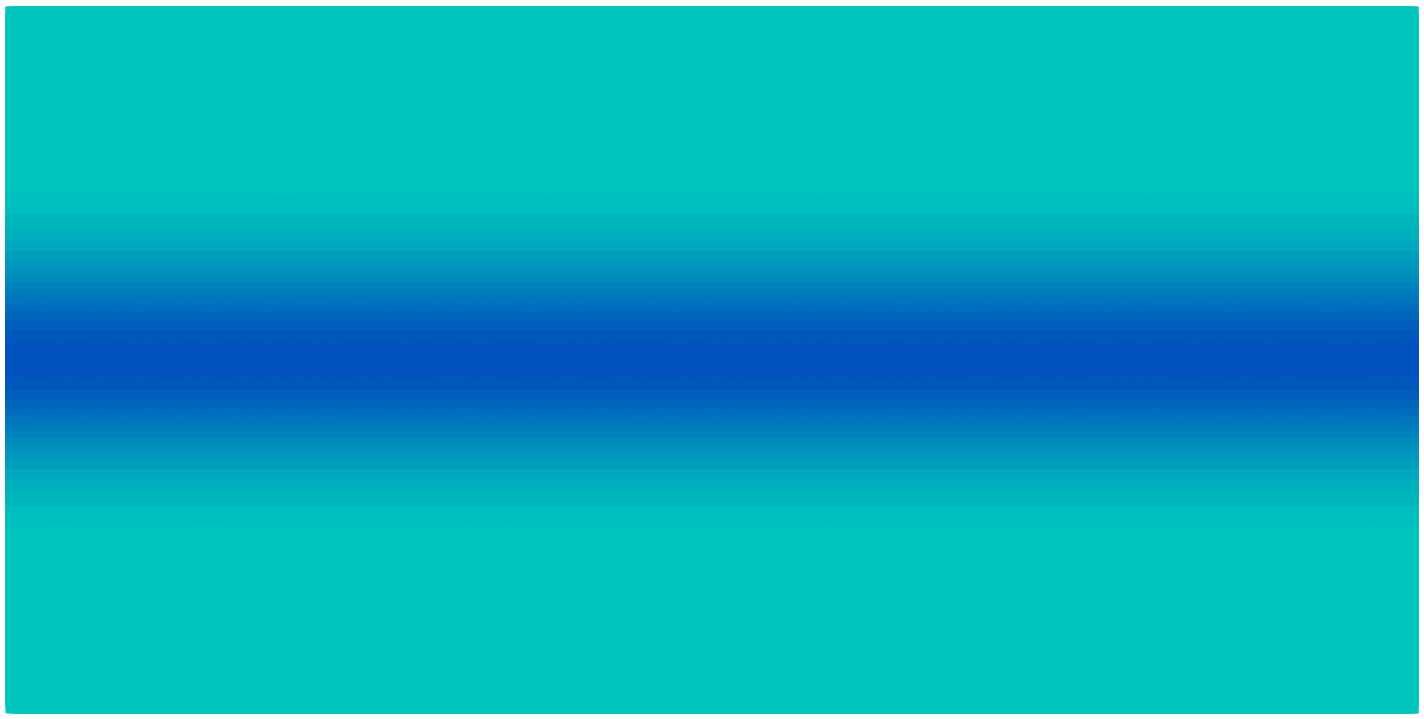} 
        \includegraphics[width=0.23\textwidth]{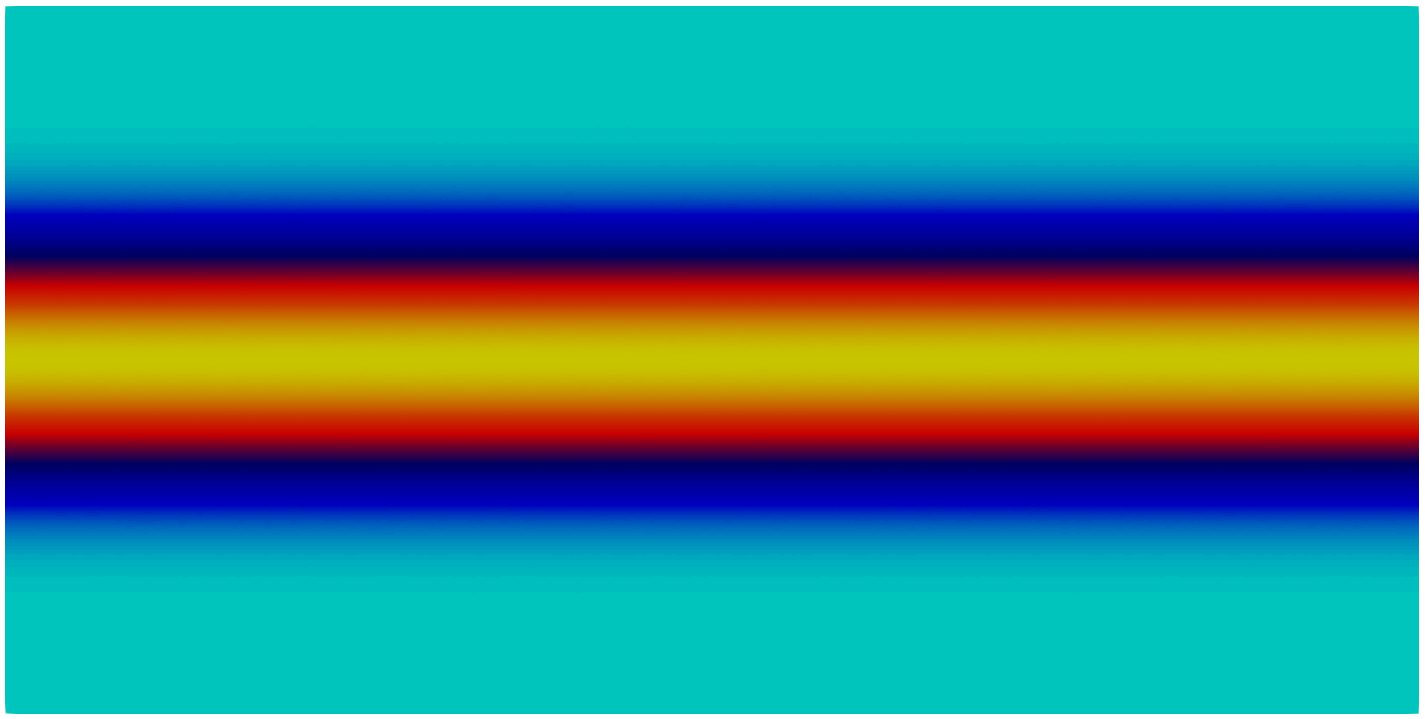}
         \includegraphics[width=0.23\textwidth]{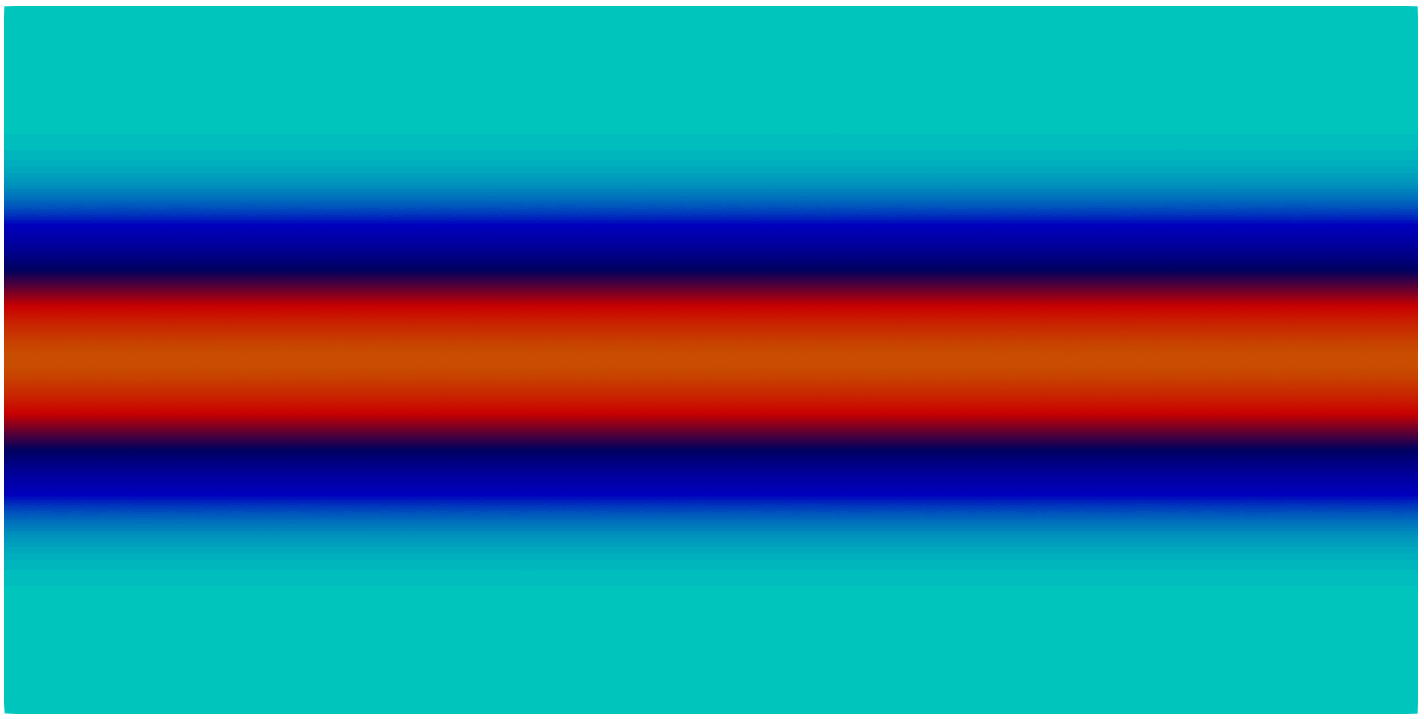} 
         \includegraphics[width=0.23\textwidth]{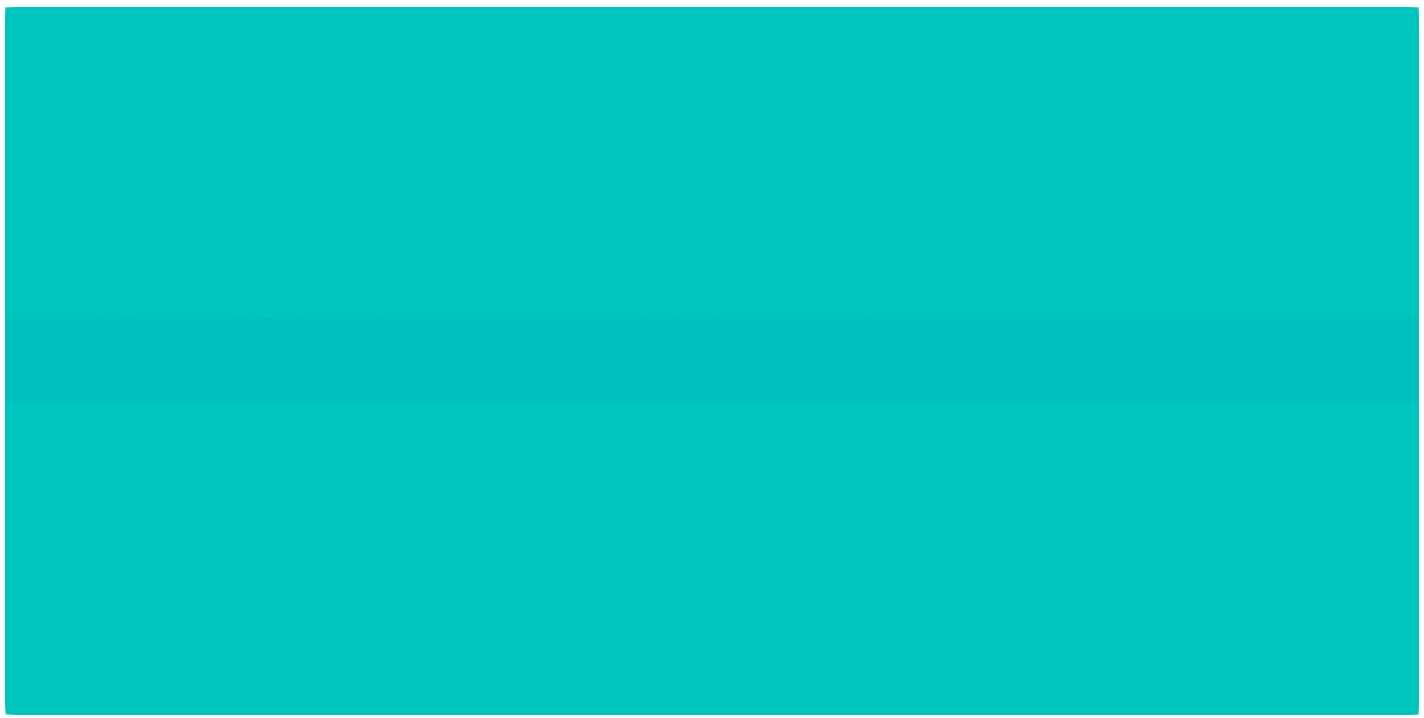}
         \includegraphics[width=0.035\textwidth]{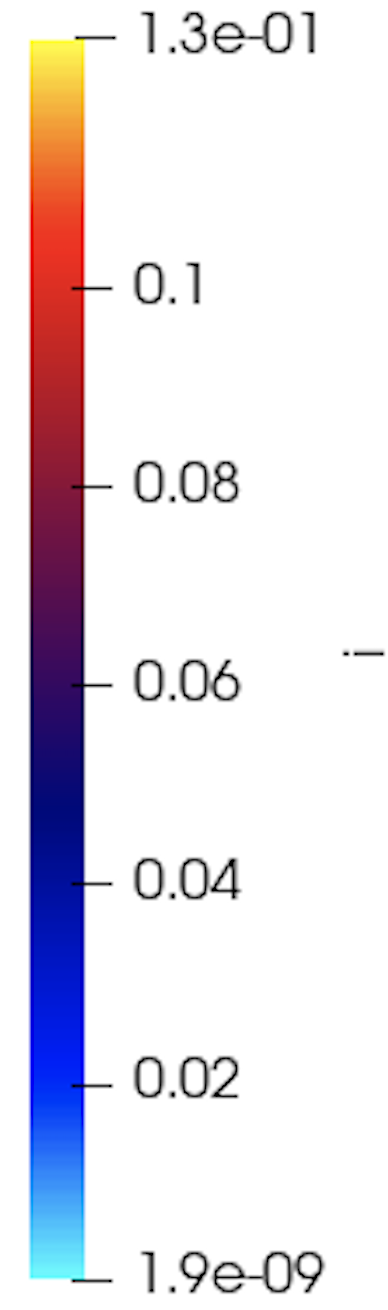}
\caption{Infectious density of PDE model (\ref{PDE-model}) for times $t \in \{0,29,30,59\}$ in rectangular domain. The Gaussian-shaped spread of infection is evident due to the Gaussian-shaped population density.} 
\label{fig:rectangle_infectious_density_PDE}
\end{figure} 
\noindent

\begin{figure}[h!]
    \begin{minipage}{0.49\textwidth}
        \includegraphics[width=\textwidth]{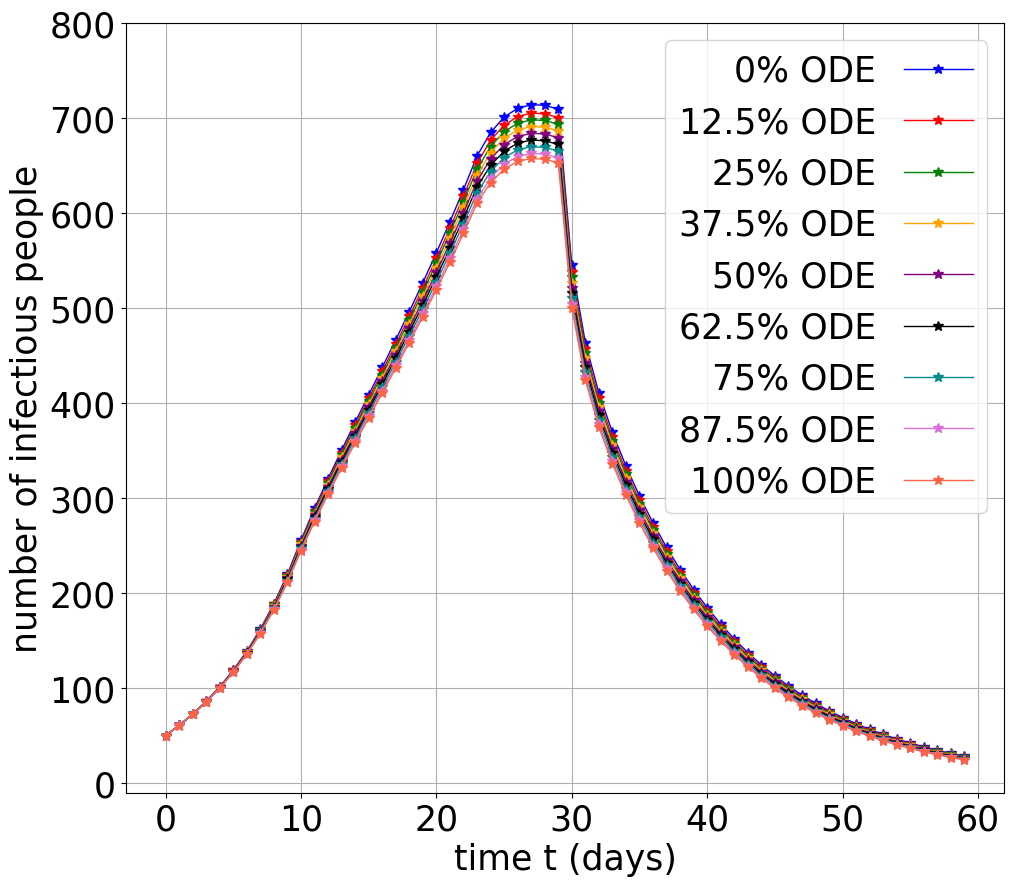} 
    \end{minipage}
    \begin{minipage}{0.49\textwidth}
        \includegraphics[width=\textwidth]{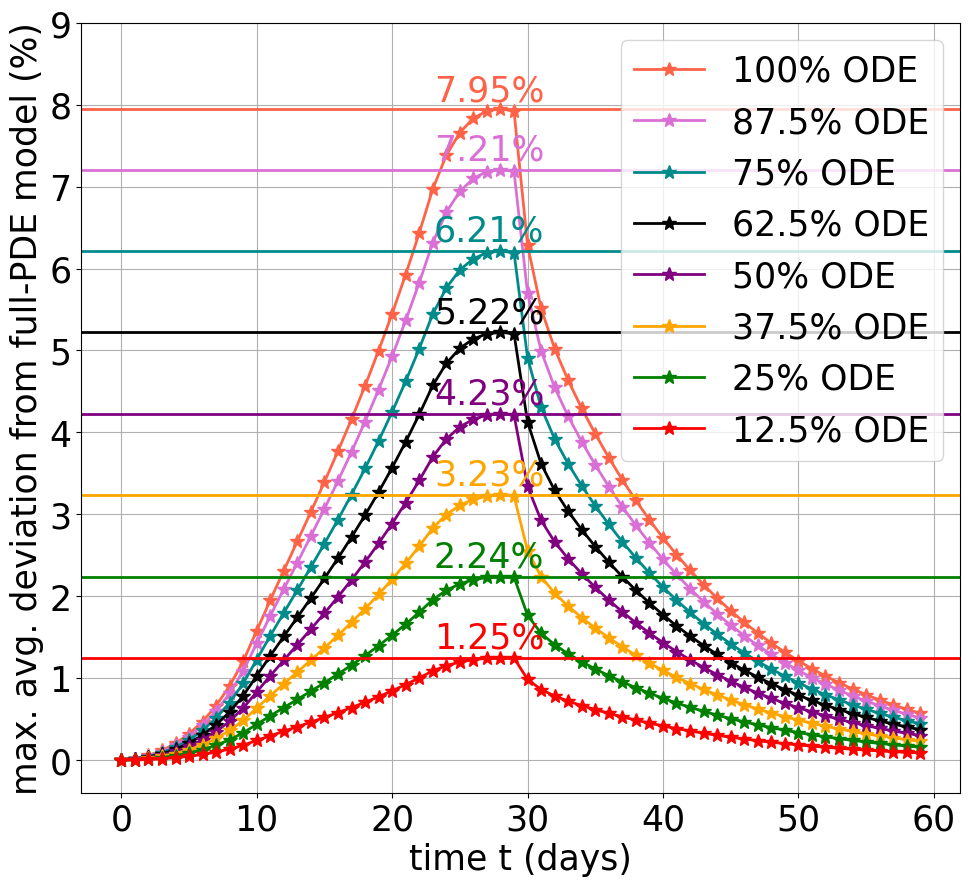}
    \end{minipage}
\caption{Number of infectious individuals of full-PDE and hybrid models in rectangular domain for times $t \in [0,59]$ (left) and maximum average deviation from full-PDE model (right). The infectious count decreases with a higher percentage of the ODE domain. This is attributed to the inclusion of the Allee term, which spatially modifies the infection rate in the PDE domain while remaining constant in the ODE domain.} 
\label{fig:rectangle_varying_ode_percentages}
\end{figure}

\begin{figure}[h]
    \centering
    \begin{minipage}{0.49\textwidth}
    \vspace{-0.5cm}
        \includegraphics[width=\textwidth]{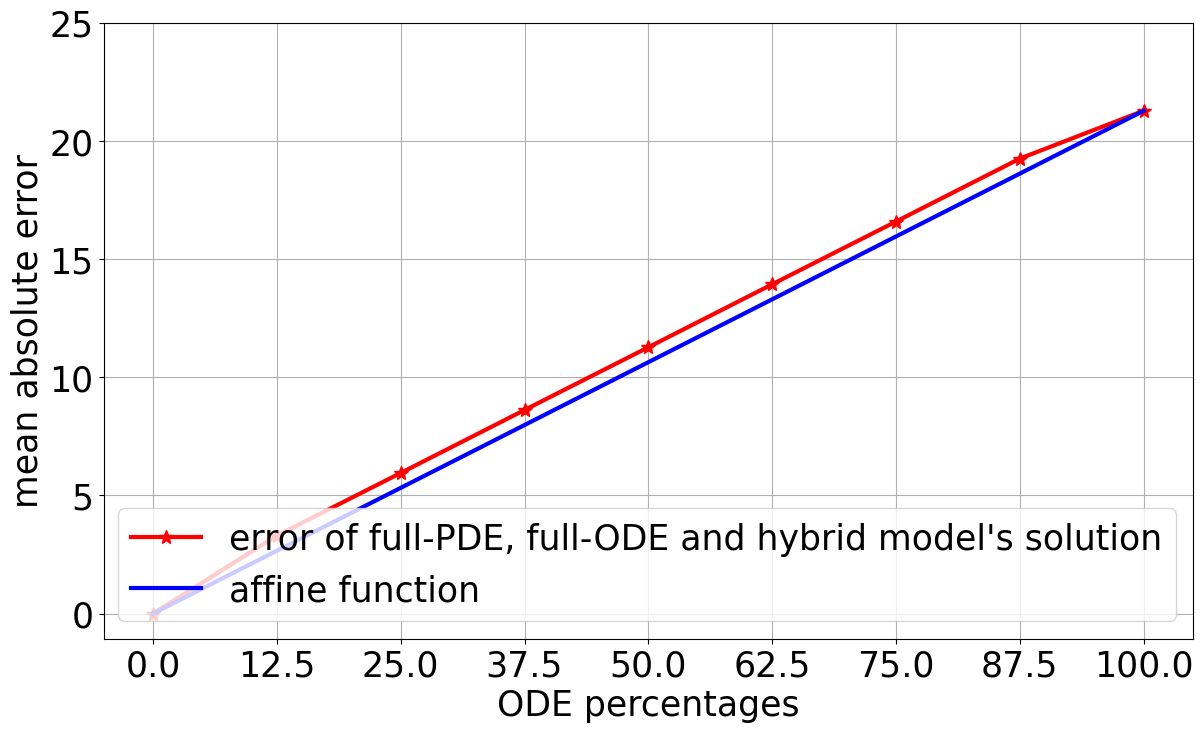}
        \caption{Accuracy (mean absolute error) of the full-PDE model (\ref{PDE-model}) and hybrid model (\ref{Eq:hybrid_model}) in rectangular domain.  We can observe an approximately linear growth.} 
        \label{fig:accuracy_rectangle}
    \end{minipage}
\begin{minipage}{0.49\textwidth}
\includegraphics[width=\textwidth]{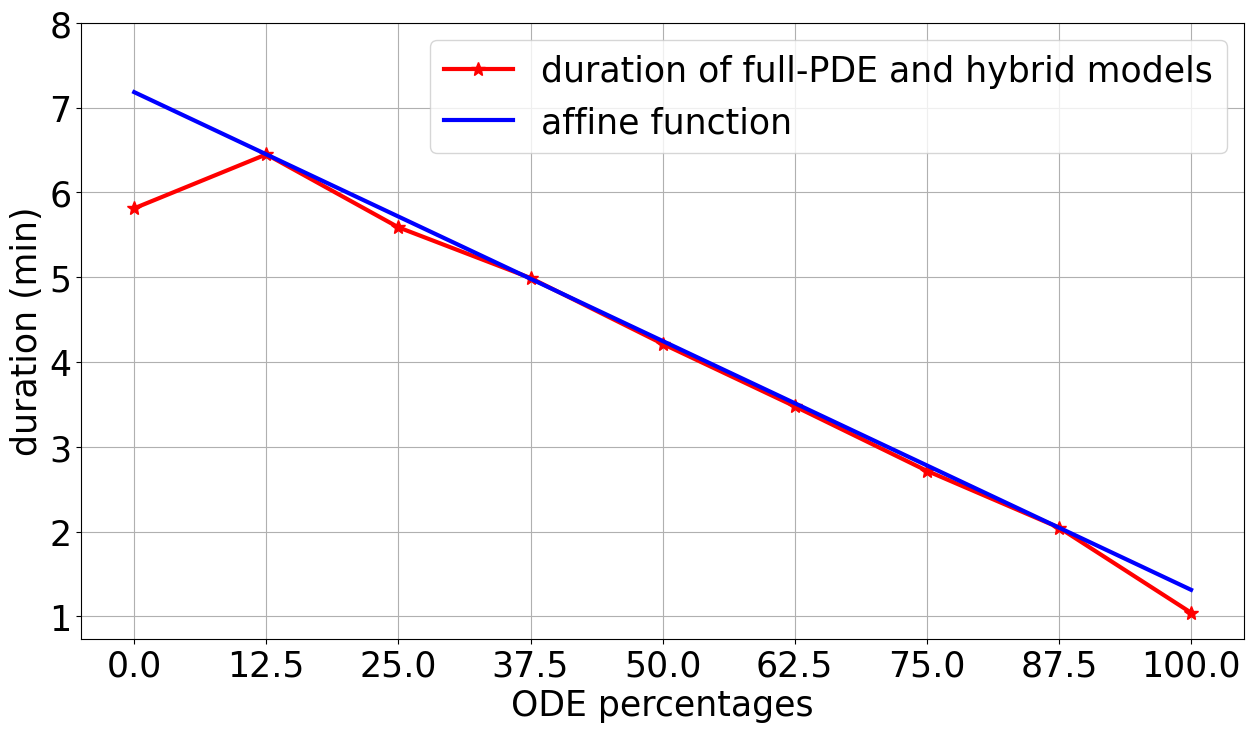}
\caption{Runtimes of the full-PDE (\ref{PDE-model}), hybrid (\ref{Eq:hybrid_model}) and full-ODE models in rectangular domain.
The runtime initially increases (/finally decreases) due to the additional (/absence of) implementation costs associated with the boundary conditions of the hybrid model.}
\label{fig:times_rectangle}
\end{minipage}
\end{figure}

Figure \ref{fig:rectangle_varying_ode_percentages} shows the number of infectious individuals for each time step, along with the maximum average deviation in the number of infectious individuals between the hybrid model and the full-ODE model when compared to the full-PDE model.
The maximum average deviation measures imply that the variations between hybrid models are around 0.99\%, as calculated by the differences in the right plot. In other words, a 12.5\% difference in area corresponds to an approximate 1\% error when comparing hybrid models.
A more notable error difference is observed for the hybrid model in the right plot with an ODE domain size of 87.5\% compared to the full-ODE model, which is around $7.95\%-7.21\%=0.74\%$.
Overall, these error differences suggest that the hybrid models tend to shift towards the direction of the ODE model, indicating that the ODE component of the hybrid model seems to exert a stronger weighting or influence on the boundary conditions between the regions, consequently affecting the model results of the PDE component of the hybrid model to a greater extent.
%
\noindent
\\ \\
Measuring the absolute deviation to the full-PDE model (see Fig. \ref{fig:accuracy_rectangle} and Table \ref{table:duration_accuracy_rectangle}), we can observe an approximately linear growth.
If the focus is on enhancing predictions or attaining greater precision in the results, it is recommended to choose the full-PDE model or consider a more advanced model, which may include additional compartments \cite{Fošnarič2022Extended}
or even an ABM \cite{Grossmann2021WhyODE}.

Figure \ref{fig:times_rectangle} illustrates the runtimes for the full-PDE, hybrid and full-ODE models in Kaskade7. On the left side, 100\% of the domain uses the PDE solver (0\% ODE), and on the right side, 100\% of the domain uses the ODE solver (0\% PDE). In the middle of the graph, we can see the results of the hybrid model with varying percentages of ODE and PDE domains. The results indicate that when only the hybrid model is considered, the runtime is approximately linear with respect to the percentage of ODE use, as the problem shifts from a spatially resolved PDE approach to a simpler ODE-based approach. 
Notably, the runtime initially increases due to the additional implementation costs associated with the boundary conditions of the hybrid model. For the same reason, it drops significantly at the far right of the graph, where the full-ODE model is applied and does not require boundary conditions.                                         
It should be noted that these runtimes are highly dependent on the implementation, and significant optimization could reduce these values. This figure is intended to provide an impression of how the domain division affects performance.
In this context, when comparing the program runtime with the accuracy of the results, it becomes evident that sacrificing a small amount of accuracy -- provided that at least 25\% of the area employs an ODE component -- can be a viable option for solving the model more efficiently. \\ \\
It is noticeable that the infectious count decreases with a higher percentage of the ODE domain, a trend typically not observed \cite{Abboubakar2023Reaction}, \cite{Schneckenreither2008Simulation}. 
This is attributed to the inclusion of the Allee term, which spatially modifies the infection rate in the PDE domain while remaining constant in the ODE domain.
Setting $A=0$ results in almost identical outcomes for all models, which is why we deliberately retained the Allee term. \\
\begin{table}[H] 
\centering %
\begin{tabular*}{0.4\textwidth}{@{\extracolsep{\fill}} *{2}{c}}
\toprule
 \; percentage of ODE-domain & error \\
  \midrule %
  \textcolor{white}{0}$0$              & \textcolor{white}{0}$0.00$ \\
  \textcolor{white}{0.}$12.5$          & \textcolor{white}{0}$3.30$ \\
  $25$                                 & \textcolor{white}{0}$5.95$ \\ 
  \textcolor{white}{0.}$37.5$          & \textcolor{white}{0}$8.61$ \\
  $50$                                 & $11.27$ \\
  \textcolor{white}{0.}$62.5$          & $13.93$ \\
  $75$                                 & $16.59$ \\
  \textcolor{white}{0.}$87.5$          & $19.25$ \\
  $100$                                & $21.27$ \\
\bottomrule 
\end{tabular*}
\caption{Accuracy (mean absolute error) of the full-PDE model (\ref{PDE-model}) and hybrid model (\ref{Eq:hybrid_model}).
} 
\label{table:duration_accuracy_rectangle}%
\end{table}%

\subsubsection{Extreme Cases} 
\label{sec:ExtremeCases}
\begin{figure}[h!]
    \centering
    \begin{minipage}{0.49\textwidth}
        \includegraphics[width=\textwidth]{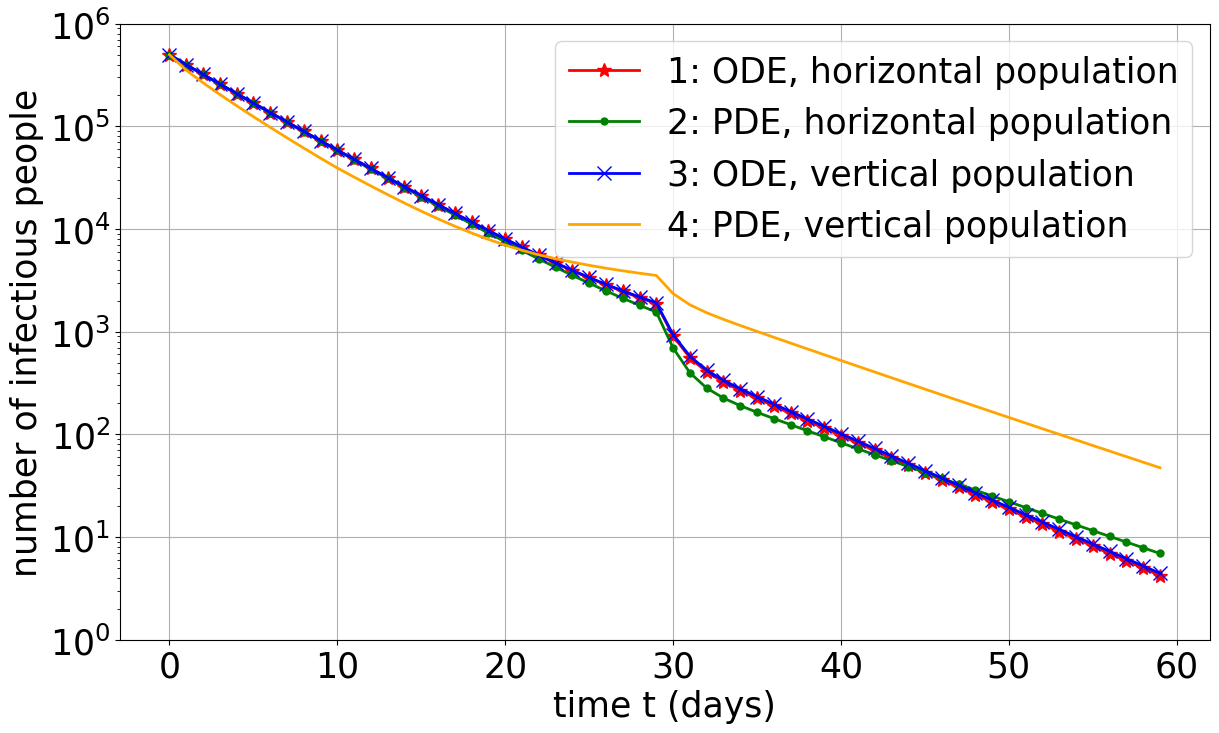}
    \end{minipage}
    \begin{minipage}{0.5\textwidth}
        \includegraphics[width=\textwidth]{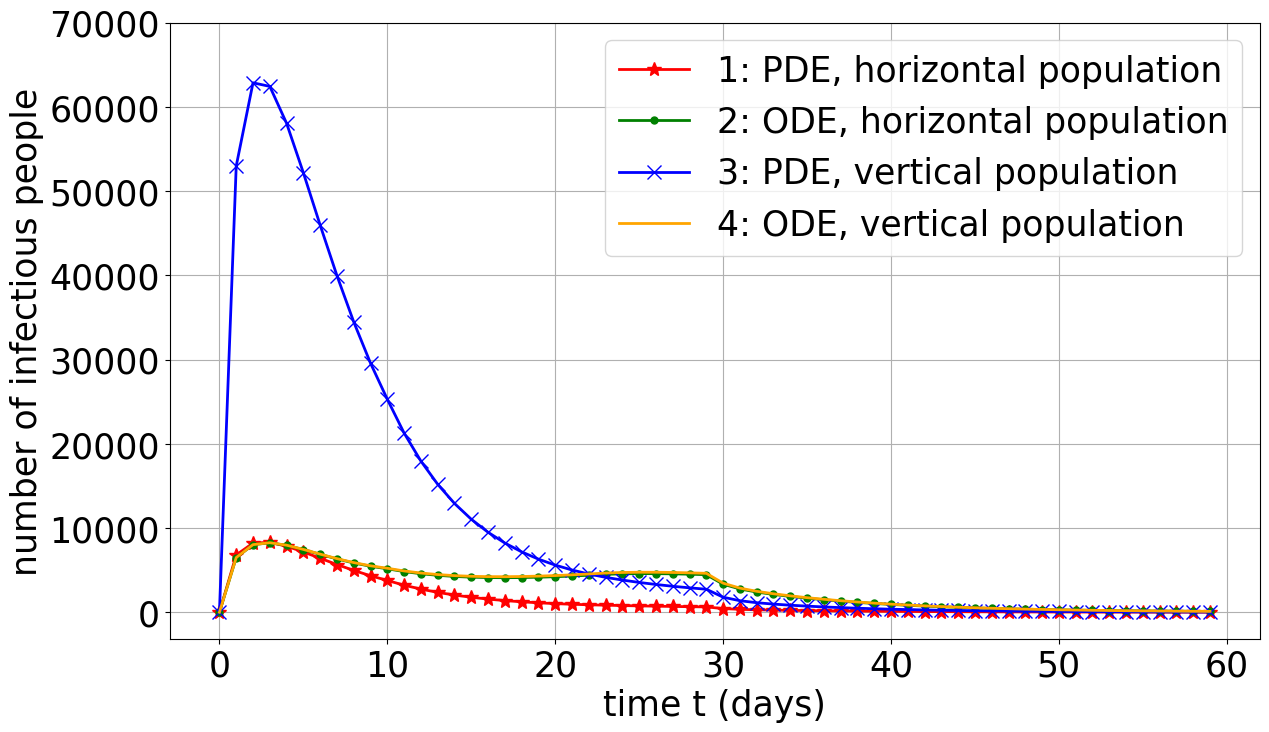}
    \end{minipage}
    \caption{Extreme cases of the hybrid model (\ref{Eq:hybrid_model}) in rectangular domain: the number of infectious individuals is initially equal to the total population number (left) and zero (right). The outcomes for the ODE region appear visually unaffected by the location of the population in the PDE domain, whereas the outcomes for the PDE region vary significantly.}
    \label{fig:extreme_cases_rectangle}
\end{figure}

\begin{figure}
    \centering
    \includegraphics[width=0.23\textwidth]{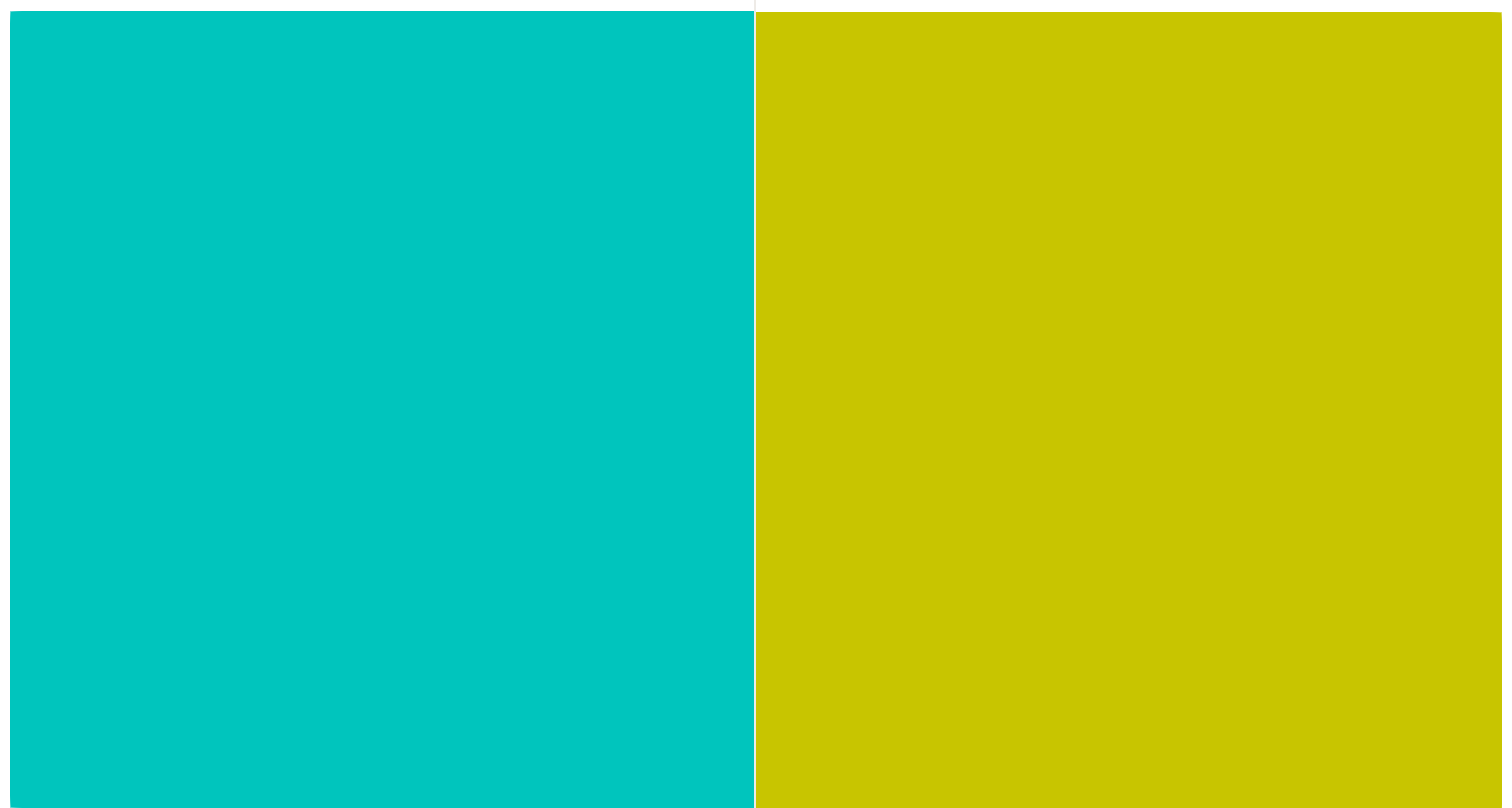}
    \includegraphics[width=0.23\textwidth]{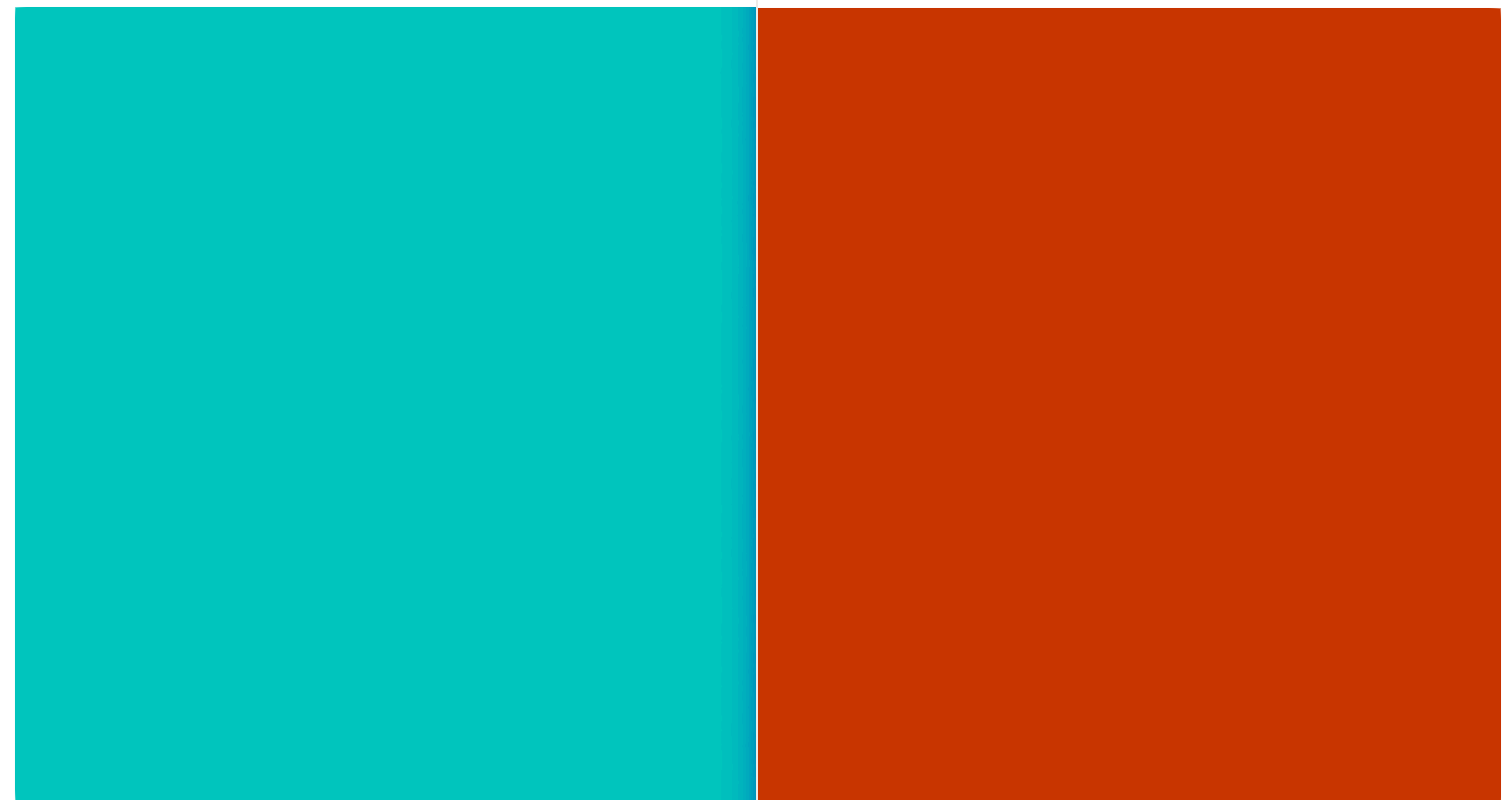}
    \includegraphics[width=0.23\textwidth]{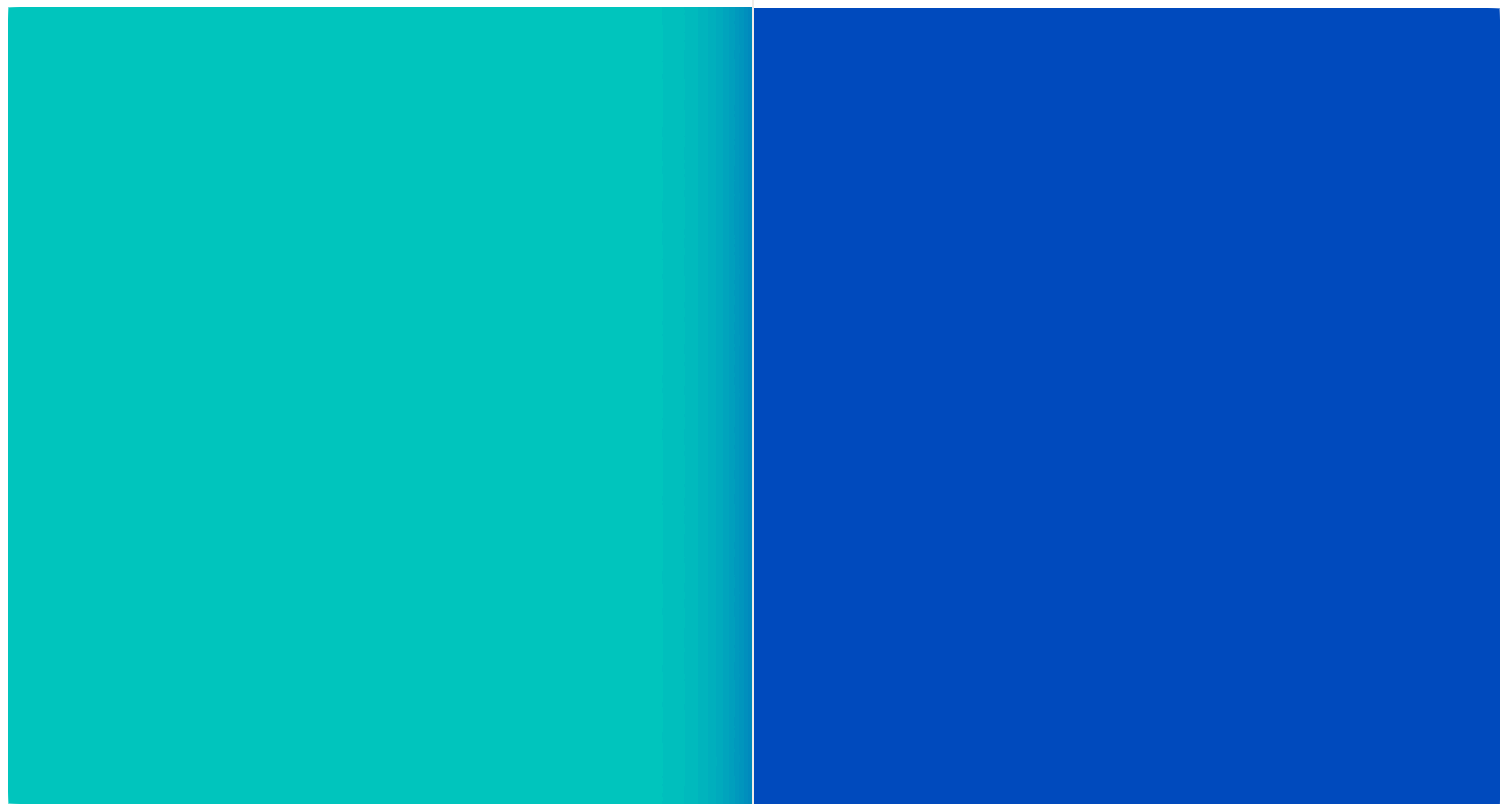}
    \includegraphics[width=0.23\textwidth]{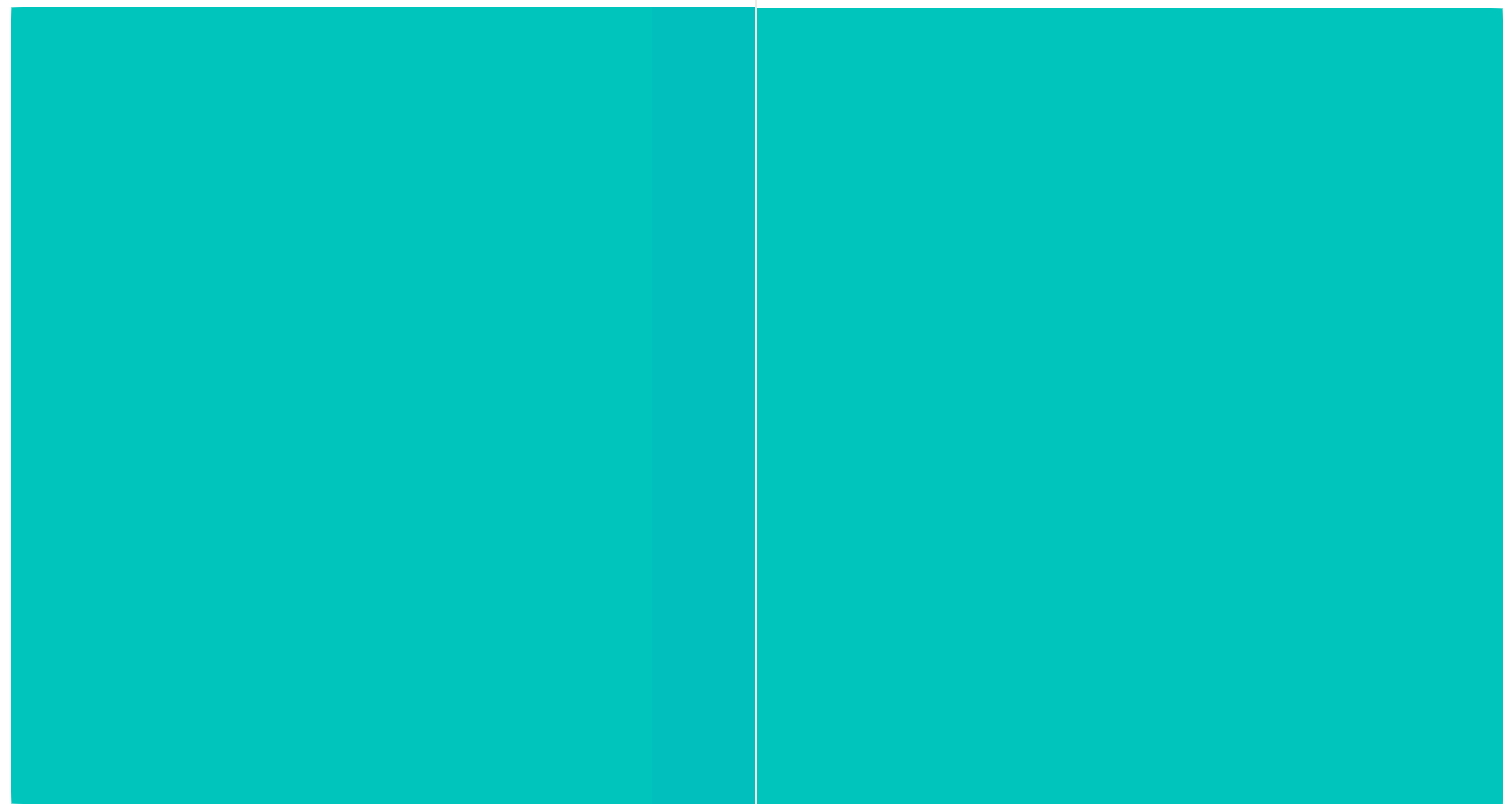}
    \includegraphics[width=0.06\textwidth]{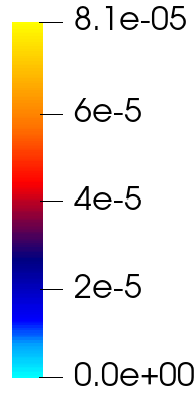}
    \caption{Fractional infectious density of hybrid model (\ref{Eq:hybrid_model}) for times $t \in \{0,2,10,29\}$ in rectangular domain. Initially, the number of infectious individuals for the PDE contribution of the hybrid model is equal to zero, while for the ODE contribution, it is equal to the population number. Due to the boundary conditions between the two regions,
there is an infection flow from the ODE region (right) to the PDE region (left).}
    \label{simulation_extreme_case_rectangle}
\end{figure}

Exploring extreme cases within a 50-50 domain division, we assess the model’s resilience by adjusting the total number of infectious individuals in the PDE contribution of the hybrid model to match the total population number. Simultaneously, we set it to zero in the ODE, and vice versa. Additionally, we introduce a Gaussian circular function that remains constant along the y-axis directly at the boundary $\Gamma$ at $x=1.0$. 
%
%
The results for each domain are shown in Figure \ref{fig:extreme_cases_rectangle}. Due to significant differences between the models, the infectious densities of the two models cannot be effectively visualized for comparison. Simulation plots of the fractional infectious density for the hybrid model were created for an extreme case (see Fig. \ref{simulation_extreme_case_rectangle}), serving the purpose of visualization. These simulation plots provide no additional value, hence, the exclusion of spatial visualization for other extreme cases.  
It is evident that the outcomes for the ODE region appear visually unaffected by the location of the population in the PDE domain. However, it is noticeable that there is a flow of infections into the ODE, indicating that the boundary conditions influence the ODE. This is evident as infection numbers increase in the right plot from zero to values as high as 8,256.5 (yellow) and 8,246.43 (green). %
Due to their simplicity, ODE models may exhibit a higher sensitivity to specific changes in parameters compared to PDE models. 
Starting at time $t=10$, a new time interval begins, introducing updated parameters. Visually, this adjustment has no apparent impact on the infection dynamics within the PDE, in contrast to the ODE, which experiences a second minor wave. This wave rapidly diminishes with the adoption of the new parameters at time $t=29$, leading to a noticeable decline in infection numbers within the PDE region as well. 


When examining the right plot of Figure \ref{fig:extreme_cases_rectangle}, where the initial number of infectious individuals is zero, the PDE with horizontal population density (red) exhibits a significantly lower peak in the infection curve compared to the PDE with vertical population density (blue). Due to the boundary conditions between the two regions, there is an infection flow from the ODE region to the PDE region (see e.g. Fig. \ref{simulation_extreme_case_rectangle}). This is more pronounced because the population density is increased precisely at the boundary and diminishes, even approaching zero, in the rest of the region. Consequently, the infection flow has a stronger impact. In the case of horizontal population density, either the population density is very low where the fractional infectious density is higher, or the fractional infectious density is low where the population density is high. As a result, the infection flow has a lesser effect. 
\\ \\
A similar principle applies to the PDE results in the left plot. When examining the left plot (and the corresponding data), where the initial number of infectious individuals corresponds to the population size, the PDE with horizontal population density (green) has fewer total infectious than the ODE (red, blue) for all time steps except the last few. 
Due to the boundary conditions and the absence of infectious in the ODE, the fractional infectious density of the PDE is reduced at the boundary, once again exerting a stronger influence on the vertical population density than the horizontal. Meanwhile, the infection numbers in the ODE region increase, and this, in turn, has a more pronounced impact on the rising infection numbers in the PDE for the vertical population density (yellow) than the horizontal (green).
\\ \\
Understanding the infection dynamics in a hybrid PDE-ODE model with Dirichlet boundary conditions is crucial for unraveling the interplay between different regions. The chosen boundary conditions significantly influence the infection dynamics within this hybrid framework. For instance, if the fractional infectious density in the ODE region is higher than at the boundary of the PDE region, the simulation 
exhibits a migration of infected individuals from the interior of the PDE region toward the boundary to fulfill the Dirichlet conditions. This 
leads to an accumulation of infected individuals at the boundary. %
Conversely, in cases where the fractional infectious density in the ODE region is lower than at the boundary of the PDE region, the system 
demonstrates a diffusion of infected individuals into the interior of the PDE region. In such scenarios, the infection dynamics would be predominantly influenced by the PDE model, with minimal impact observed in the ODE model.

To amplify the influence of the Dirichlet boundary conditions and observe heightened effects in the ODE region, adjusting the penalty term could be employed. By assigning a higher weight to the boundary conditions through an elevated penalty term, the model emphasizes the importance of meeting the specified conditions at the boundary, potentially leading to more pronounced effects in the ODE solution.

\subsection{Real-World Examples} \label{sec:Real-World Examples}

In the following section, we evaluate the practical applicability of our hybrid PDE-ODE model through its application to actual epidemiological scenarios. This examination not only validates the model against empirical data, but also highlights its potential for real-world utility in public health management. We illustrate the model's effectiveness by applying it to two distinct geographical settings: Lombardy, Italy, and Berlin, Germany. These locations were chosen due to their varied population densities and distinct public health challenges during the COVID-19 pandemic, offering a robust framework for testing and demonstrating the model's capabilities. For each region, we carefully calibrated the model using available epidemiological data, adapting our simulations to reflect the specific disease dynamics observed in each area. The results from these case studies provide insights into the model's ability to simulate and predict the complexities of infectious disease dynamics in the context of diffusive processes.

\subsubsection{Lombardy, Italy}

\begin{figure}[h]
        \centering
        \includegraphics[width=0.5\textwidth]{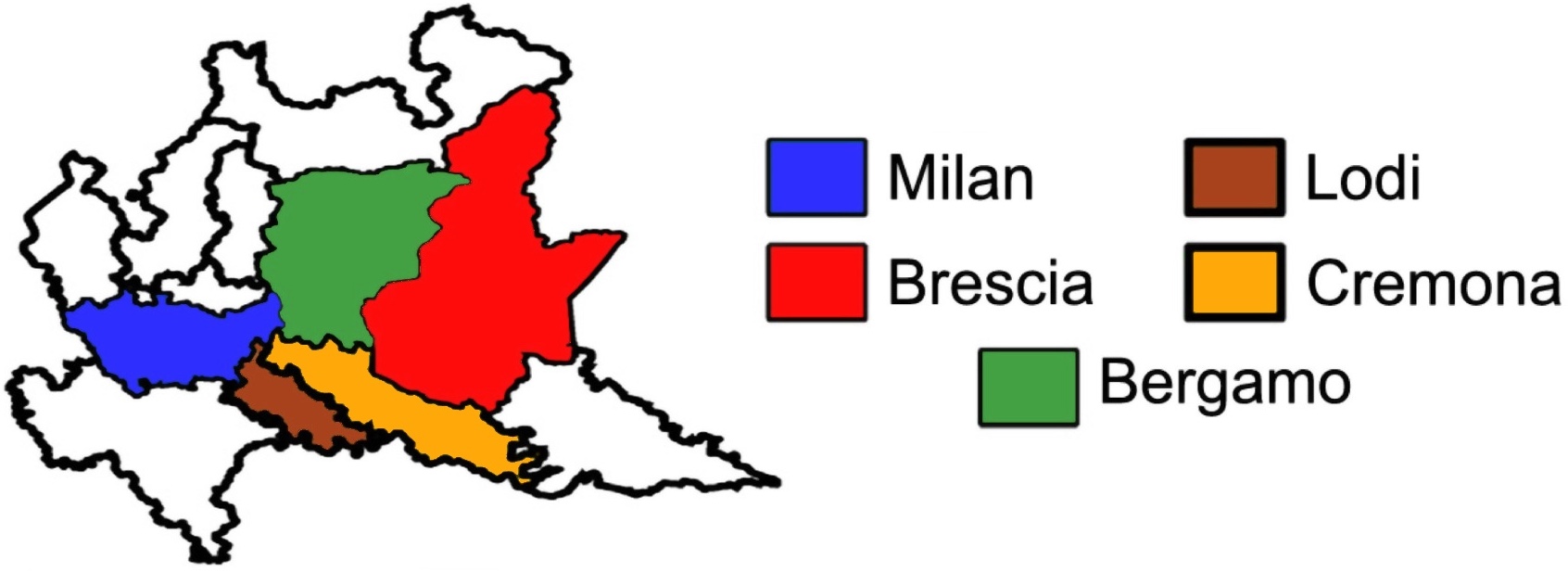}
        \caption{Main provinces affected by the pandemic in Lombardy: Lodi, Milan, Cremona, Bergamo, Brescia.}
        \label{fig:lombardy_provinces}
\end{figure}

In the first example, we explore the application of our hybrid PDE-ODE model to the Italian region of Lombardy (see Fig. \ref{fig:lombardy_provinces}). We follow the approach by Viguerie et al. \cite{Quelle5} and focus on the five largest provinces, which were mainly affected by the pandemic. We compare the performance of our hybrid model with a full-PDE model. Specifically, in the hybrid model, the province of Milan, which is relatively small but has a high population, is simulated using an ODE due to its high population density, facilitating faster and more uniform population exchanges, making it a more suitable candidate for ODE simulation compared to the other provinces. The aim of this example is to demonstrate the predictive power of our model with respect to the spread of Covid-19. The geographical data, i.e., (1) latitude and longitude of provinces, (2) boundary points of Lombardy, and (3) boundary points of Milan were taken from public repositories \cite{Quelle8,GithubMilan}\footnote{Note that in the case of longitude, Lombardy is close to the equator, i.e., we can simply scale the location data by $111.3 \; \nicefrac{\text{degrees}}{\text{km}}$ to convert \cite{Wiley2020Navigation}. The same holds for the latitude, since each degree of latitude is approximately $111.3$ km apart. This is not necessary in the calculations and will ultimately occur before the dynamics are simulated.}. We calibrated our model using documented case data from Lombardy \cite{Quelle8}, calculating a 7-day average from these cases as target for parameter identification.

To simulate the dynamics within this region, the mesh setup consists of 63,017 triangles, 99,084 edges and 36,068 vertices. 
When adapted for the coupled model, this results in a mesh configuration of 31,953 triangles, 52,912 edges, and 20,960 vertices. For the initial conditions, we use circular Gaussian functions with a variance of $0.05$.  
The Lombardy specific model parameters are listed in Table \ref{table:parameters_Lombardy} and are adjusted across four distinct time intervals. Start- and end-times are taken from Viguerie et al. \cite{Quelle5}. Each interval corresponds to different phases of pandemic response, reflecting changes in public health policies and restrictions. These parameters are finely tuned for each period and remain constant over the time interval. During the initial time interval, no restrictions are in effect. However, as each subsequent time interval begins, restrictions are introduced, and changes in parameter values show how policy changes affect disease dynamics. 
%
%
%
%
%
%
\begin{table}[h!] 
\centering %
\begin{tabular*}{\textwidth}{
@{\extracolsep{\fill}} *{5}{c}} 
\multicolumn{1}{c}{parameter}    &
\multicolumn{4}{c}{value}    \\
\cmidrule{1-1}
\cmidrule{2-5}
           & February 27th, 2020 & March 9th, 2020 & March 22nd, 2020 & March 28th, 2020 \\
  \hline
  \midrule %
  $\sigma$              & $2.6676e-02$  
  & - & - & $1.8747e-01$ \\
  $\phi_e$              & $0.0000e+00$ 
  & - & - & $3.2655e-11$ \\
  $\phi_i$              & $2.3310e-01$  
  & - & - & $2.1907e+00$ \\
  $\beta_i$, $\beta_e$  & $4.4202e-01$ & $2.0588e-01$ & $6.0352e-08$ & $2.2168e-01$ \\
  $D$  & $\tfrac{0.0435}{111.3^2} \cdot 1.0e+02 $ 
  & $\tfrac{0.0198}{111.3^2}  \cdot 1.0e+2$ & $\tfrac{0.0090}{111.3^2} \cdot 1.0e+2$ & $\tfrac{0.0075}{111.3^2} \cdot 1.0e+2$ \\
  $A$                   & $3.0e+07$ & - & - & - \\
  $n_0$                 & $4.5e+07$ & - & - & - \\
\bottomrule 
\end{tabular*}  %
\caption{Parameters of PDE model (\ref{PDE-model}) and of the hybrid model (\ref{Eq:hybrid_model}) for Lombardy. 
The symbol “-” means that the parameter has not changed over time. } 
\label{table:parameters_Lombardy}
\end{table}%
\noindent
\begin{figure}[h]
    \begin{minipage}{0.5\textwidth}
        \includegraphics[width=\textwidth]{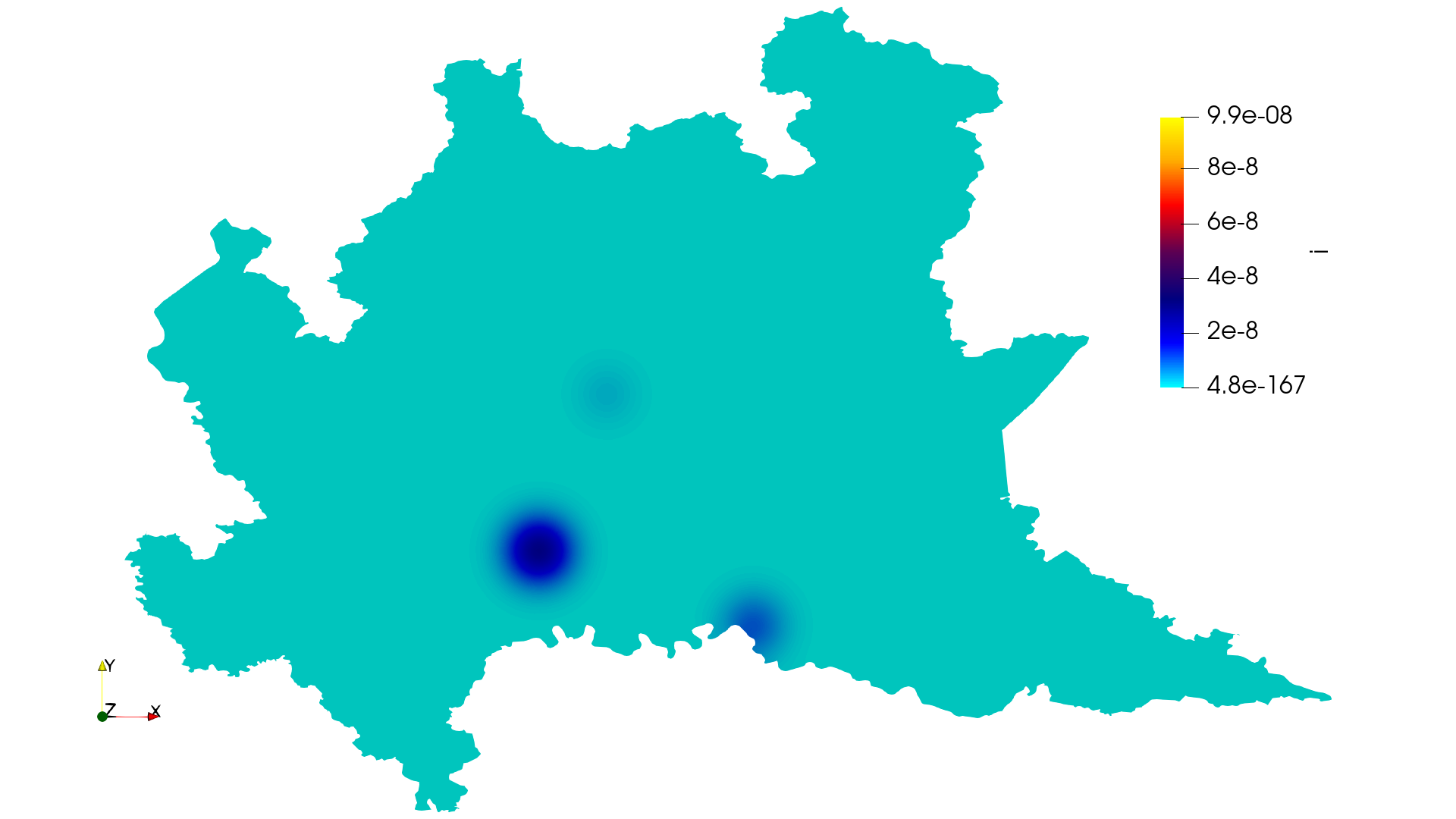}
    \end{minipage}
    \begin{minipage}{0.49\textwidth}
        \includegraphics[width=\textwidth]{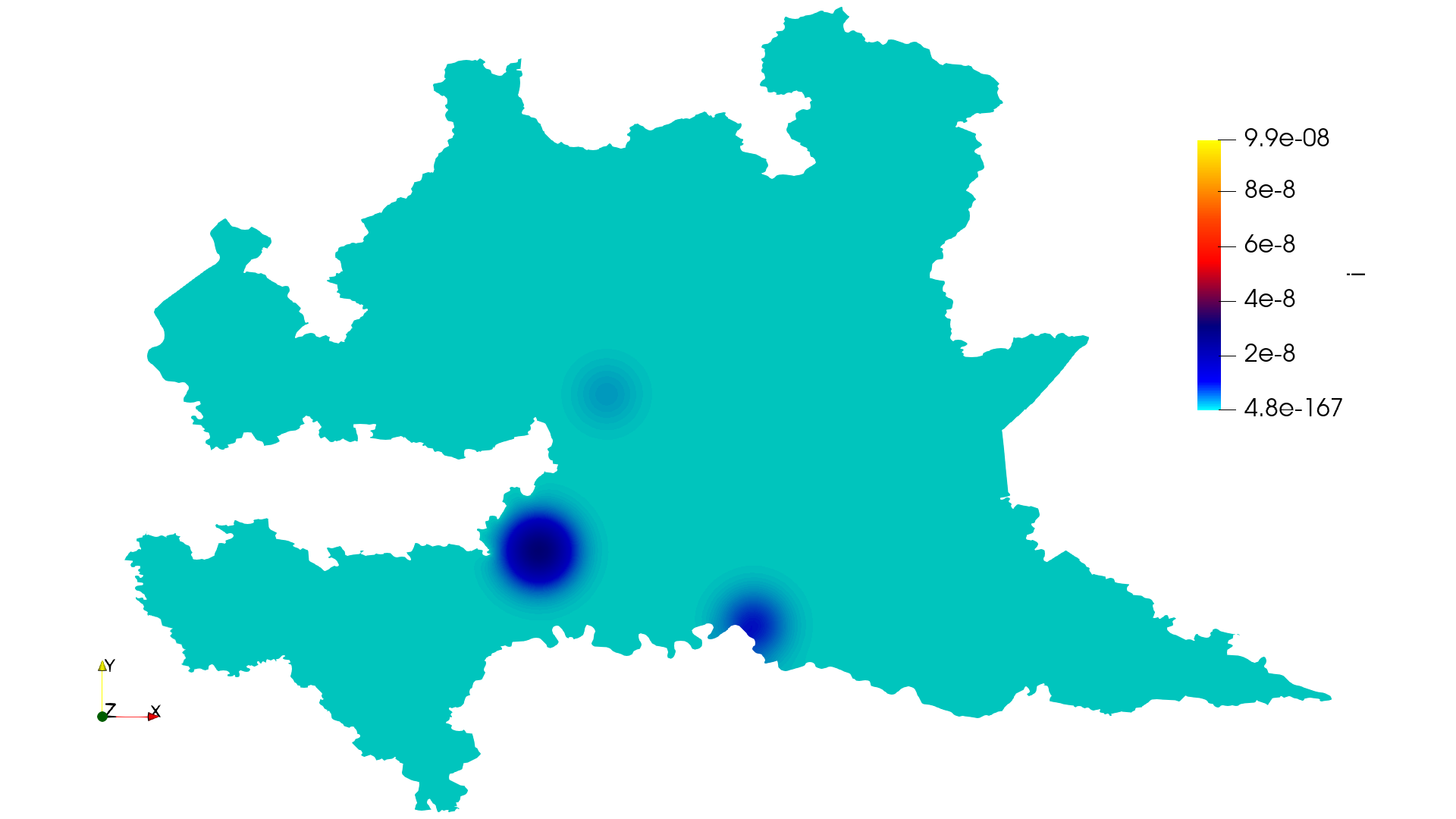}
    \end{minipage}
    \\
    \begin{minipage}{0.5\textwidth}
         \includegraphics[width=\textwidth]{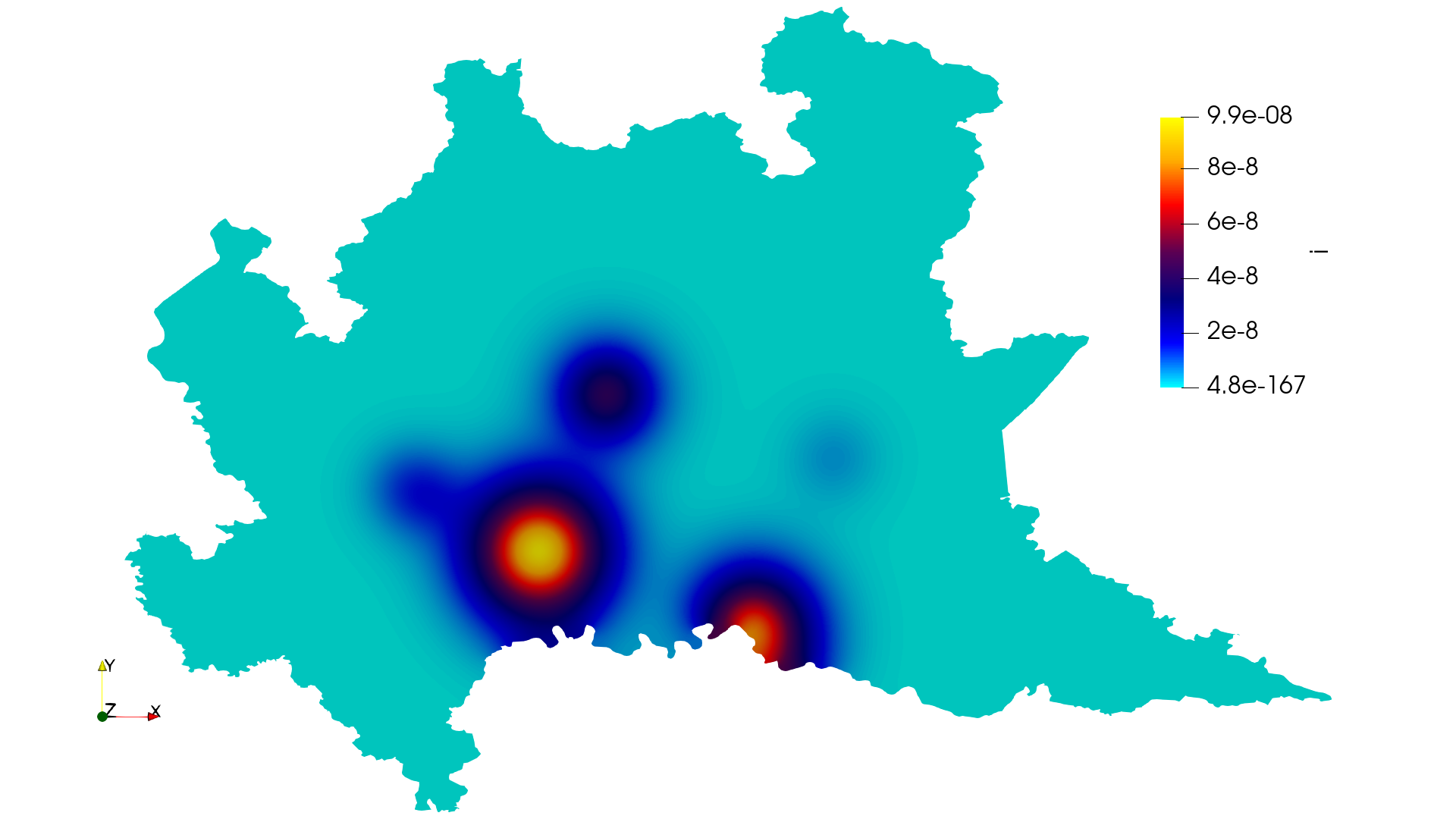}
    \end{minipage}
    \begin{minipage}{0.49\textwidth}
         \includegraphics[width=\textwidth]{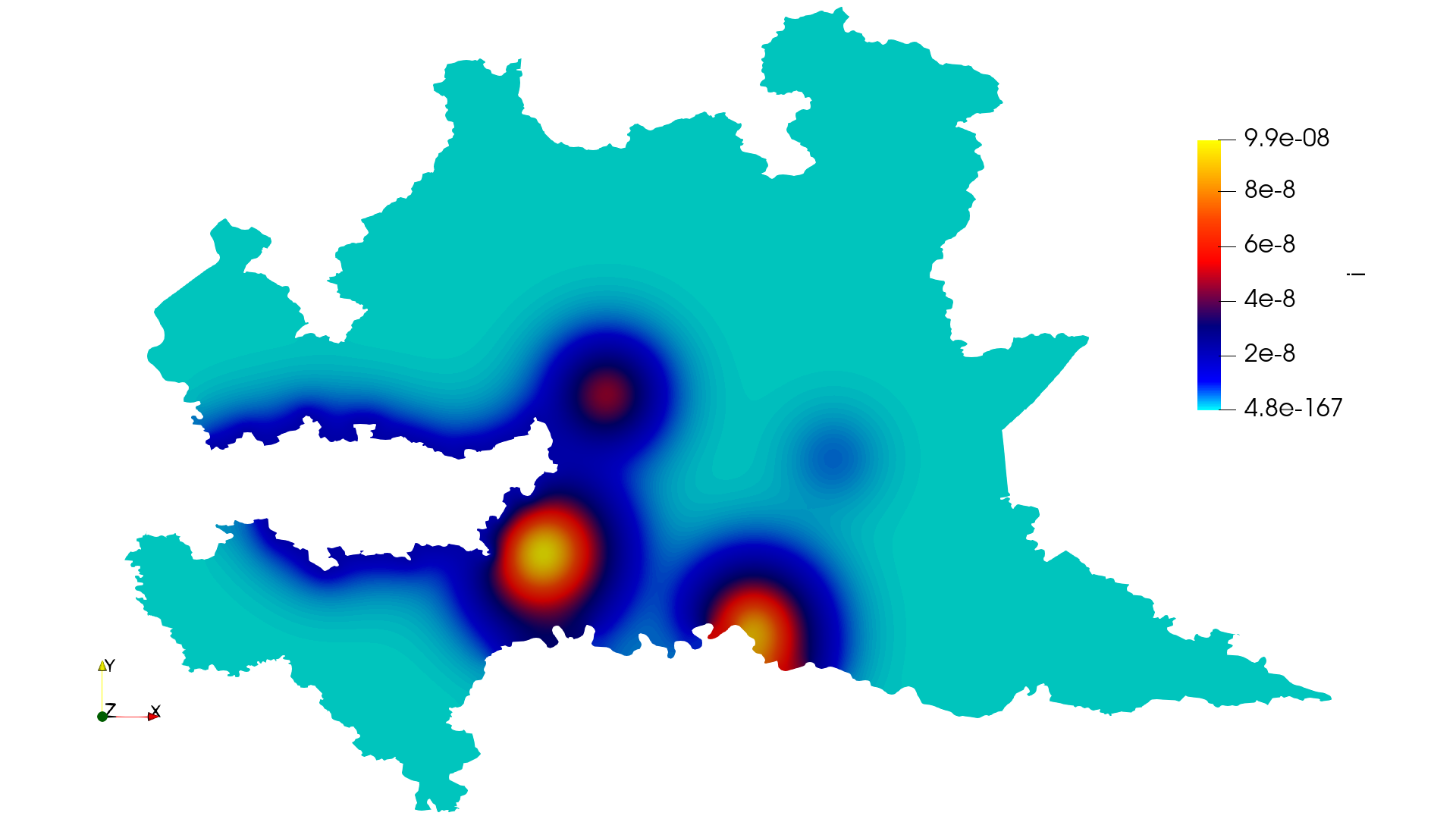}
    \end{minipage}
    \\
    \begin{minipage}{0.5\textwidth}
        \includegraphics[width=\textwidth]{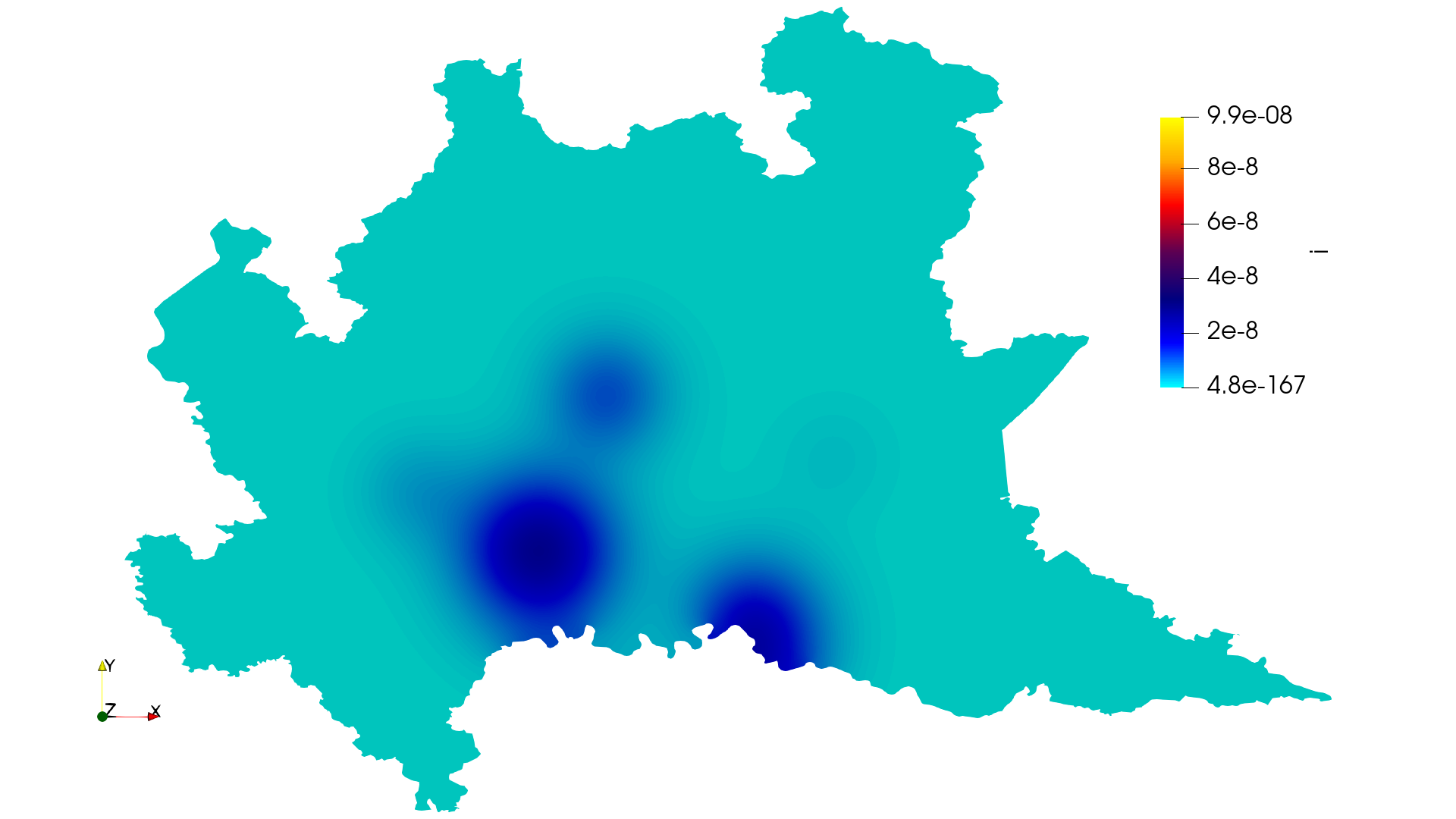}
    \end{minipage}
    \begin{minipage}{0.49\textwidth}
        \includegraphics[width=\textwidth]{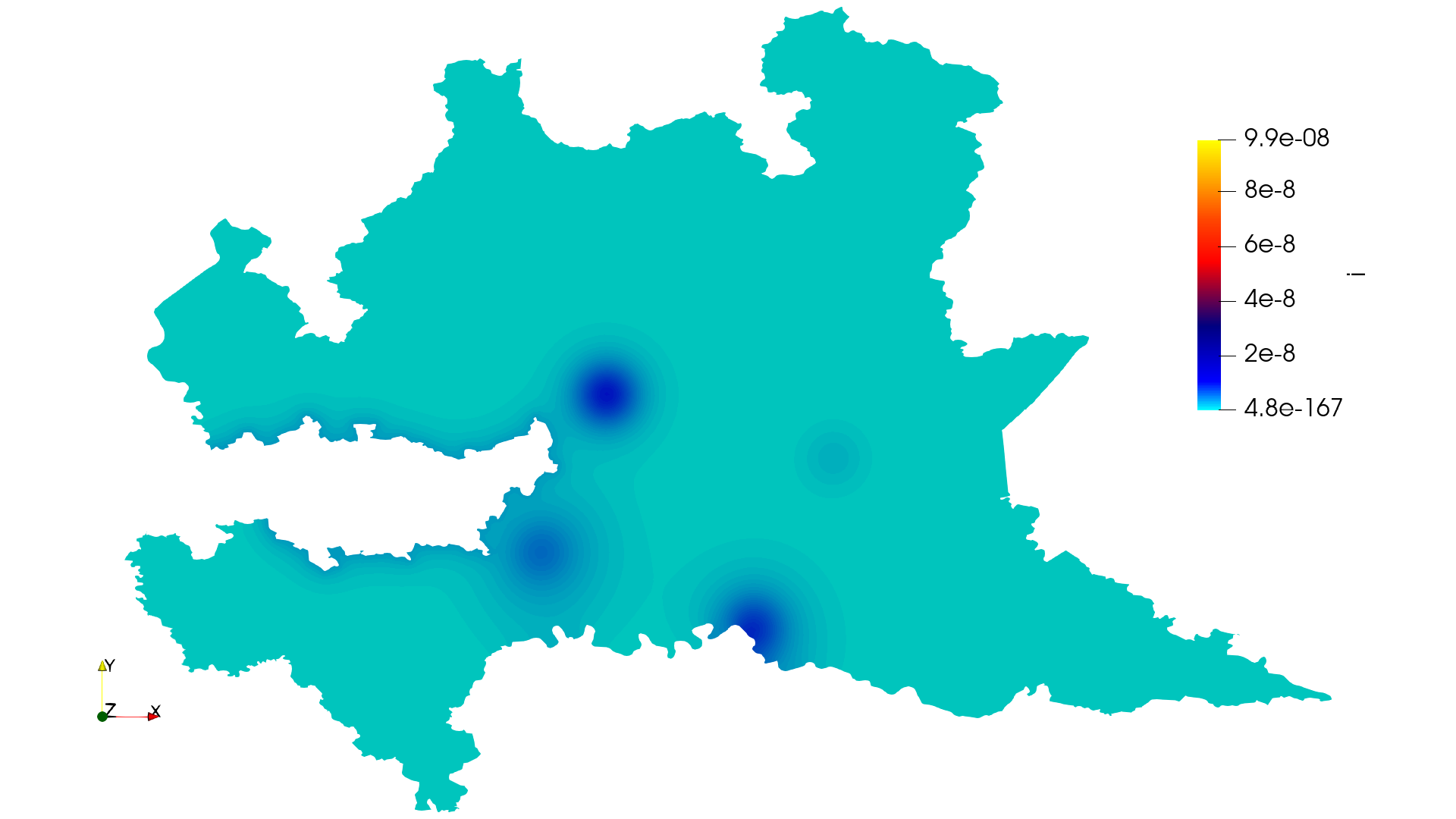}
    \end{minipage}
\caption{Fractional infectious density in km$^2$ of PDE model (\ref{PDE-model}) (left) and PDE contribution of hybrid model (\ref{Eq:hybrid_model}) (right) in Lombardy for times $t \in \{1,28,59\}$. The boundary conditions between the PDE and ODE domains ensure, among other things, a uniform flow of infectious individuals originating from Milan in all directions.} 
\label{fig:Lombardy_fractional_infectious_density_both_models}
\end{figure} 
\noindent
\begin{figure}[h!]
    \begin{minipage}{\textwidth}
        \centering
        \includegraphics[width=0.8\textwidth]{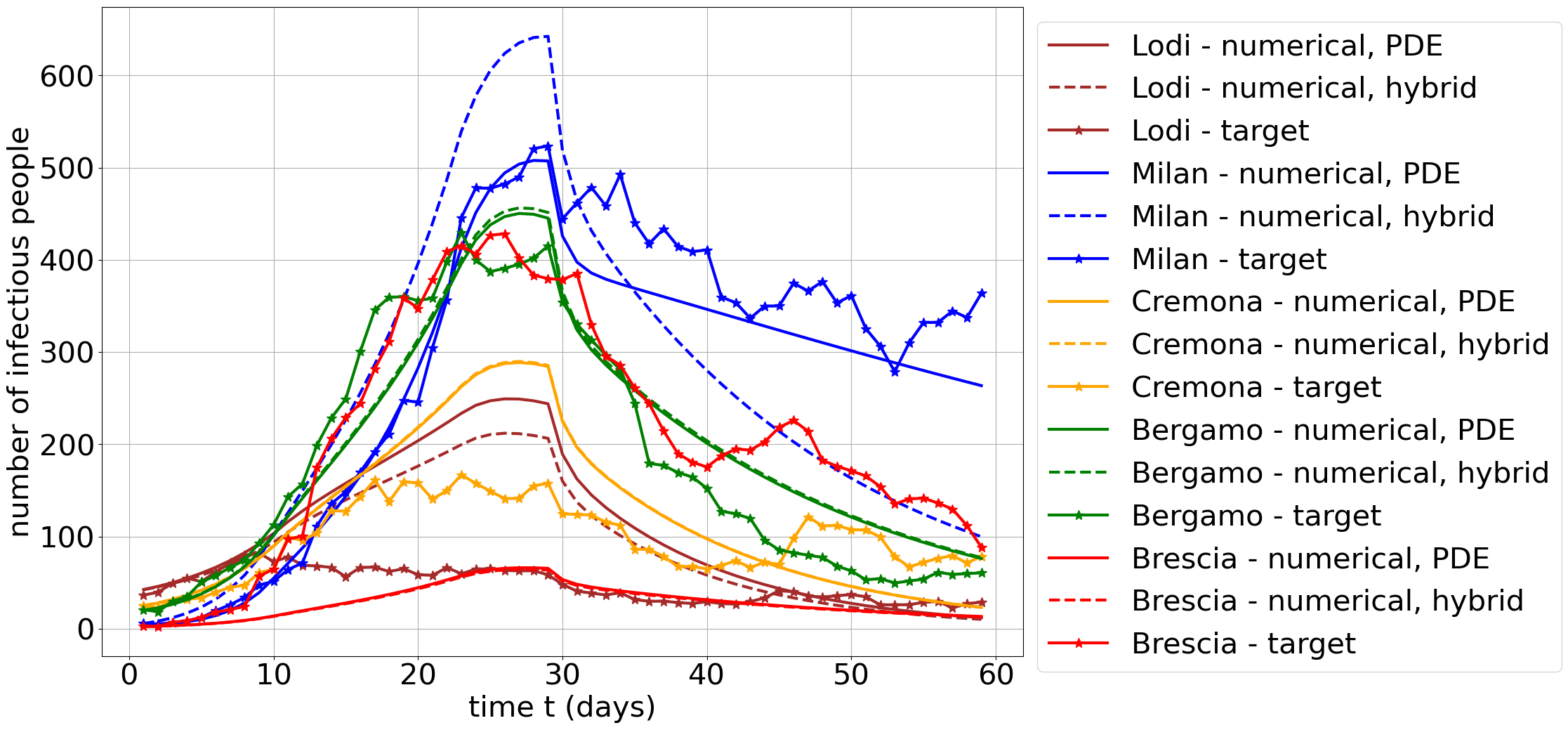}
        \caption{Number of infectious people of PDE (\ref{PDE-model}) and hybrid models (\ref{Eq:hybrid_model}) in each province (simulated and real 7-day average data). In Cremona, Bergamo, and Brescia, both models showed similar trends, but the hybrid model overestimated infections in Milan and underestimated them in Lodi.} 
    \label{fig:number_of_infectious_pde_and_hybrid_lombardy}
    \end{minipage}
\end{figure} 
\noindent
\\
The dynamics are simulated and shown in Fig. \ref{fig:Lombardy_fractional_infectious_density_both_models}, and the total number of infectious individuals are depicted in Fig. \ref{fig:number_of_infectious_pde_and_hybrid_lombardy}.
A notable feature of our simulations is the incorporation of inter-provincial interactions, particularly the bidirectional transmission of infections between Lodi and Milan, facilitated by specified boundary conditions. This dynamic interaction results in a significant impact on infection trends in both provinces. \\ \\
Our findings revealed a consistent pattern across both modeling approaches in provinces Cremona, Bergamo, and Brescia (see Fig. \ref{fig:number_of_infectious_pde_and_hybrid_lombardy}). The simulation results for these provinces showed similar trends, indicating that both models provide consistent and reliable predictions, regardless of the chosen modeling approach, for these provinces located farther from the boundary of the ODE-modeled area.
However, a noteworthy deviation emerged in the hybrid model's predictions for Milan and its neighboring province, Lodi (see Fig. \ref{fig:number_of_infectious_pde_and_hybrid_lombardy}). The hybrid model tended to overestimate the infectious cases in Milan. Simultaneously, it underestimated the infection rates in the adjacent province of Lodi.
The boundary conditions between the PDE and ODE domains ensure, among other things, a uniform flow of infectious individuals originating from Milan in all directions, a feature that is clearly visible in Fig. \ref{fig:Lombardy_fractional_infectious_density_both_models}. This leads to a slight increase in infectious numbers in Bergamo. Simultaneously, the infectious flow from Lodi towards Milan is restrained due to the boundary conditions. Instead, infectious cases tend to migrate more prominently towards the neighboring province of Cremona. \\ \\
The parameter fitting for the full-PDE model naturally yields results closer to the target data than those obtained from the hybrid model. An alternative approach involves directly fitting the parameters of the hybrid model.
\newpage %

\subsubsection{Berlin, Germany} 
For the second example, we apply our hybrid PDE-ODE model to the urban context of Berlin, Germany. Berlin presents an interesting case with its diverse and dense urban structure, which offers a challenging environment for epidemiological modeling due to varied population distribution and mobility patterns.

The geographical data, i.e., latitude and longitude of Berlin, are publicly available\footnote{\url{http://opendatalab.de/projects/geojson-utilities/}}. Utilizing this data, as before, we created a mesh for the full-PDE model, consisting of 73,294 triangles, 110,471 edges and 37,178 vertices. For our hybrid model, we divided the domain vertically. This resulting reduced domain consists of 36,633 triangles, 55,249 edges and 18,617 vertices. 

However, opposed to the Lombardy study, this time we used data generated by an agent-based model (ABM) specifically calibrated for Berlin (see \cite{Muller2021ABM} for details). This ABM effectively simulates individual-level behaviors that influence disease transmission dynamics. The model is based entirely on mobility data, derived from mobile phone data, and incorporates a 25\% sample\footnote{\url{https://doi.org/10.14279/depositonce-11495.2}} of the Berlin population. This allows capturing typical movements of persons as well as their activity reductions during the various public health restrictions. While the foundational ABM parameters were sourced from existing literature, the infection rate parameter was specifically calibrated to align with the observed data during the COVID-19 pandemic. 

While Lombardy's model focused on the dynamics across larger provincial areas, the Berlin model incorporates the granularity of individual agents, representing residents with specific demographic and health status characteristics. This approach allows us to capture detailed spatial information while maintaining a constant population density over time, providing a comprehensive understanding of infection dynamics in urban environments, such as Berlin.

\begin{figure}[ht]
    \centering
    \includegraphics[width=0.49\textwidth]{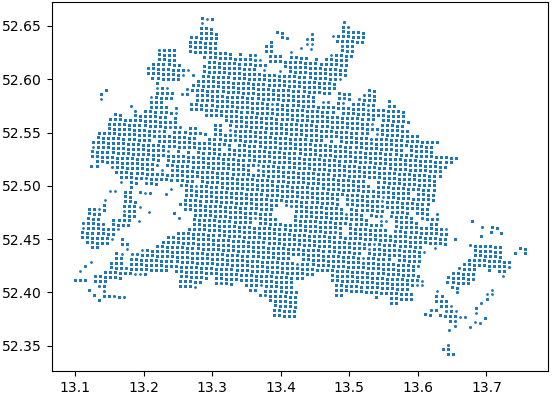}
    \includegraphics[width=0.49\textwidth]{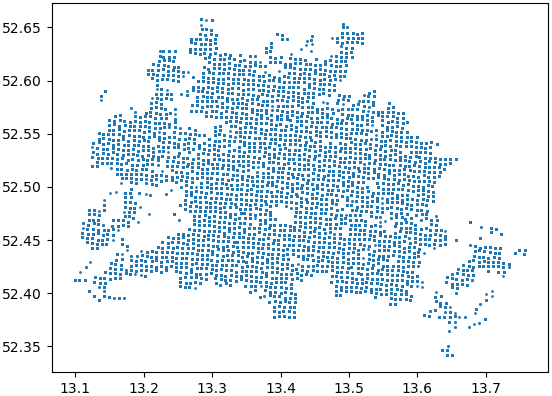}
    \caption{Original home coordinates of agents in Berlin, Germany (left) and adjusted locations to grid (right).}
    \label{fig:berlin_home_coordinates}
\end{figure}

Each agent is assigned home coordinates (see Fig. \ref{fig:berlin_home_coordinates}). We map these home coordinates to grid points based on the shortest distance. This enables us to determine the number of agents at each grid point with specific health statuses and how many agents are present there. Following this, we have smoothed the (sub)populations on the mesh. Subsequently, we calculate the density at each grid point. Initially, we determine the population density $n$, followed by the fractional densities $s$, $e$, $i$, $r$.
Since not all locations have residences in our area, resulting in zero population density in those areas, we fill the gaps with the smallest representable value to minimize their impact on the overall population. Filling the gaps ensures our condition $s+e+i+r=1$ holds. Additionally, this population fraction is labeled as \textit{removed} to ensure that this additional population does not affect the infection dynamics.
\noindent

Although Berlin is not situated close to the equator, for the purpose of our simulation, we have chosen a simplified approach by scaling location data using  $111.3 \; \nicefrac{\text{degrees}}{\text{km}}$. While this method may introduce some margin of error due to the deviation from the equator, it allows for a straightforward representation of spatial relationships in our model.
%
The restrictions were introduced in Berlin at different times and to a different extent than in Lombardy \cite{Bosa2022Response}. We chose March 15th, 2020, and March 23rd, 2020\footnote{\url{https://taz.de/Corona-Tagebuch-der-Hauptstadt/!5671055/}}, as the pivotal dates for the introduction of (additional) restrictions. We could not find details about the last restriction in Lombardy and, therefore, omitted it in the case of Berlin.
In Table \ref{table:parameters_Berlin}, you can find the parameters corresponding to different time intervals for Berlin.
\begin{table}[t!] 
\centering %
\begin{tabular*}{\textwidth}{
@{\extracolsep{\fill}} *{4}{c}} 
\multicolumn{1}{c}{parameter}    &
\multicolumn{3}{c}{value}    \\
\cmidrule{1-1}
\cmidrule{2-4}
           & February 25th, 2020 & March 15th, 2020 & March 23rd, 2020 \\
  \hline
  \midrule %
  $\sigma$              & $1.8513e+01$ & $1.8513e+01$ & $1.8513e+01$ \\
  $\phi_e$              & $2.9179e-03$ & $1.8561e-02$ & $8.3704e-03$ \\
  $\phi_i$              & $2.7734e-01$ & $3.1585e-01$ & $2.8986e-01$ \\
  $\beta_i$, $\beta_e$  & $2.9740e+00$ & $1.4135e+00$ & $1.4490e+00$ \\
  $D$  & $\tfrac{0.0435}{111.3^2} \cdot 1.0e+02 $ 
  & $\tfrac{0.0198}{111.3^2}  \cdot 1.0e+2$ & $\tfrac{0.0090}{111.3^2} \cdot 1.0e+2$ 
  \\
  $A$                   & $1.0e+09$ & - & - \\
  $n_0$                 & $1.5e+09$ & - & - \\
\bottomrule 
\end{tabular*}  %
\caption{Parameters of PDE model (\ref{PDE-model}) and of the hybrid model (\ref{Eq:hybrid_model}) for Berlin. 
The symbol “-” means that the parameter has not changed over time.} 
\label{table:parameters_Berlin}
\end{table}%
\noindent
\begin{figure}[t!]
    \begin{minipage}{0.49\textwidth}
        \includegraphics[width=\textwidth]{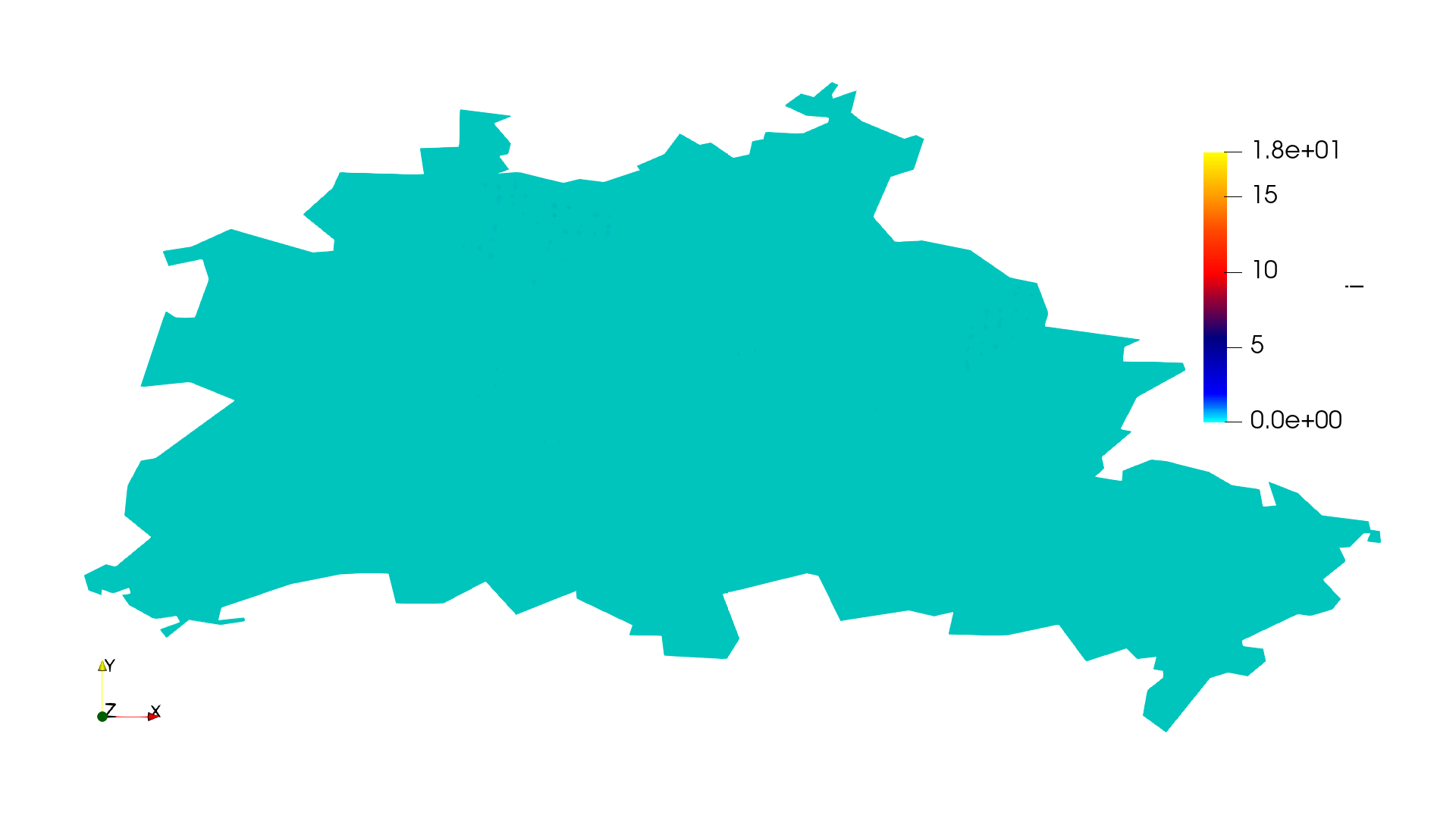} \\
        \includegraphics[width=\textwidth]{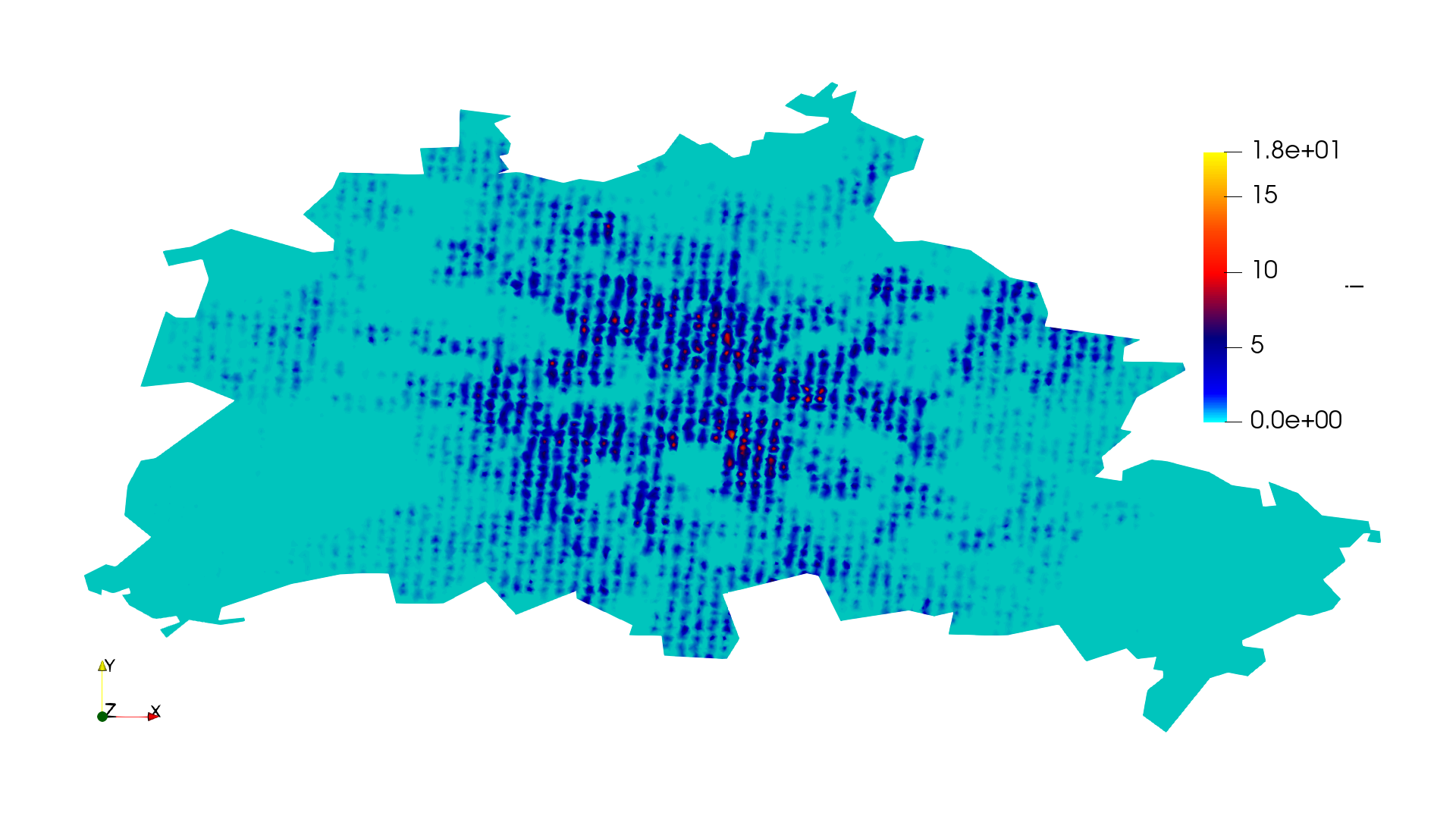}
    \end{minipage}
    \begin{minipage}{0.49\textwidth}
         \includegraphics[width=\textwidth]{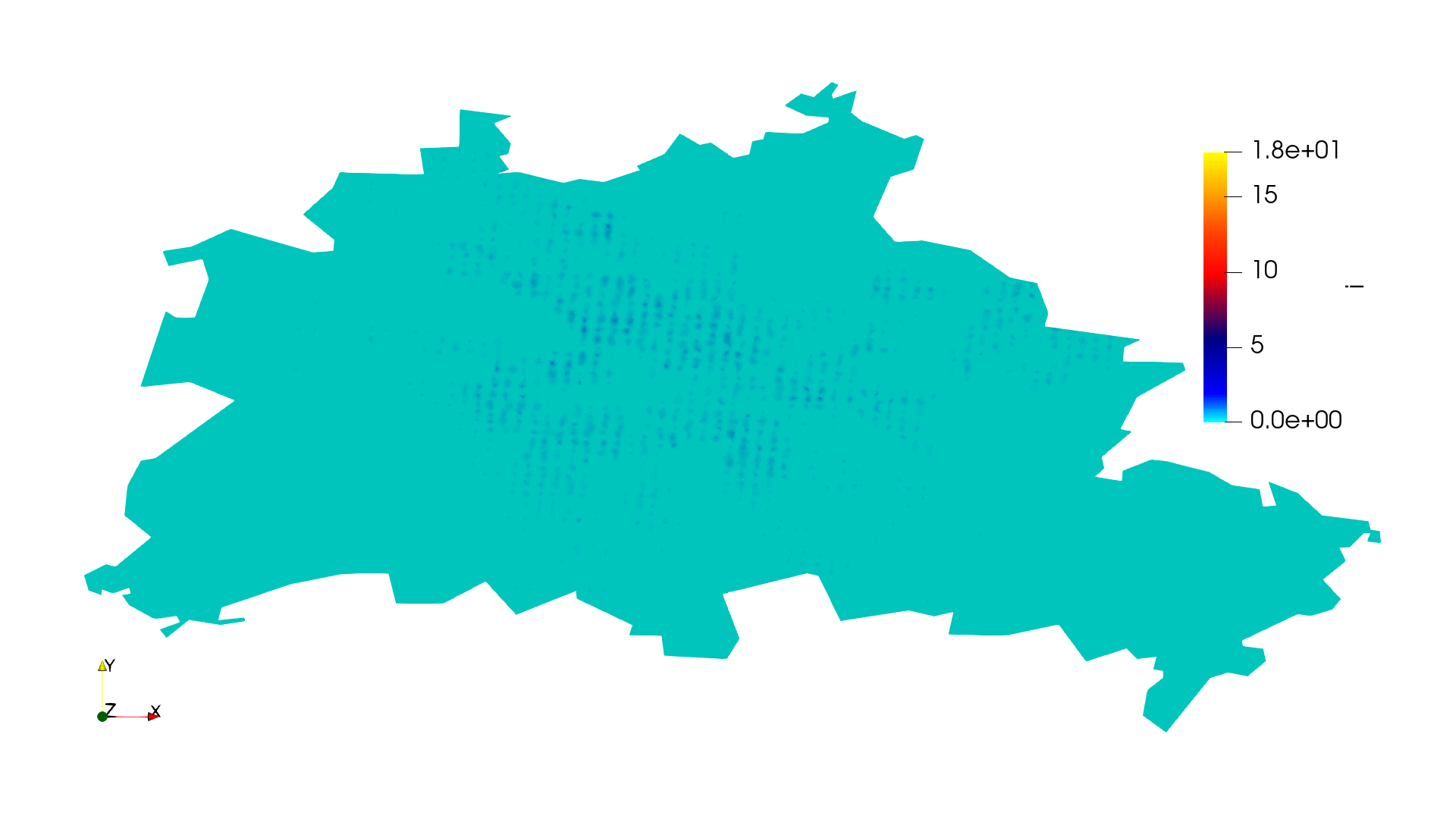} \\
         \includegraphics[width=\textwidth]{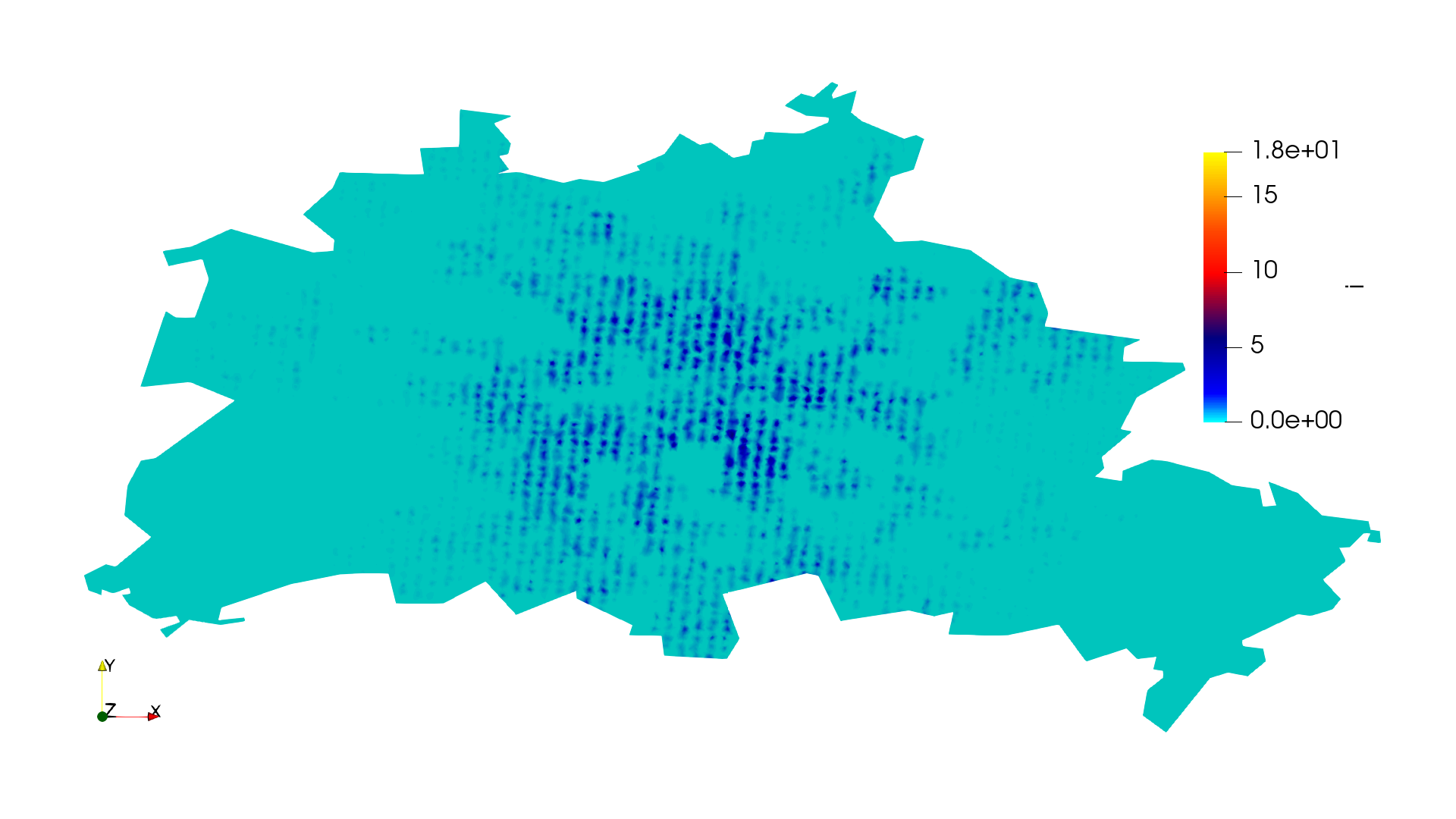}
    \end{minipage}
\caption{Infectious density of PDE model (\ref{PDE-model}) for times $t \in \{1,8,18,60\}$. Here, the spread can be observed particularly well at the level of home locations.} 
\label{fig:Berlin_infectious_density_PDE}
\end{figure} 
\noindent
%

\begin{figure}[h]
    \begin{minipage}{0.49\textwidth}
        \includegraphics[width=\textwidth]{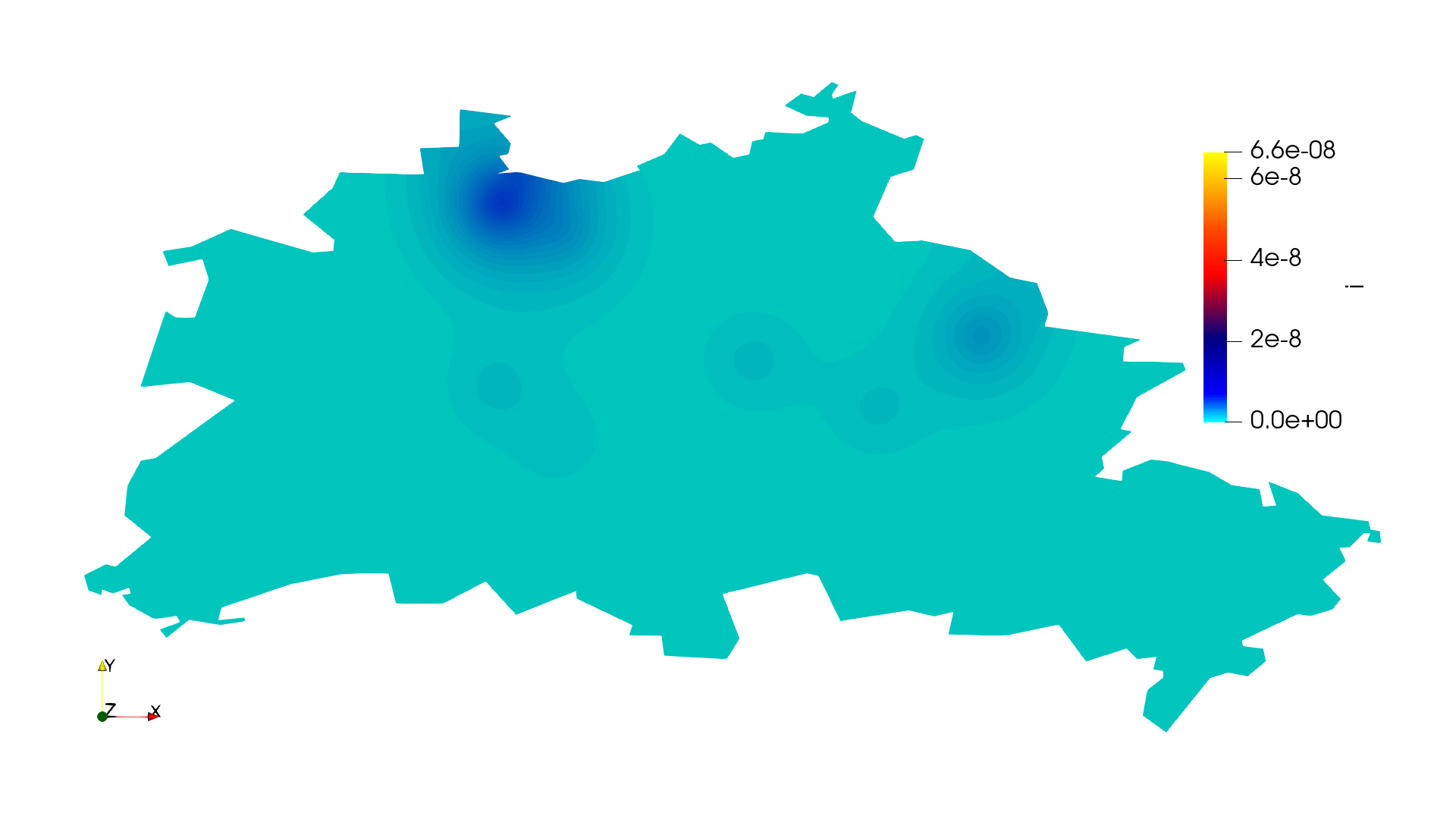} \\
        \includegraphics[width=\textwidth]{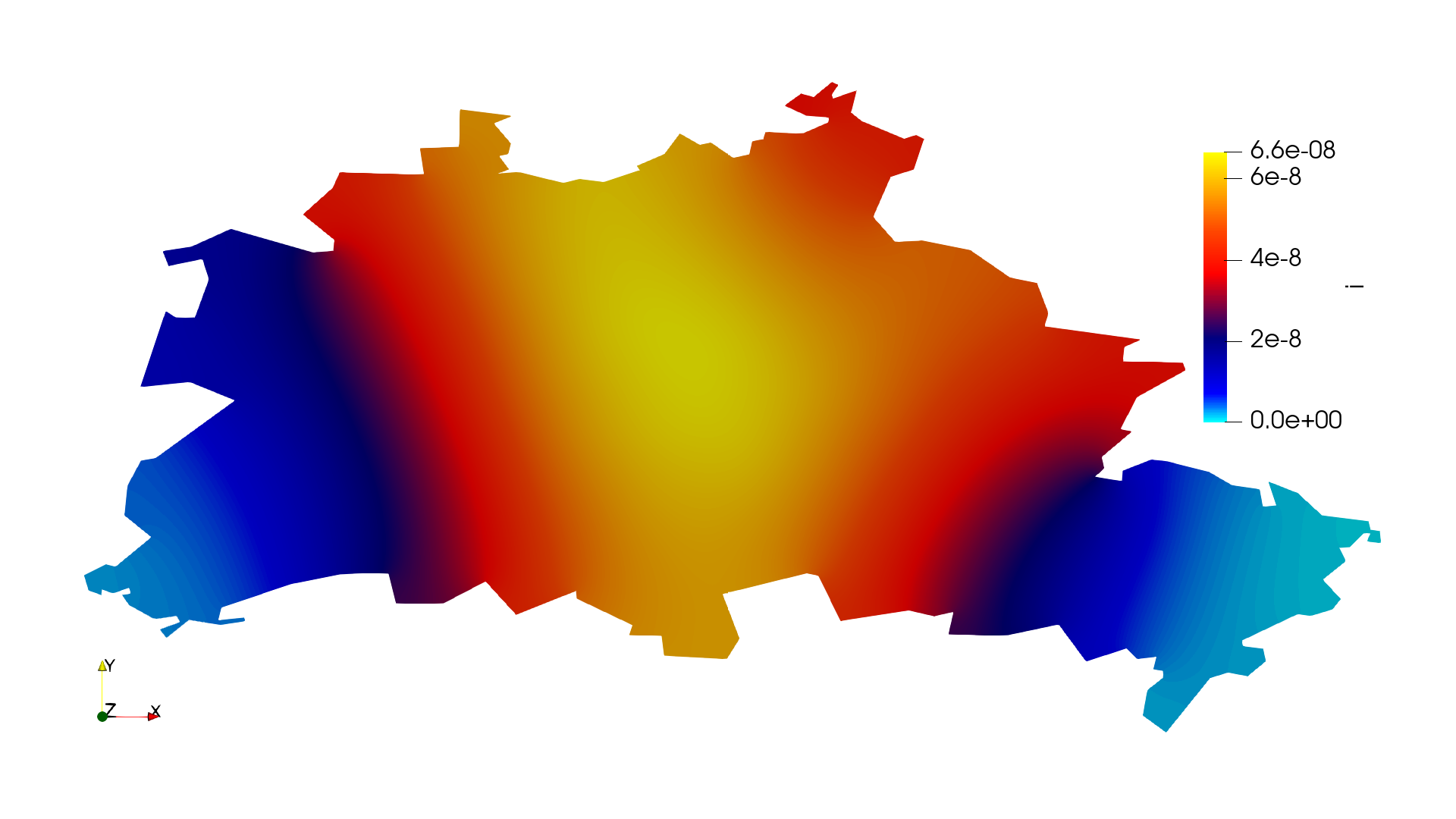}
    \end{minipage}
    \begin{minipage}{0.49\textwidth}
         \includegraphics[width=\textwidth]{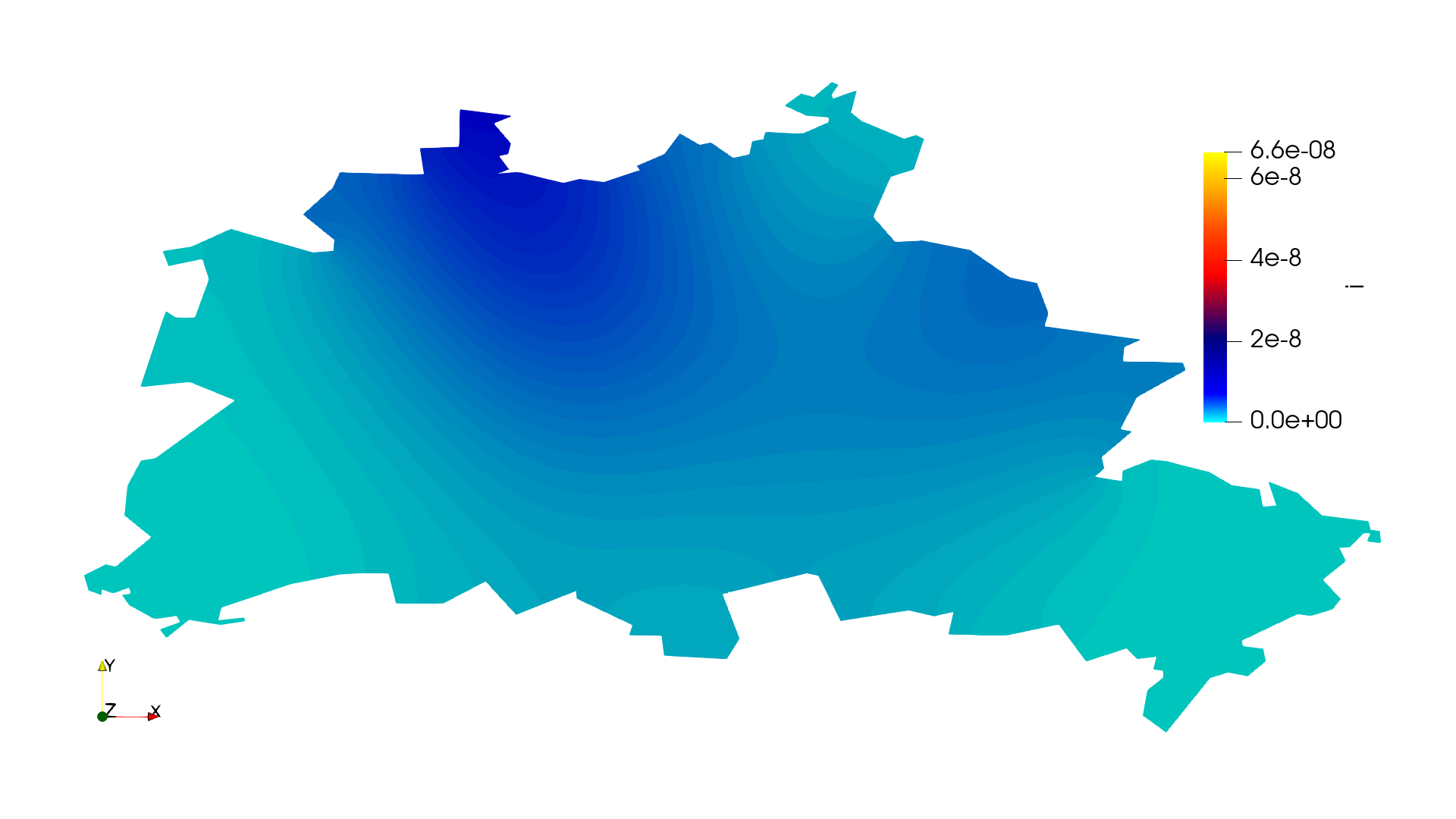} \\
         \includegraphics[width=\textwidth]{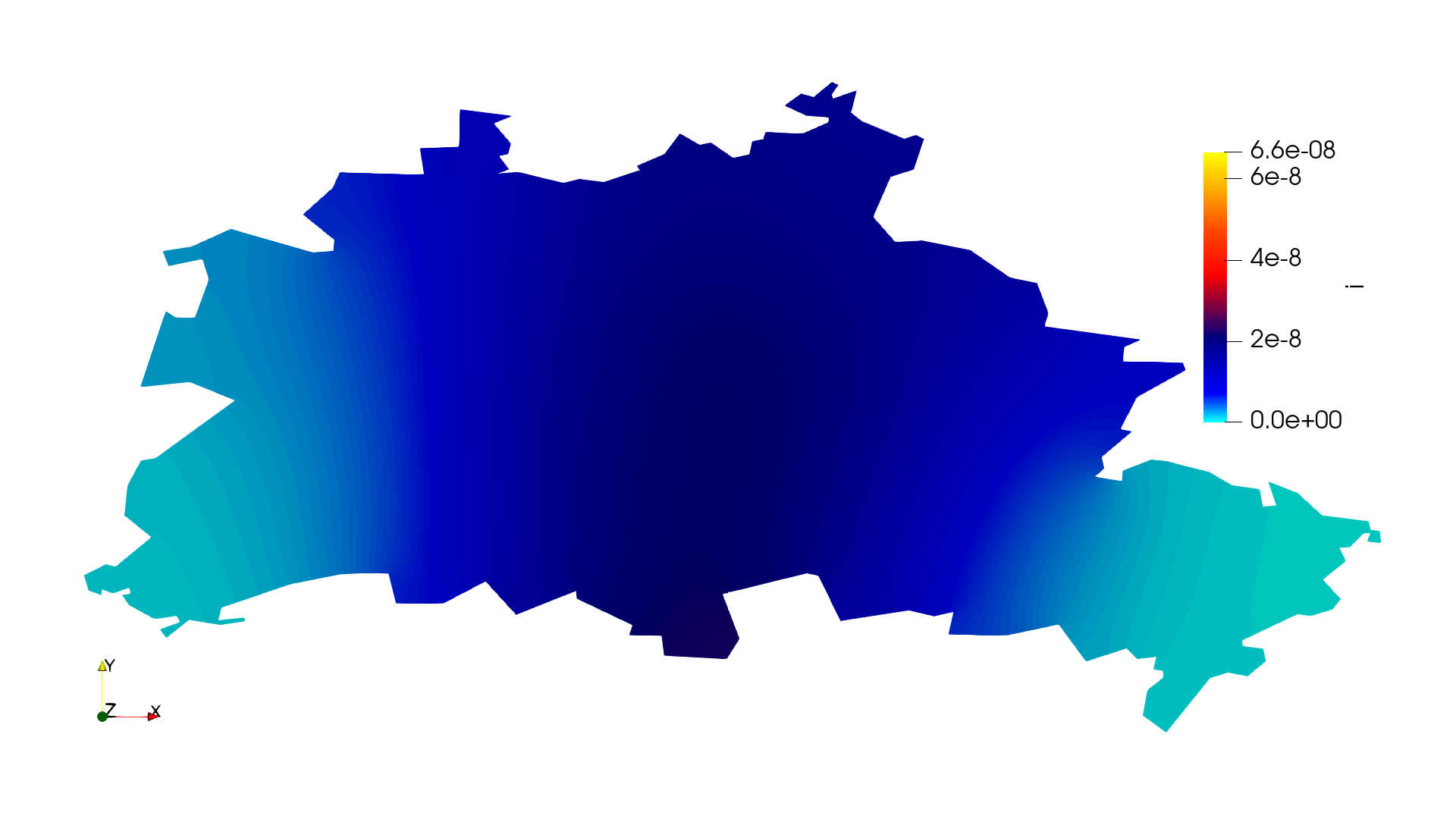}
    \end{minipage}
\caption{Fractional infectious density of PDE model (\ref{PDE-model}) for times $t \in \{1,8, 18,60\}$. Due to Berlin's significantly smaller area, the spread in Berlin appears faster than in Lombardy.} 
\label{fig:Berlin_fractional_infectious_density_PDE}
\end{figure}

\begin{figure}[t!]
\begin{minipage}{0.49\textwidth}
    \includegraphics[width=\textwidth]{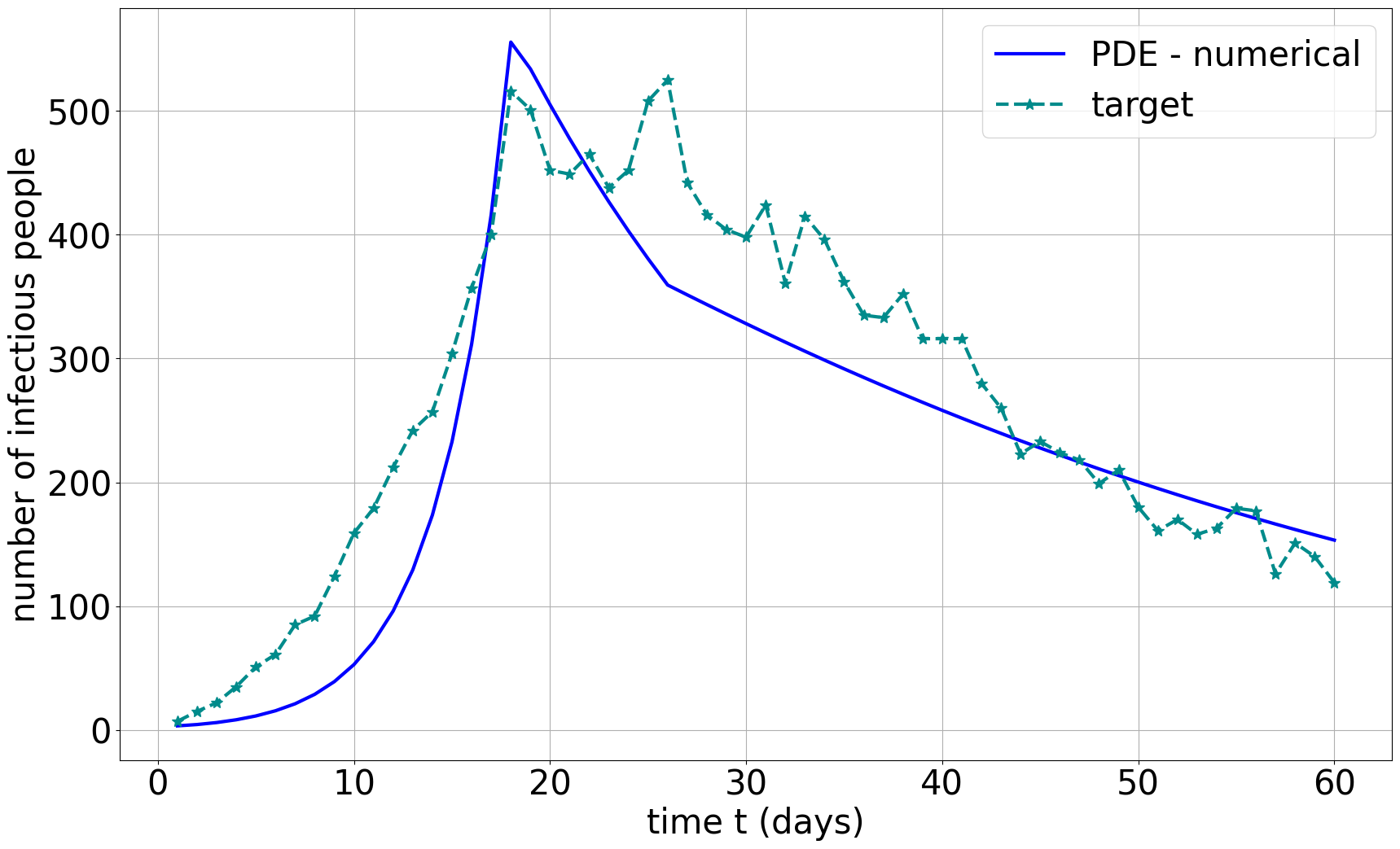}
\end{minipage}
\begin{minipage}{0.49\textwidth}
\includegraphics[width=\textwidth]{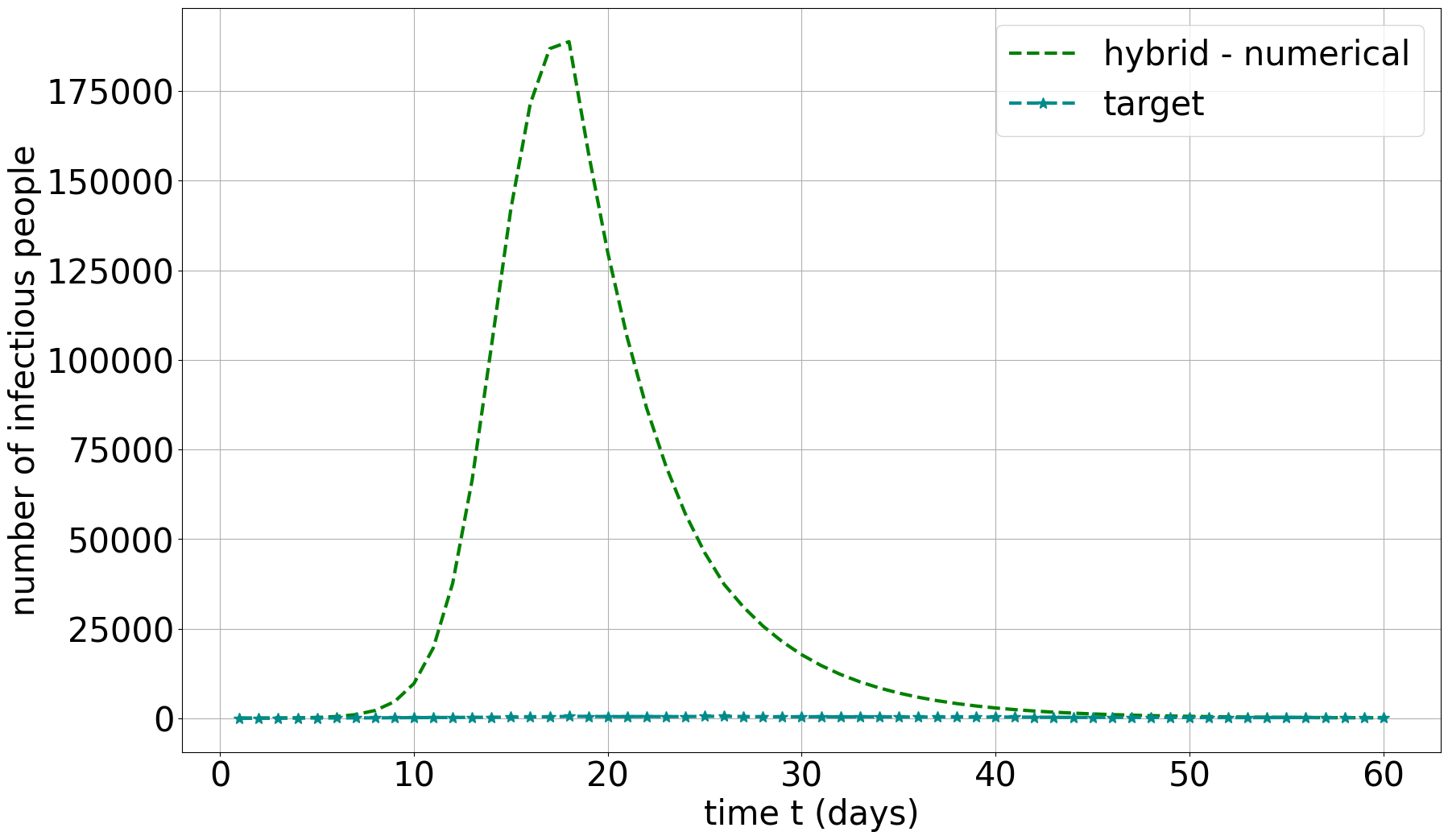}
\end{minipage}
    \caption{Number of infectious people of full-PDE model (\ref{PDE-model}) (left) and hybrid model (\ref{Eq:hybrid_model}) (right) in Berlin (simulated and ABM data). The primary reason for the varying number of infectious cases is likely the non-uniform distribution of the population within the PDE domain.} 
\label{fig:Berlin_number_of_infectious}
\end{figure}

\begin{figure}[ht]
    \begin{minipage}{0.49\textwidth}
        \includegraphics[width=\textwidth]{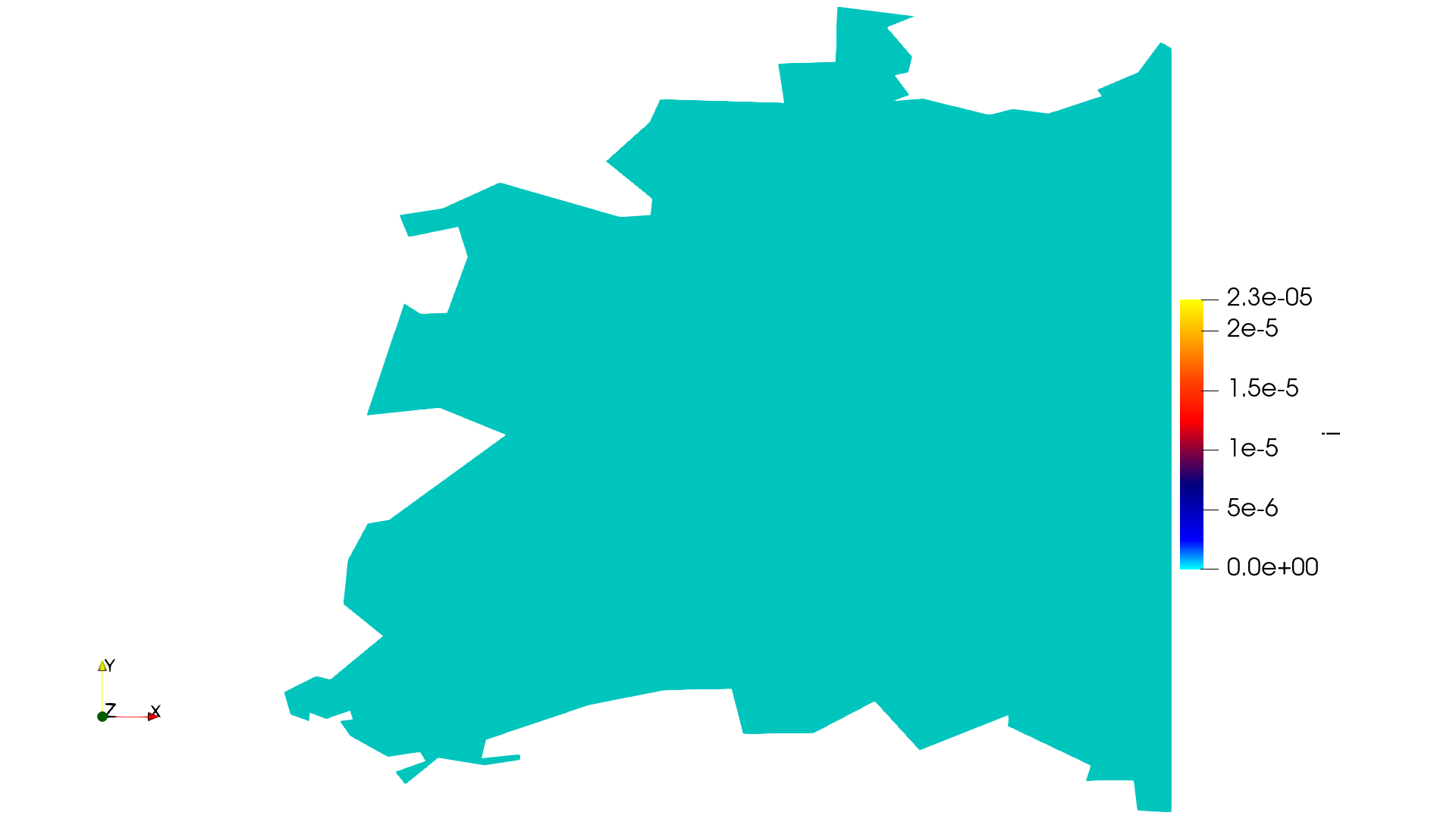} \\
        \includegraphics[width=\textwidth]{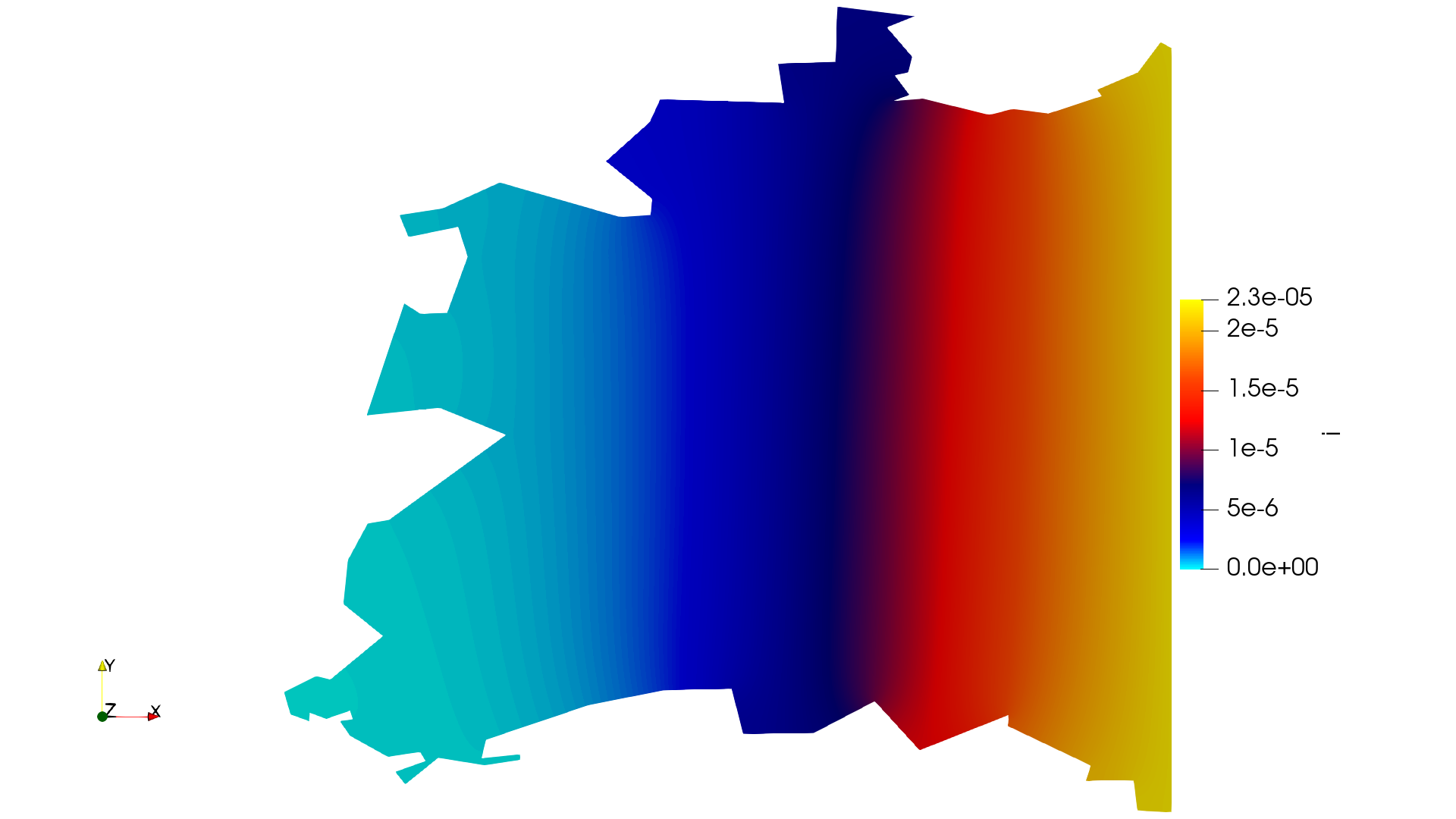}
    \end{minipage}
    \begin{minipage}{0.49\textwidth}
         \includegraphics[width=\textwidth]{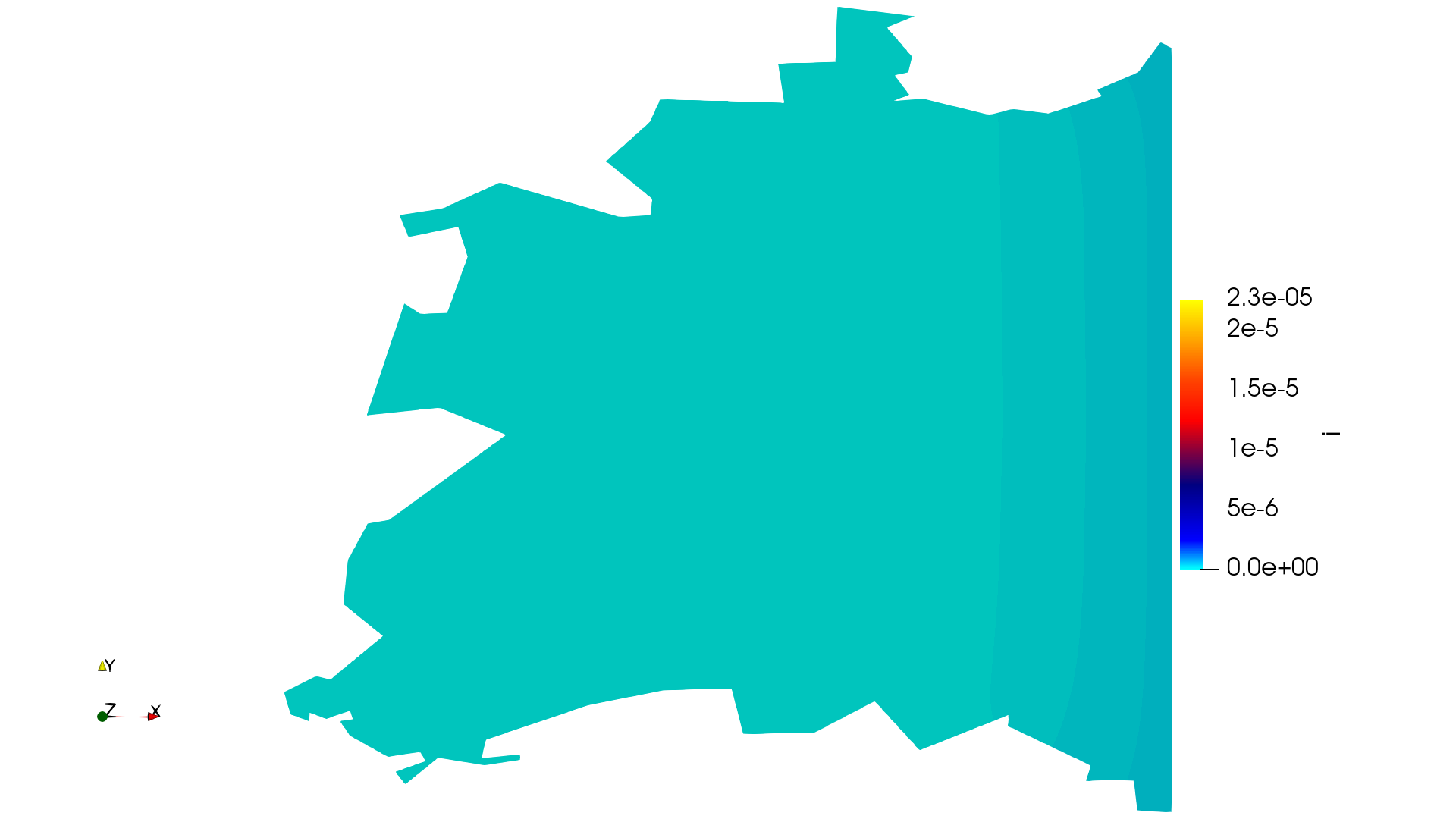} \\
         \includegraphics[width=\textwidth]{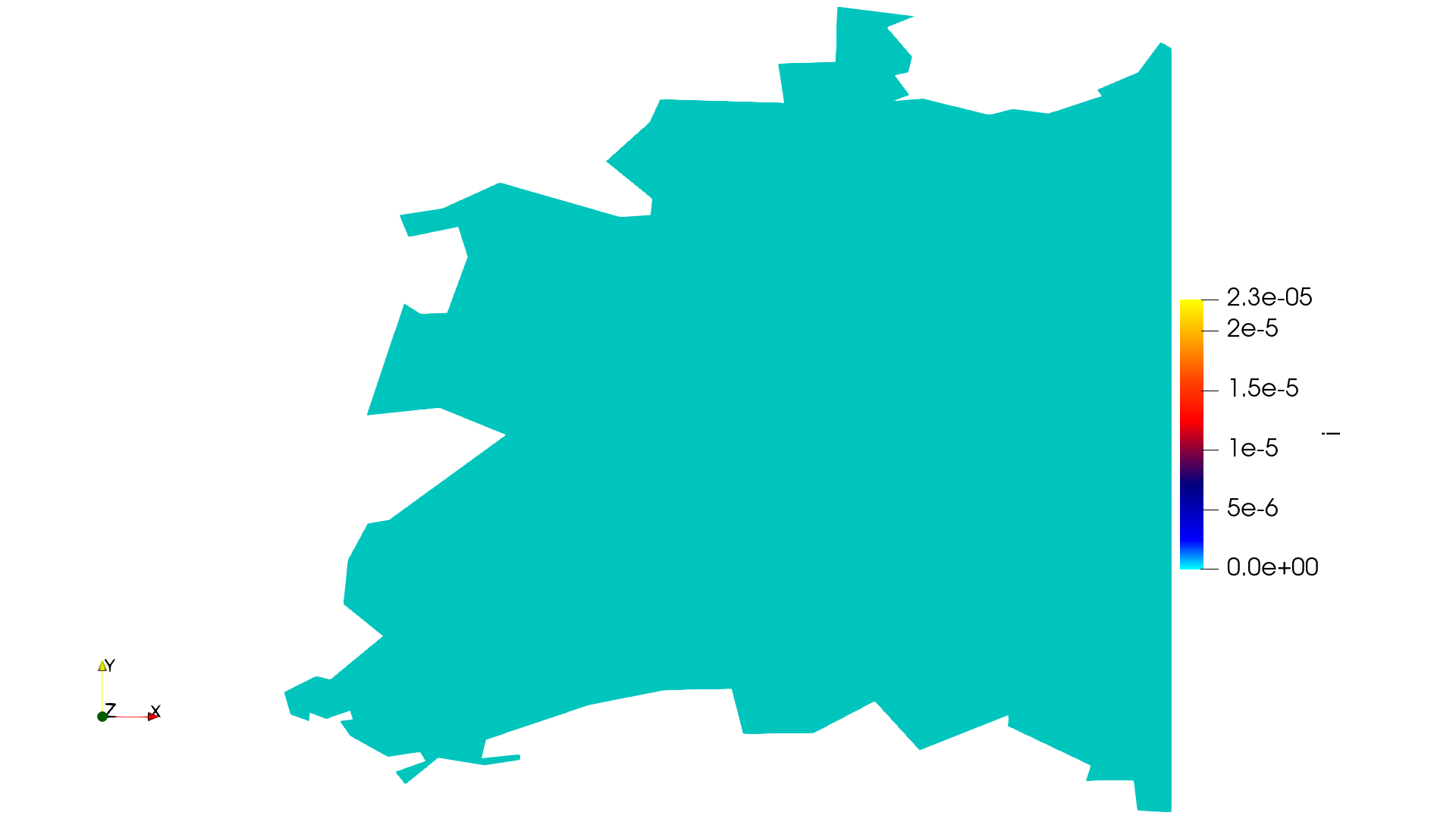}
    \end{minipage}
\caption{Fractional infectious density of hybrid model (\ref{Eq:hybrid_model}) for times $t \in \{1,8,18,60\}$. Due to Berlin’s significantly smaller area, the spread in Berlin appears faster than in Lombardy.} 
\label{fig:Berlin_fractional_infectious_density_hybrid}
\end{figure}

\begin{figure}[h]
\begin{minipage}{0.49\textwidth}
    \includegraphics[width=\textwidth]{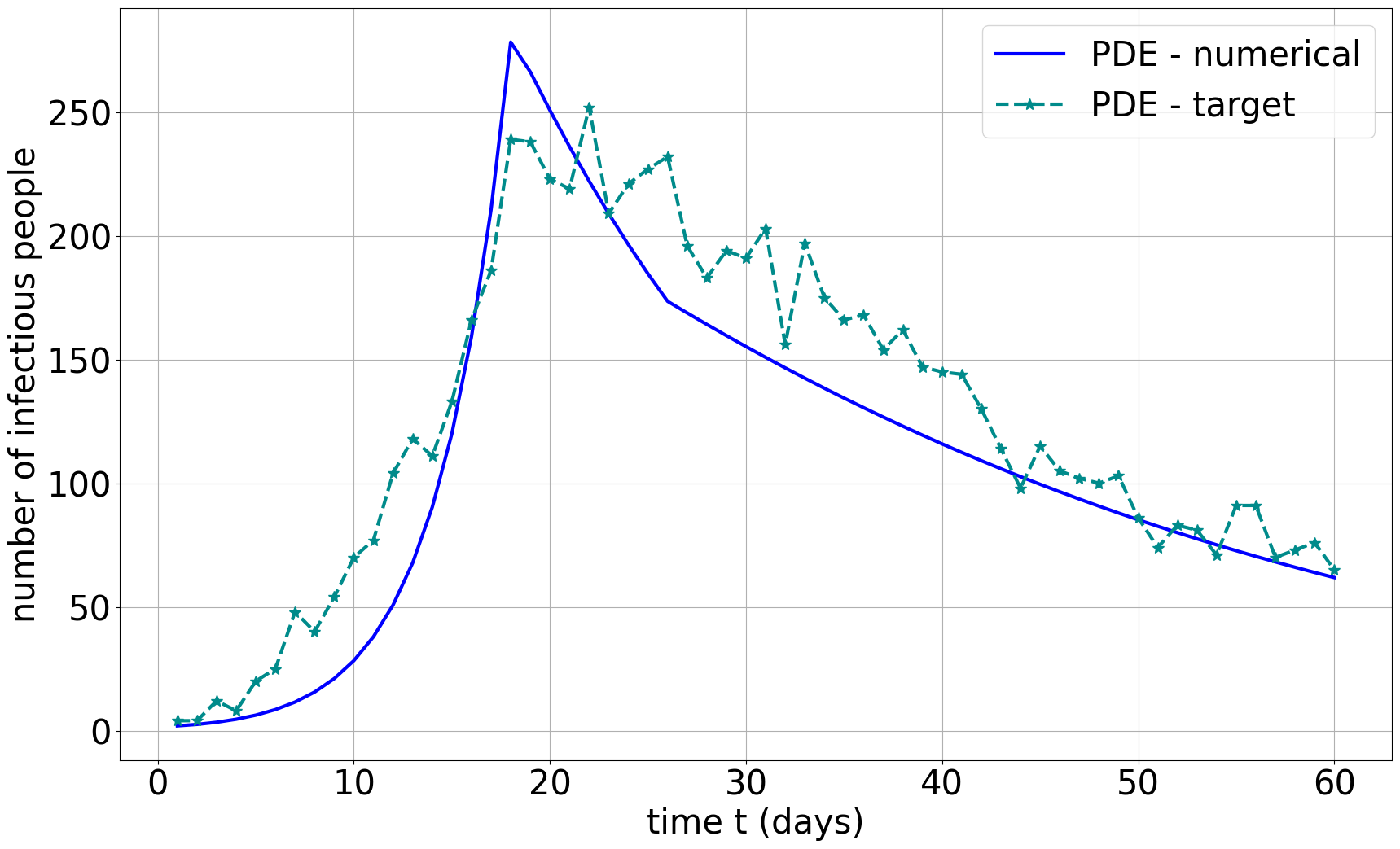}
\end{minipage}
\begin{minipage}{0.49\textwidth}
\includegraphics[width=\textwidth]{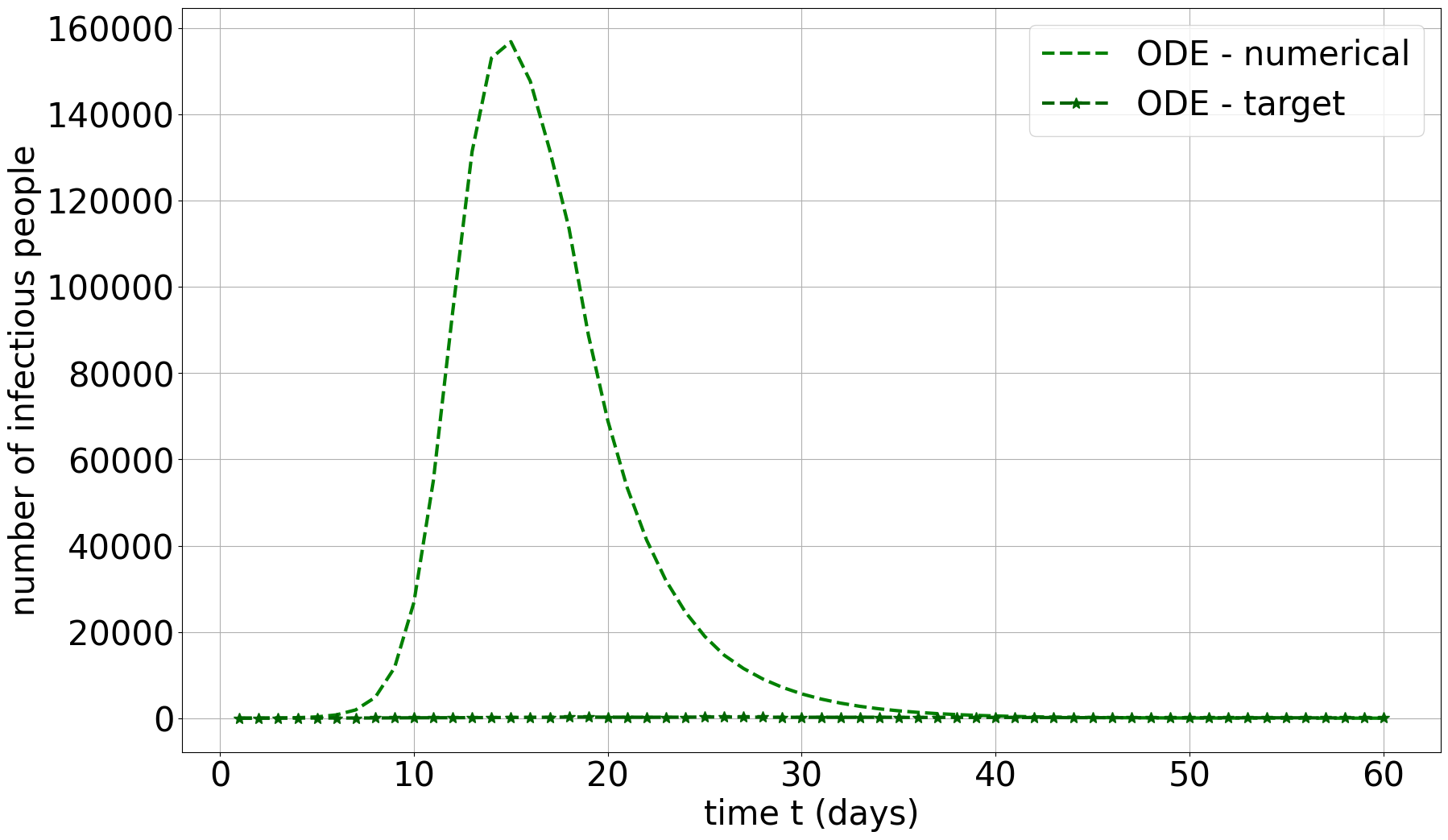}
\end{minipage}
    \caption{Number of infectious people of PDE contribution (left) and ODE contribution (right) of hybrid model (\ref{Eq:hybrid_model}) in Berlin (simulated and ABM data) with zero boundary conditions. Once the flow across the boundary is prohibited, it becomes apparent that the infectious numbers only explode within the ODE region.} 
\label{fig:Berlin_number_of_infectious_zero_BCs}
\end{figure}

\begin{figure}[h]
    \begin{minipage}{0.49\textwidth}
        \includegraphics[width=\textwidth]{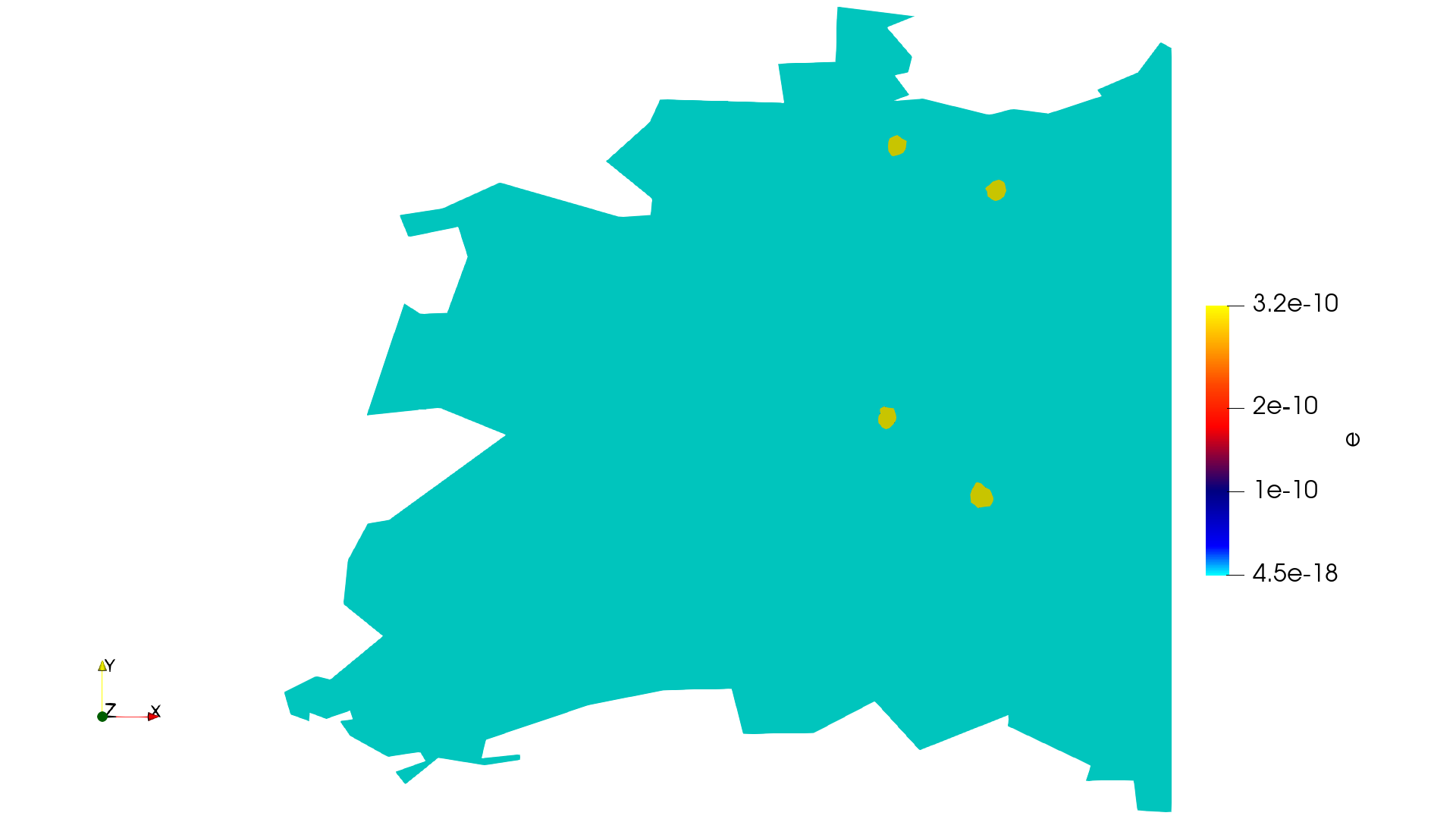} \\
        \includegraphics[width=\textwidth]{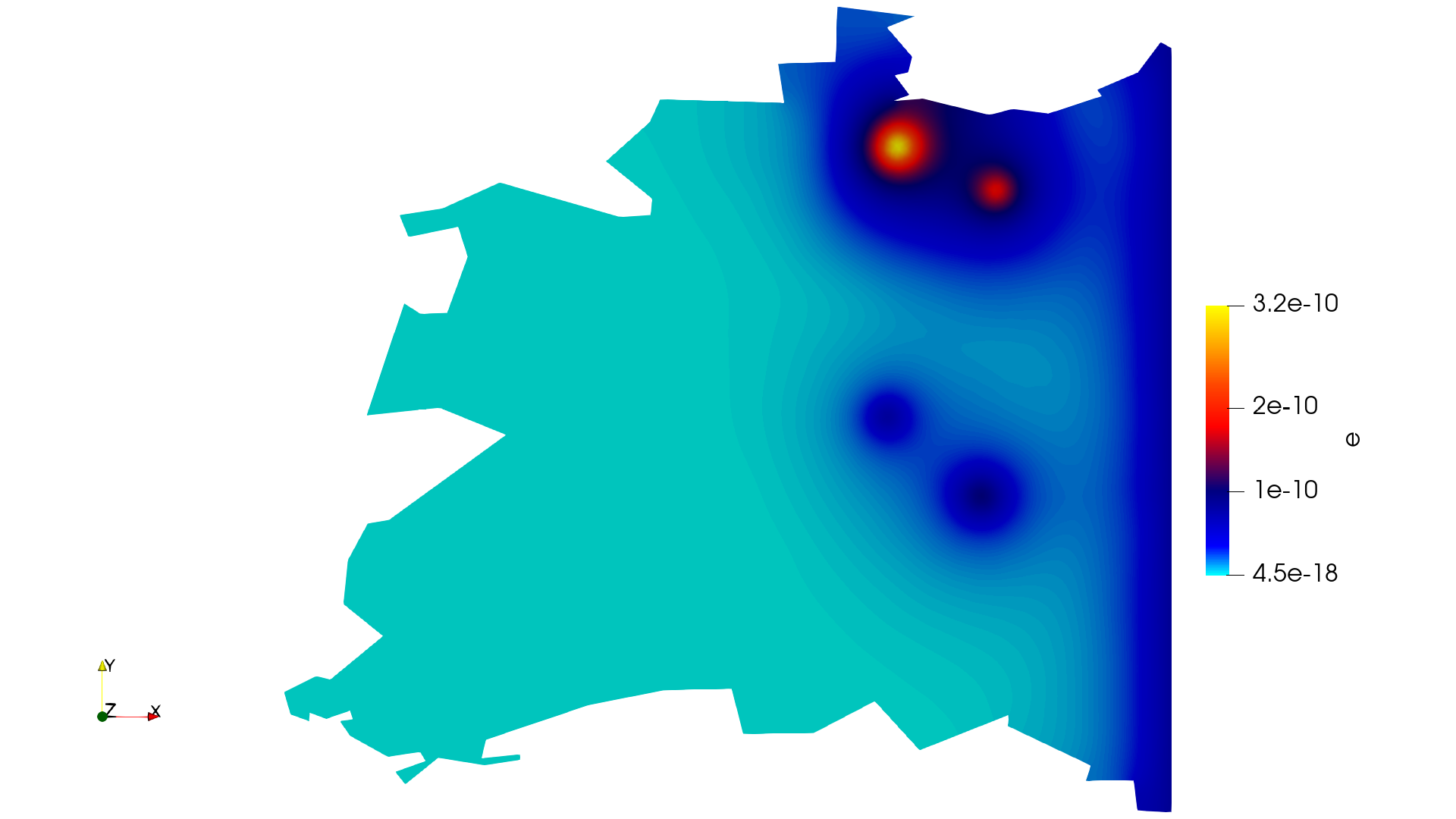} \\
        \includegraphics[width=\textwidth]{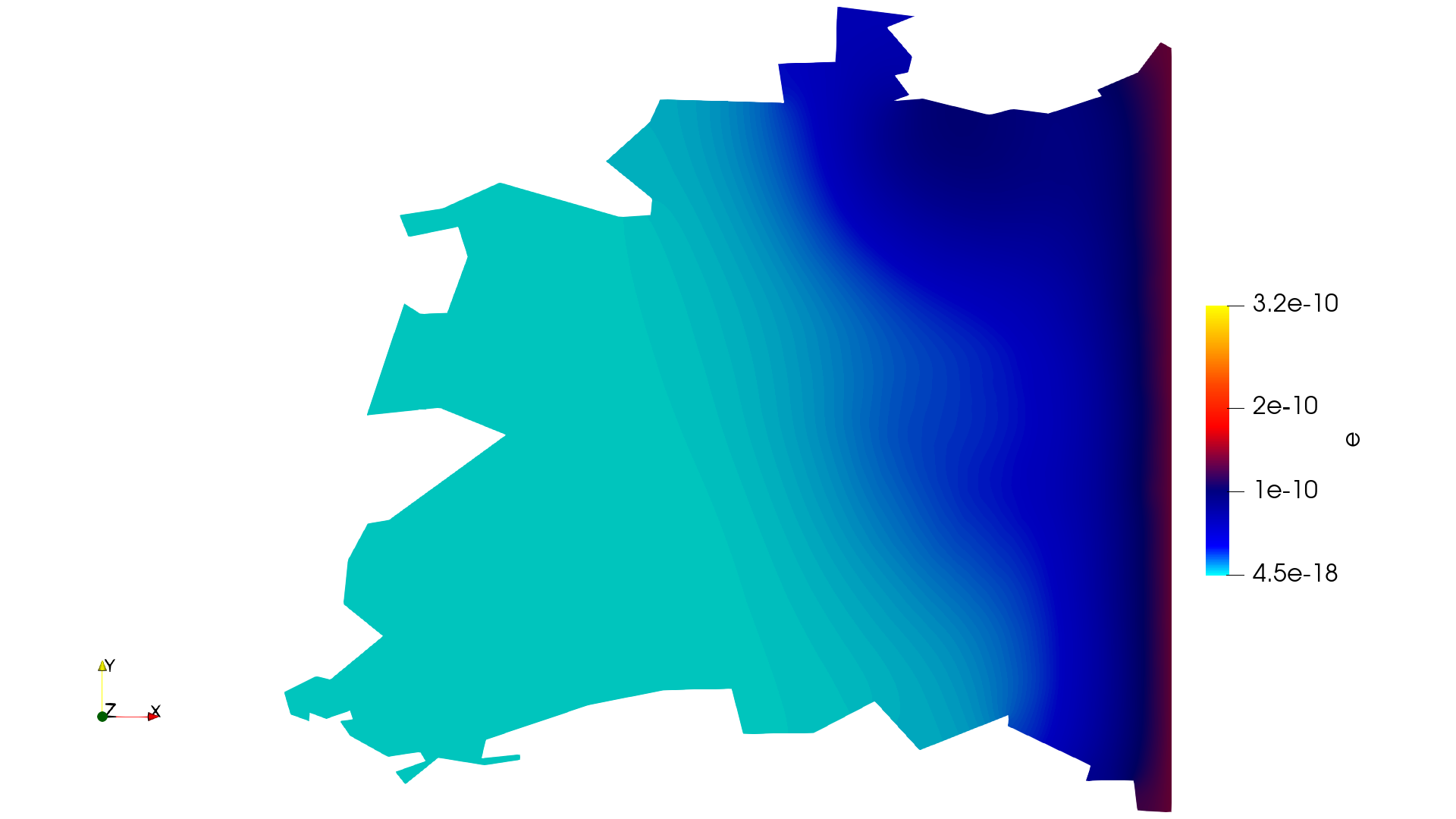}
    \end{minipage}
    \begin{minipage}{0.49\textwidth}
         \includegraphics[width=\textwidth]{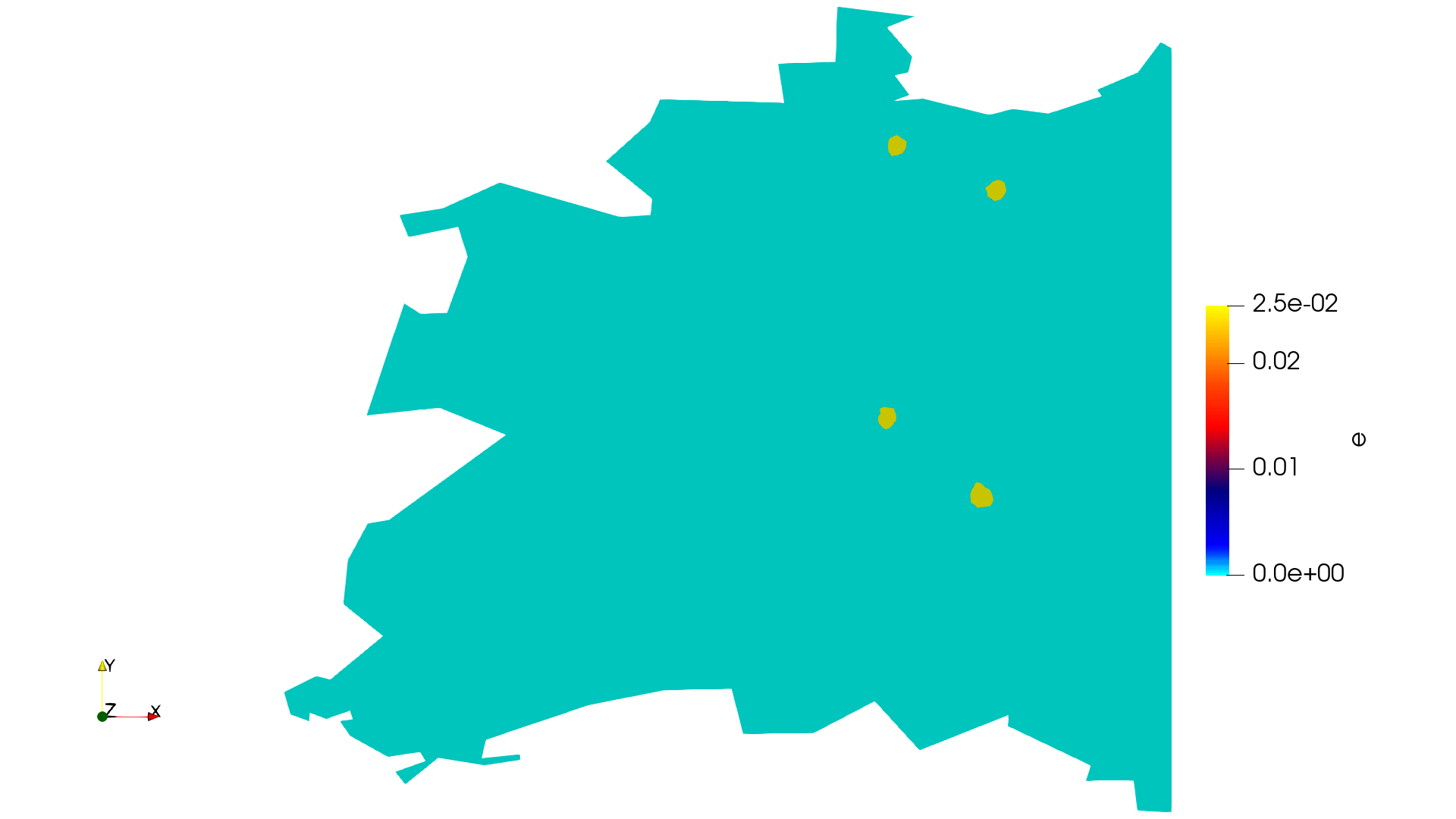} \\
        \includegraphics[width=\textwidth]{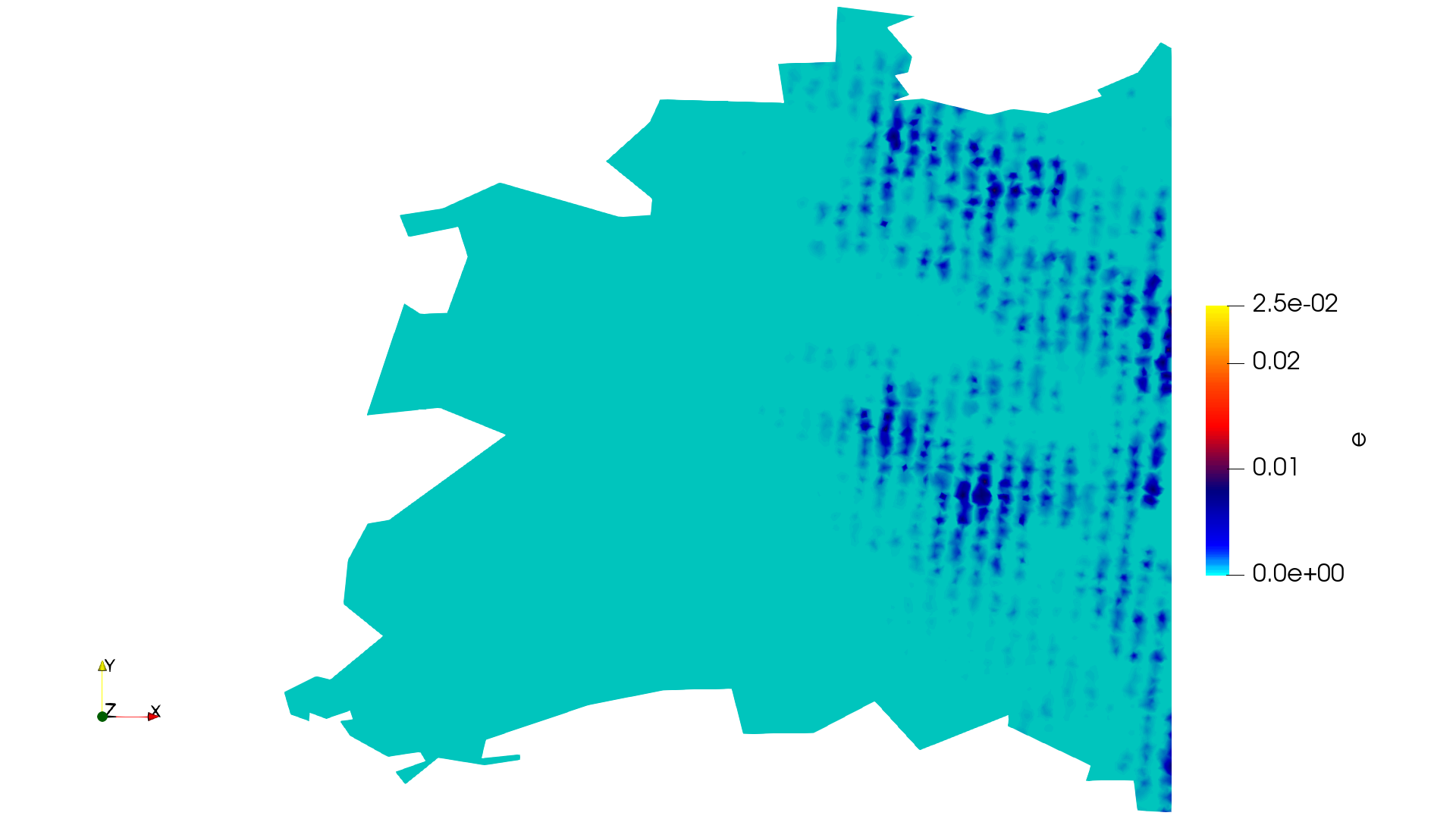} \\
        \includegraphics[width=\textwidth]{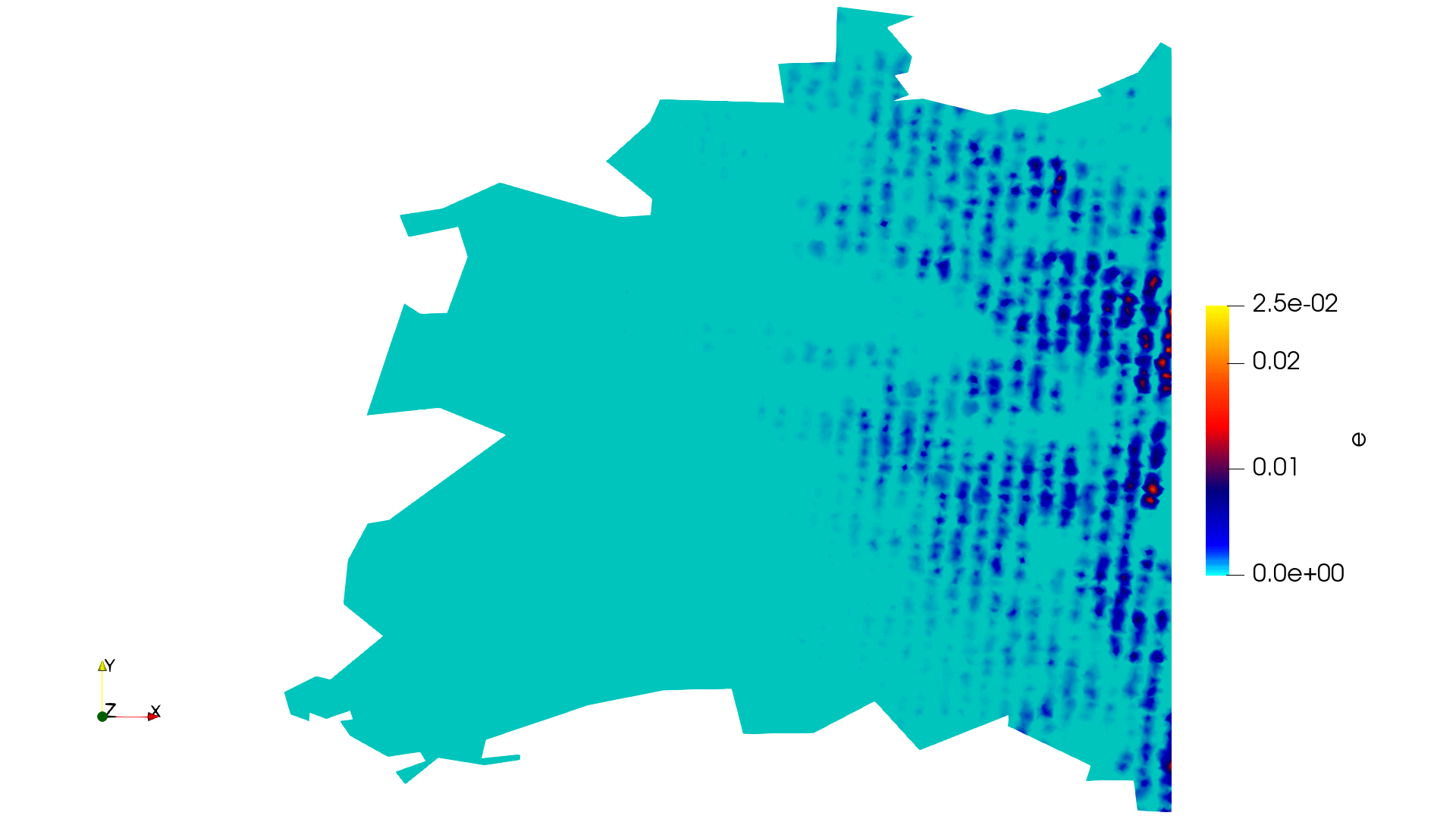}
    \end{minipage}
\caption{(Fractional) Exposed density of hybrid model (\ref{Eq:hybrid_model}) for times $t \in \{0,1,2\}$. The scale for the initial time step is incorrect, as the results in the first image exceed the range of the scale. The initial time step was included solely to illustrate the initial position of exposed individuals. By time step two, the infection flow originating from the boundary has already begun to have a greater impact than the exposed individuals.} 
\label{fig:Berlin_(fractional)_exposed_density_hybrid}
\end{figure}

\begin{figure}[h]
    \centering
    \includegraphics[width=0.6\textwidth]{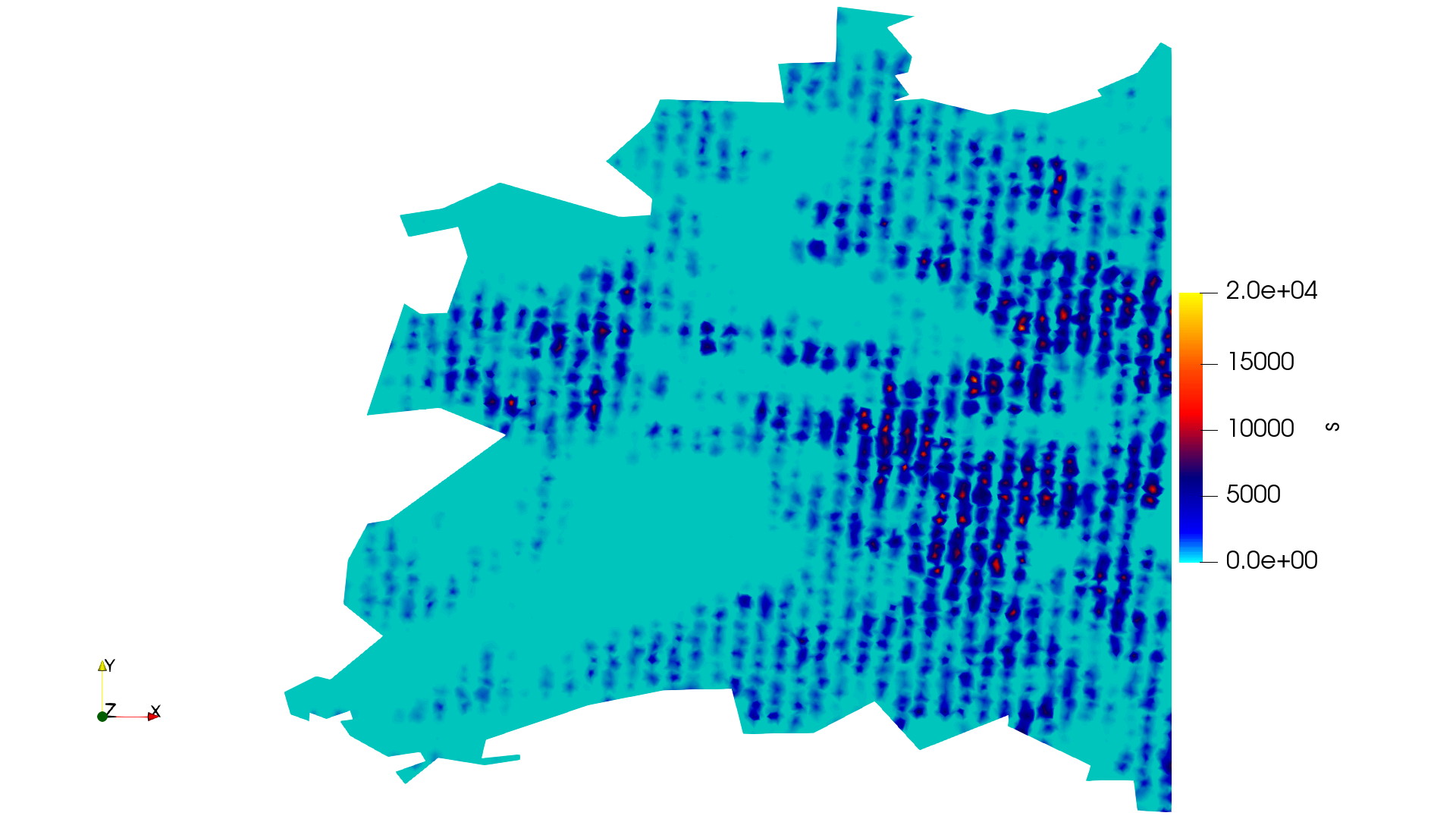} 
    \caption{Susceptible density of hybrid model (\ref{Eq:hybrid_model}) for time $t = 0$. The population density roughly corresponds to the susceptible density initially.} 
    \label{fig:Berlin_susceptible_density_hybrid}
\end{figure}

We adopted the diffusion coefficient for Berlin in a manner consistent with that used for Lombardy, although it's worth noting that different time intervals were employed, and Berlin was simulated for an overall longer duration. However, due to Berlin's significantly smaller area, the spread in Berlin appears faster in relation to the geographical size. This is clearly visible in Figs.  \ref{fig:Berlin_fractional_infectious_density_PDE}, \ref{fig:Berlin_fractional_infectious_density_hybrid}, where we have plotted the fractional infectious density. The spread of infections can be observed particularly well at the level of home locations in Fig. \ref{fig:Berlin_infectious_density_PDE}. Unsurprisingly, once again we observe that 
the fitting of the PDE model contributes to the closer alignment with the target data (see Fig. \ref{fig:Berlin_number_of_infectious}). 
We observe a quite similar issue here, which we have already seen in experiments conducted in a rectangular area (see Sec. \ref{sec:ExtremeCases}). The uniform flow of infection from the ODE region, combined with a high population density directly at the boundary between the differently modeled regions, dominates, even though the population density is also high in other areas of Berlin (compare to susceptible density in Fig. \ref{fig:Berlin_susceptible_density_hybrid}). 
This results in a sharp increase in infection numbers in the hybrid model. These two exposed agents in the north have a stronger effect on the spread than the other two exposed agents further south (see Fig. \ref{fig:Berlin_(fractional)_exposed_density_hybrid}). This is likely because we have no-flux Neumann boundary conditions in the north, preventing the infection from spreading out of Berlin. Note that the population density roughly corresponds to the susceptible density initially, and in Figure \ref{fig:Berlin_susceptible_density_hybrid}, it can be seen that the population density for the exposed agents in the north is lower than those further south. 
Directly at the boundary, the fractional infectious density is slightly higher than in the adjacent regions. This is again due to the Dirichlet boundary conditions we enforce adjacent to the ODE region. 
The infection flow slowly reaches the densely populated region in western Berlin, as we have a large area in the region with very low population density, and many areas between the home locations with very low population density, thus resulting in lower infection rates in these areas. Additionally, we only have an infection flow originating from the boundary between the PDE and ODE regions, along with a few initial COVID cases that are not located in western Berlin. Including infection flows from outside Berlin would be advantageous here.
Once the flow across the boundary is prohibited, by setting the boundary conditions between the two regions to zero, it becomes apparent that the infectious numbers only explode within the ODE region (see Fig. \ref{fig:Berlin_number_of_infectious_zero_BCs}). Originally, the boundary conditions were intended to maintain the fractional densities on both sides of the boundary as similar as possible, which is why the fractional infectious density in the PDE region has increased significantly. This, in turn, can cause an escalating effect in the ODE region, as we can observe.
The primary reason for the varying number of infectious cases is likely due to the non-uniform distribution of the population within the PDE domain, with gaps between the home locations. Achieving a uniform distribution of home locations is not feasible in this experiment, as we specifically aimed to test precise data and its impact on the infection dynamics of the hybrid model. What we did in this regard was to assign home locations to a grid point and distribute them more evenly within the adjacent triangles whose vertices touch the very same grid point. Additionally, one could consider a non-adaptive coarser grid, which would make the solution less accurate. Further experiments are required in this regard and will be explored in the future.
\noindent
\noindent
\begin{figure}[ht]
    \begin{minipage}{0.49\textwidth}
        \includegraphics[width=\textwidth]{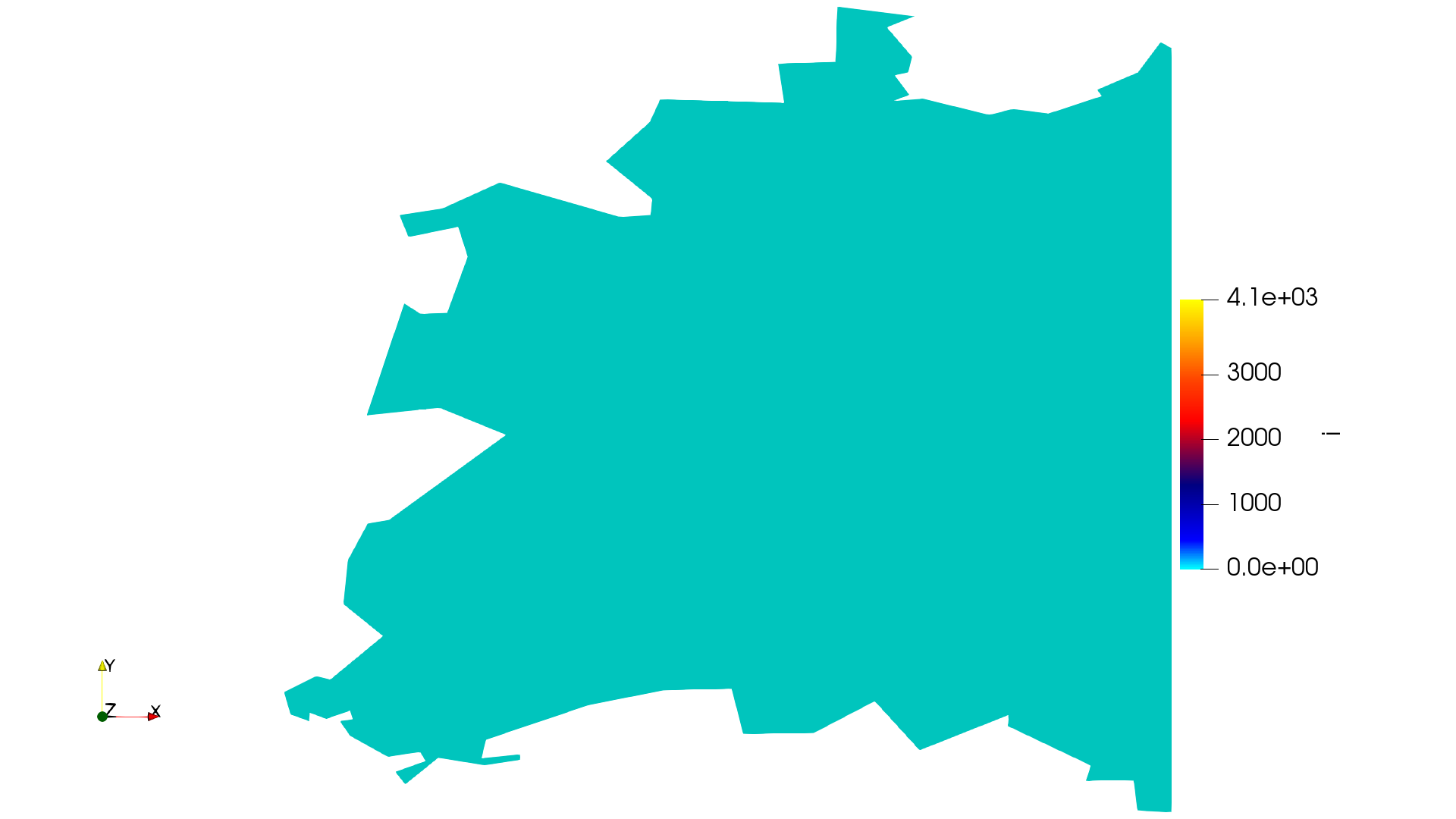} \\
        \includegraphics[width=\textwidth]{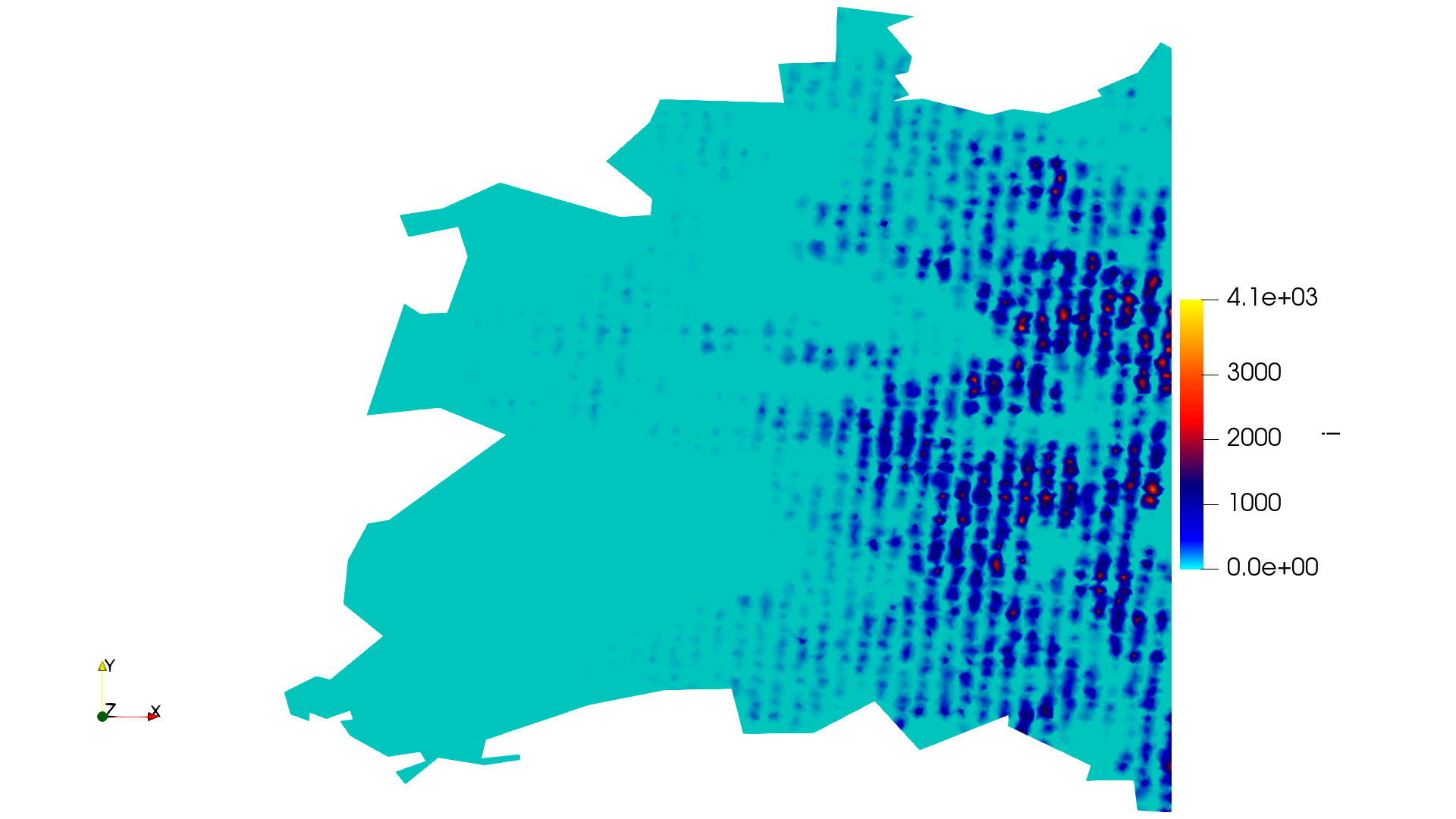}
    \end{minipage}
    \begin{minipage}{0.49\textwidth}
         \includegraphics[width=\textwidth]{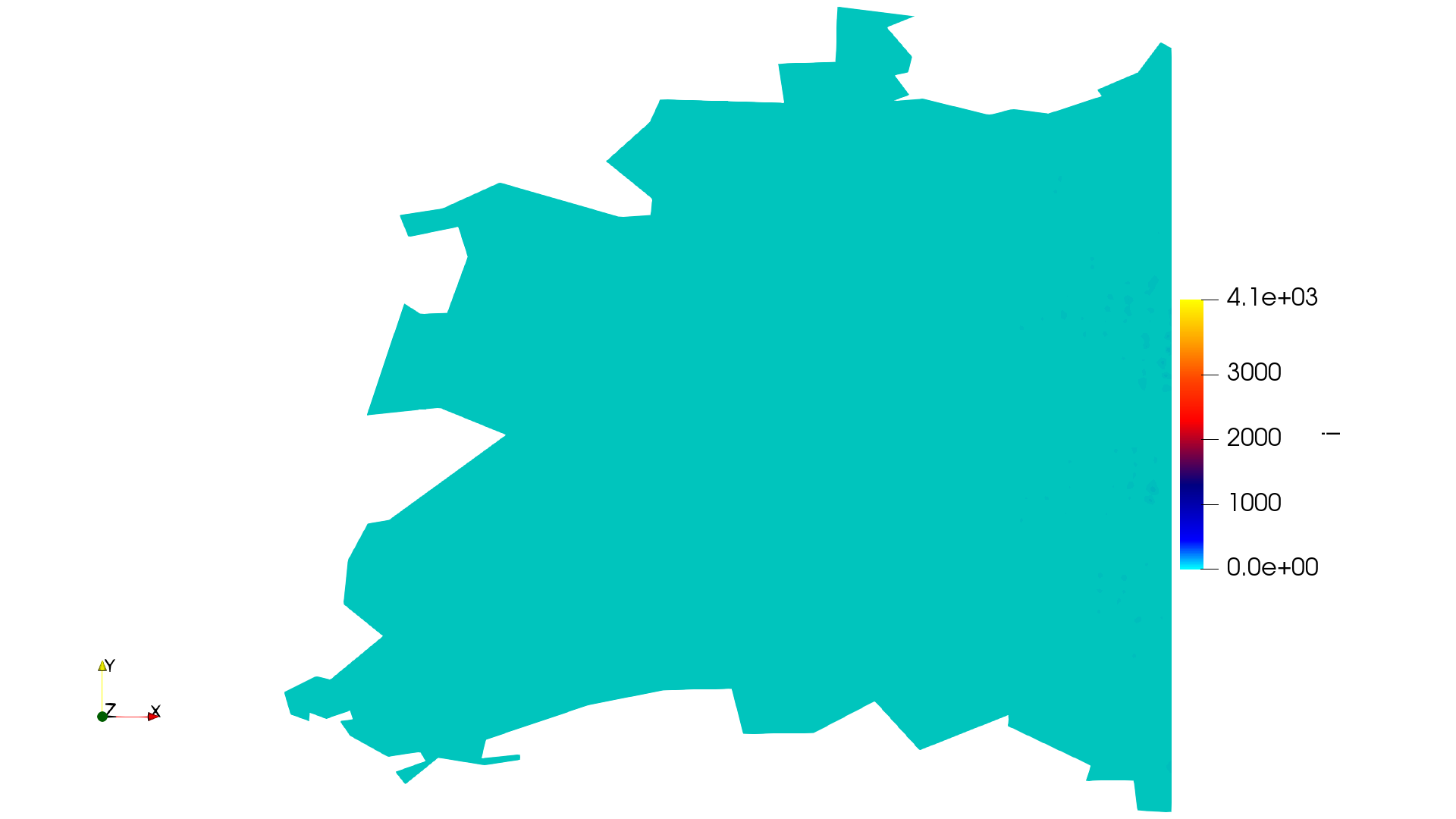} \\
         \includegraphics[width=\textwidth]{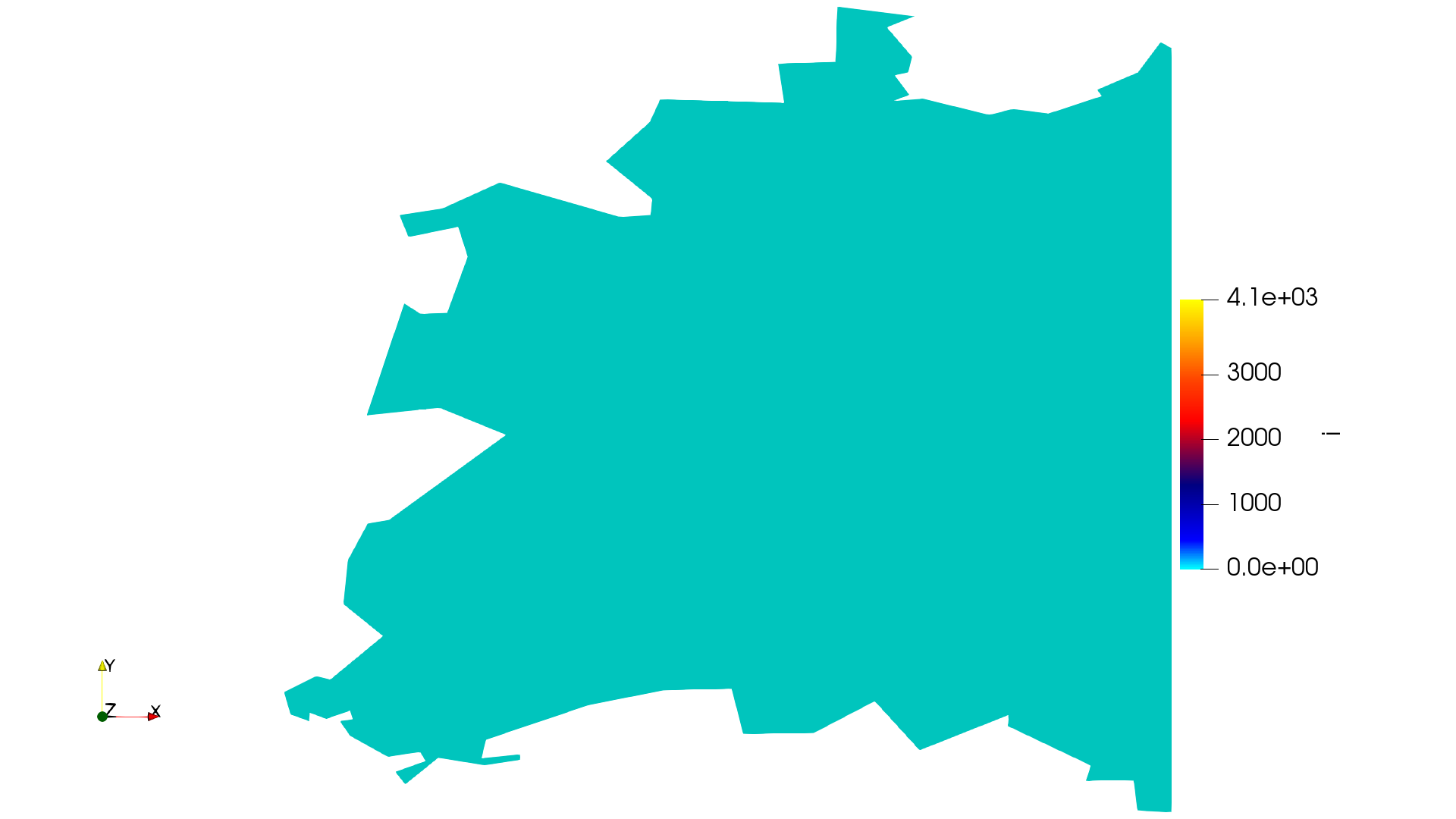}
    \end{minipage}
\caption{Infectious density of hybrid model (\ref{Eq:hybrid_model}) for times $t \in \{1,8,18,60\}$. Here, the spread can be observed particularly well at the level of home locations.} 
\label{fig:Berlin_infectious_density_hybrid}
\end{figure} 
\noindent
\noindent
\noindent
\noindent
\noindent
\newpage
\textcolor{white}{}
\newpage
\textcolor{white}{}
\newpage
\section{Conclusions}

In this paper, we introduced a novel hybrid mathematical modeling approach that couples Partial Differential Equations (PDEs) with Ordinary Differential Equations (ODEs) to simulate epidemiological processes effectively across geographic landscapes. Our model exemplifies the integration of the spatially detailed representation of disease dynamics provided by PDEs with the computational efficiency of ODEs, allowing for the timely assessment of public health interventions and potential impacts of infectious diseases across large populations. Through the implementation and evaluation of this hybrid PDE-ODE model in a simple setting using synthetic data and in two distinct real-world settings, specifically Lombardy and Berlin, we explored how interactions between ODE and PDE regions can significantly influence infection dynamics. The results from these case studies highlight the model’s adaptability and its ability to reflect complex epidemiological behaviors in various settings.

In addition to these considerations, incorporating the flow of infections from neighboring areas is crucial for the efficacy of our hybrid model. In the case of Milan, the initial number of infectious individuals was too low to account for the significant increase in infection numbers. Moreover, the accurate placement of boundary conditions is essential; improper boundaries can lead to misleading results and inadequate public health responses, as they may fail to reflect the true dynamics of infection spread. Thus, ensuring appropriate boundary placements not only enhances the model's accuracy but also strengthens its utility in real-world epidemiological applications.

Building on these insights, our findings suggest the necessity for strategic adjustments to minimize potential risks, particularly the steady flow of infections that has a heightened impact on neighboring areas with high population density. To manage the spread of infectious diseases more effectively, we recommend modifying the size and location of ODE regions or creating buffer zones within PDE domains. For example, it is advisable to avoid splitting well-connected areas like Berlin at their center. Instead, the ODE region can be shifted or resized to minimize its impact on densely populated areas. If direct changes to the ODE regions are not possible or too complicated, implementing transition zones within the PDE domains could serve as an effective alternative. These zones can be established without altering the proportionate sizes of the regions. For instance, in the Lombardy case-study, a transition zone was indirectly implemented by not initially placing individuals from Milan in this area. Instead, a few individuals from Lodi or Bergamo were distributed through the Gaussian distribution in the transition zone. 
Furthermore, on the rectangular domain, a more uniform population distribution led to better control of infection dynamics. However, the significant rise in infection numbers in the Berlin ODE region highlights the need for further experiments to refine this strategy.

Directly fitting the parameters of the hybrid model, rather than adapting those from the PDE model, would enhance its performance by allowing precise adjustments tailored to the combined dynamics of the ODE and PDE regions. This adjustment ensures that the infection numbers for both models align closely with the target data, providing a more accurate reflection of the epidemiological situation.

The exploration of infection dynamics beyond the confines of Lombardy and Berlin introduces additional complexities and opportunities for our modeling approach. By expanding our analysis to include other regions, the hybrid model proves very useful for simulating extensive geographical areas and understanding the broader implications of disease spread. Incorporating dynamics from external areas could further deepen our understanding and enhance the model's applicability.

Finally, considering the urban spread dynamics, introducing a location-dependent diffusion coefficient could control infection flow at the periphery, potentially speeding up the spread within cities and modifying or replacing the Allee term’s role. To further enhance our model, we plan to incorporate transportation dynamics, including trains and airplanes, using stochastic differential equations with discrete jumps. By expanding the scope of our hybrid model to include more detailed interactions and broader geographic applications, we aim to make ongoing contributions to the field of epidemiological forecasting and analysis, providing more precise and effective tools for managing public health crises.

\section*{Acknowledgments}
We acknowledge Kai Nagel for valuable discussions on the agent-based model. We also thank Jakob Benjamin Rehmann and Sydney Cornelia Paltra for providing essential data for the Berlin experiments and for assistance with their software.

The work on this paper was funded by the German Ministry of research and education (BMBF) (project ID: 01KX2022A) and by the Deutsche Forschungsgemeinschaft (DFG, German Research Foundation) under Germany’s Excellence Strategy via MATH+: The Berlin Mathematics Research Center (EXC-2046/1, project ID: 390685689). 

\newpage
\clearpage

\bibliography{references}{}
\bibliographystyle{plain}



\newpage
\appendix

{\Large Appendices}

\section{Levenberg-Marquardt Algorithm}
\label{appendix:Levenberg_Marquardt}

\begin{algorithm}
\caption{Levenberg-Marquardt algorithm} 
\label{alg:Levenberg_Marquardt}
\begin{algorithmic}[1]
\Procedure{LM}{$F(p^{(0)}), \text{DF}(p^{(0)}), p^{(0)}$} 
\State $l \gets 0$
\State $\lambda_2 \gets 1$
\State $\lambda_2^{\text{positive}} \gets 1$
\State $||F(p^{(-1)})||_2 \gets -1$
\Repeat
    \State choose $\lambda_1 \geq 0$ so that $\text{DF}(p^{(l)})^T \text{DF}(p^{(l)}) + \lambda_1 I$ is strictly diagonally dominant
    \State $\Delta p^{(l)} \gets \text{solve}\bigl( \text{DF}(p^{(l)})^T \text{DF}(p^{(l)}) + \lambda_1 I, - \text{DF}(p^{(l)})^T F(p^{(l)}) \bigr)$
    \\
    \State decreased\_resi $\gets \text{true}$
    \If{$||F(p^{(l-1)})||_2 \geq 0$}
        \If{$||F(p^{(l)})||_2 < ||F(p^{(l-1)})||_2$}
        \State $\lambda_2 \gets 1$ \Comment{reset damping coefficient }
        \ElsIf{$||F(p^{(l)})||_2$ sufficiently small after some runs} 
        \State $\lambda_2 \gets 1$  \Comment{reset because $\lambda_2$ too small}
        \Else 
        \State decreased\_resi $\gets \text{false}$
        \State $\lambda_2 \gets \frac{1}{2} \lambda_2$ \Comment{choose smaller damping coefficient}
        \EndIf 
    \EndIf \\
    \If{not decreased\_resi}
        \State $l \gets l-1$ \Comment{go back and repeat with new choice of parameters}
    \EndIf \\
    \State $\lambda_2^{\text{positive}} \gets \lambda_2$
    \State $p^{(l+1)} \gets p^{(l)} + \lambda_2^{\text{positive}} \Delta p^{(l)}$ 
    \While{$p^{(l+1)} \leq 0$} \Comment{in case new parameters are not positive}
    \State $\lambda_2^{\text{positive}} \gets \frac{1}{2} \lambda_2^{\text{positive}}$
    \State $p^{(l+1)} \gets p^{(l)} + \lambda_2^{\text{positive}} \Delta p^{(l)}$ 
    \EndWhile \\
    \If{not decreased\_resi} 
    \State $\lambda_2 \gets \lambda_2^{\text{positive}}$ 
    \EndIf
    \State $p^{(l)} \gets p^{(l+1)}$ \Comment{continue with new iteration}
    \State $l \gets l+1$
    %
    %
\Until{$||F(p^{(l)})||_2$ is sufficiently small or number of maximum iterations is reached}
\EndProcedure
\end{algorithmic}
\end{algorithm}

\end{document}